\documentclass[reqno,11pt]{amsart}
\usepackage[colorlinks=true, linkcolor=blue, citecolor=blue]{hyperref}
\usepackage{cite}   
\usepackage{amssymb}
\usepackage{amsmath, graphicx, rotating}
\usepackage{color}
\usepackage{soul}
\usepackage[dvipsnames]{xcolor}

\usepackage{ifthen}
\usepackage{xkeyval}
\usepackage{todonotes}
\usepackage{amsmath}
\usepackage{amsfonts}
\usepackage[T1]{fontenc}
\usepackage{lmodern}
\usepackage[english]{babel}

\usepackage{ upgreek }
\usepackage{stmaryrd}
\SetSymbolFont{stmry}{bold}{U}{stmry}{m}{n}
\usepackage{amsthm}
\usepackage{float}

\usepackage{ bbm }
\usepackage{ stmaryrd }
\usepackage{ mathrsfs }
\usepackage{ frcursive }
\usepackage{ comment }

\usepackage{pgf, tikz}
\usetikzlibrary{shapes}
\usepackage{varioref}
\usepackage{enumitem}
\usepackage{ulem}

\definecolor{rouge}{rgb}{0.7,0.00,0.00}
\definecolor{vert}{rgb}{0.00,0.5,0.00}
\definecolor{bleu}{rgb}{0.00,0.00,0.8}
\usepackage[margin=1in]{geometry}
\newtheorem{theorem}{Theorem}[section]
\newtheorem*{theorem*}{Theorem}
\newtheorem{lemma}[theorem]{Lemma}

\newtheorem{corollary}[theorem]{Corollary}
\newtheorem{proposition}[theorem]{Proposition}

\labelformat{hypothesis}{\textbf{M\kern-0.1mm#1}}

\newtheorem{condition}{Condition}
\newtheorem{hypo}{}

\labelformat{conditionA}{\textbf{A\kern-0.1mm#1}}

\labelformat{conditionB}{\textbf{B\kern-0.1mm#1}}

\labelformat{conditionC}{\textbf{C\kern-0.1mm#1}}

\theoremstyle{definition}

\newtheorem{remark}[theorem]{Remark}

\numberwithin{equation}{section}

\def\geq{\geqslant}
\def\leq{\leqslant}

\allowdisplaybreaks

\begin{document}
	
	\title[Strong Convergence Rate]
	{Strong Convergence Rate for Non-Autonomous Slow-Fast SDEs Driven by Brownian Motions and $\alpha$-stable processes}
	\author{XIFENG SU}
	\address{Su.X: School of Mathematical Sciences, Laboratory of Mathematics and Complex Systems (Ministry of Education)\\
		Beijing Normal University,
		No. 19, XinJieKouWai St., HaiDian District, Beijing 100875, P. R. China}
	\email{xfsu@bnu.edu.cn, billy3492@gmail.com}
	
	\author{TIANXIANG WANG}
	\address{Wang.T: School of Mathematical Sciences, Beijing Normal University,
		No. 19, XinJieKouWai St., HaiDian District, Beijing 100875, P. R. China}
	\email{202631130013@mail.bnu.edu.cn}
	\begin{abstract}
		In this paper, we establish a quantitative strong averaging principle for a class of time quasi-periodic  slow-fast stochastic differential equations driven by both Brownian motions and multiplicative \(\alpha\)-stable L\'evy processes. We show that the slow component converges strongly to an
		autonomous averaged equation with rate
		$
		\varepsilon^{\min\{1-1/\alpha_2,\,1/(\nu+2)\}}
		$, where
		$\alpha_2\in(1,2)$ is the  fast stability index and $\nu$ comes from the quasi-periodic structure. 
		Our strategy is as follows. 
		\begin{itemize}\item We first derive the intermediate step of the averaging principle to average out the effect of fast-variable oscillations on the system, although the averaged variable still depends on $\varepsilon$.   We obtain the convergence rate \(\varepsilon^{1-1/\alpha_2}\).
			\item Furthermore, under the uniformly quasi-periodic conditions, we show that the slow component converges strongly to the solution of some autonomous averaged
			equation, whose coefficients independent of $\varepsilon$. We thus obtain the convergence rate \(\varepsilon^{1/(\nu+2)}\).
			
		\end{itemize}
		
		The main contribution of this paper is to extend the additive L\'evy noise to multiplicative L\'evy noise and extend the fast equation of the system from the autonomous setting to the non-autonomous setting. In particular, we establish an explicit convergence rate in the quasi-periodic setting, extending the corresponding results for the periodic case.
	\end{abstract}
	\date{\today}
	\maketitle
	\tableofcontents
	
	\section{Introduction}
	
	Multiscale stochastic systems arise naturally in a wide range of
	applications, including molecular dynamics,
	chemical kinetics,
	climate dynamics,
	biology,
	engineering,
	and mathematical finance, see e.g. \cite{LegollLelievre2010,KangKurtz2013,MajdaTimofeyevVandenEijnden2001,CruduDebusscheRadulescu2009,RobertsSpanos1986,FouquePapanicolaouSircarSolna2003}. 
	A common feature of such systems is the
	coexistence of "slow" and "fast" components evolving on substantially different time scales.
	When the separation between the slow and fast scales becomes large, it is
	natural to approximate the slow component by an effective dynamics in which
	the rapidly fluctuating variables have been averaged out. 
	The stochastic averaging principle goes back to the classical work of
	Khasminskii \cite{Khasminskii1968}. Since then, averaging and homogenization
	for multiscale stochastic systems have been extensively investigated; see,
	for example,
	\cite{FreidlinWentzell2012,PavliotisStuart2008,GivonKupfermanStuart2004}
	and the references therein. 
	
	For stochastic differential equations driven by
	Brownian motions, convergence of the slow component can be established by
	Khasminskii's time-discretization argument, martingale methods, or the Poisson
	equation approach. {Quantitative} averaging results have also been obtained for
	various finite and infinite dimensional stochastic systems; see e.g.
	\cite{KhasminskiiYin2004,Liu2010,Brehier2012,Brehier2013,
		LiuRocknerSunXie2020,RocknerSunXie2021,RocknerXie2021,YangYin2026}.
	In particular, Sun, Xie and Xie \cite{SunXieXie2022} studied slow-fast
	stochastic differential equations driven by additive $\alpha$-stable
	processes and obtained the optimal strong convergence order
	$
	1-1/\alpha
	$\footnote{The same convergence rate will be obtained in our intermediate averaging principle in the present multiplicative setting.}.
	
	When the driving noises contain jumps, the situation becomes substantially more delicate. 
	L\'evy processes provide a natural framework for modeling
	discontinuous random perturbations and heavy-tailed fluctuations; see
	\cite{Applebaum2009,Sato1999,PeszatZabczyk2007}. Among them,
	$\alpha$-stable processes play an important role because of their scaling
	properties and heavy-tailed distributions. 
	Especially, averaging principles for stochastic systems driven by such jump processes have  attracted considerable attention, see e.g. \cite{Liu2012,BaoYinYuan2017,SunXieXie2022,SunXie2023,SunXiaXieZhou,ChengHaoZhang2025,SunXieXie2022}.

	It is well-known that when the $\alpha$-stable noise is \emph{multiplicative}, the analysis becomes more involved. In fact, 
	when the jump coefficient depends on the state variable, the corresponding
	nonlocal part of the generator becomes a stable-like operator with a
	state-dependent jump kernel. Quantitative averaging in this setting would require
	suitable exponential  contraction and regularity estimates for the frozen
	dynamics. More concretely, in order to get exponential  contraction, one can apply
	coupling methods for stochastic differential equations driven by
	L\'evy processes; see, for
	example, \cite{Wang2018,LiuWang2019,LiangMajkaWang2021}.
	On the other hand, regularity estimates could also be obtained from the viewpoint of nonlocal operators of PDEs, see e.g. 
	\cite{PriolaZabczyk2011,Zhang2013,MenozziPesceZhang2021,BDVV2022,BDVV2025, SuValdinociWeiZhang2022,SuValdinociWeiZhang2025,PardouxVeretennikov2001,PardouxVeretennikov2003,
		PardouxVeretennikov2005} and references therein.
	
	\emph{Time-dependent} fast dynamics pose another difficulty, that is, the non-autonomous system requiring an evolution family of measures  instead of a single invariant measure.
	We note that the non-autonomous stochastic averaging has been studied from several different
	points of view. Averaging principles for time-dependent stochastic systems were investigated in
	\cite{LiuRocknerSunXie2020,DeFeo2021, DaPratoRockner2006,CerraiLunardi2017, WangXiaXieYang2024, LiSunWangXie2026} and so on. 
	
	The preceding works reveal two different difficulties in quantitative
	averaging theory. 
	\begin{itemize}
		\item On the one hand, multiplicative stable noises lead to
		state-dependent nonlocal generators and require refined coupling and
		regularity estimates. 
		\item On the other hand, non-autonomous fast dynamics requires
		an evolution system of measures and a time-inhomogeneous corrector equation.
	\end{itemize}
	
	The purpose of the present paper is to investigate a slightly more general framework driven by both Brownian motion and multiplicative $\alpha$-stable noise
	and to overcome both difficulties above simultaneously.
	More precisely, we consider the following non-autonomous slow-fast stochastic
	system:
	\begin{equation}\label{eq:system-section2}
		\left\{
		\begin{aligned}
			dX_t^\varepsilon
			={}&
			b_1(X_t^\varepsilon,Y_t^\varepsilon)\,dt
			+\sigma_1(X_t^\varepsilon)\,dW_t^1
			+g_1(X_{t-}^\varepsilon)\,dL_t^1,
			\\
			dY_t^\varepsilon
			={}&
			\frac{1}{\varepsilon}
			b_2\left(\frac{t}{\varepsilon},
			X_t^\varepsilon,Y_t^\varepsilon\right)\,dt
			+
			\frac{1}{\sqrt{\varepsilon}}
			\sigma_2\left(\frac{t}{\varepsilon},
			X_t^\varepsilon,Y_t^\varepsilon\right)\,dW_t^2
			\\
			&+
			\frac{1}{\varepsilon^{1/\alpha_2}}
			g_2\left(\frac{t}{\varepsilon},
			X_{t-}^\varepsilon,Y_{t-}^\varepsilon\right)\,dL_t^2,
		\end{aligned}
		\right.
	\end{equation}
	where $t\geq0$, $\varepsilon>0$, and the initial conditions are
	$
	X_0^\varepsilon=x\in\mathbb{R}^{d_1},
	Y_0^\varepsilon=y\in\mathbb{R}^{d_2}.
	$
	The drift coefficients are functions
	\[
	\begin{aligned}
		b_1&:\mathbb{R}^{d_1}\times\mathbb{R}^{d_2}
		\longrightarrow\mathbb{R}^{d_1},\\
		b_2&:\mathbb{R}\times\mathbb{R}^{d_1}
		\times\mathbb{R}^{d_2}
		\longrightarrow\mathbb{R}^{d_2},
	\end{aligned}
	\]
	while the Brownian and jump coefficients\footnote{One can refer Section~\ref{sec: preliminary} for more detailed assumptions of coefficients.}
	are matrix-valued functions
	\[
	\begin{aligned}
		\sigma_1,g_1&:\mathbb{R}^{d_1}
		\longrightarrow\mathbb{R}^{d_1\times d_1},\\
		\sigma_2,g_2&:\mathbb{R}\times\mathbb{R}^{d_1}
		\times\mathbb{R}^{d_2}
		\longrightarrow\mathbb{R}^{d_2\times d_2}.
	\end{aligned}
	\]
	For $i=1,2$, $W^i$ is a $d_i$-dimensional standard Brownian motion,
	and $L^i$ is a $d_i$-dimensional isotropic $\alpha_i$-stable
	L\'evy process with $\alpha_i\in(1,2)$.
	The four driving processes $W^1,W^2,L^1,L^2$ are mutually independent.	
	Note that the coefficients of the fast component depend explicitly
	on the fast time variable $t/\varepsilon$, and both the Brownian and
	stable noise coefficients of the fast dynamics may depend on the
	state variables.
	
	Our strategy is to divide the averaging process into two steps by establishing both averaged equations \eqref{eq:first-averaged} and \eqref{secondaveraging} below, from which two averaging principles follow.
	
	For every fixed slow variable $x\in\mathbb{R}^{d_1}$ and $s\in\mathbb{R}$, we have the following frozen fast equation
	\begin{equation}\label{eq:intro_frozen}
		\left\{
		\begin{aligned}
			dY_t^{s,x,y}
			={}&
			b_2(t,x,Y_t^{s,x,y})\,dt
			+\sigma_2(t,x,Y_t^{s,x,y})\,dW_t^2
			\\
			&+
			g_2(t,x,Y_{t-}^{s,x,y})\,dL_t^2,
			\qquad t\geq s,
			\\
			Y_s^{s,x,y}
			={}&y.
		\end{aligned}
		\right.
	\end{equation}
	We remark that  the averaging structure of \eqref{eq:system-section2} is
	described by an evolution family of probability measures
	$
	\{\mu_t^x\}_{t\in\mathbb{R}}
	$ of $Y_t^{s,x,y}$ as $s\rightarrow -\infty$ since \eqref{eq:intro_frozen} is non-autonomous.

	Thus, the first averaged equation is then obtained by
	\begin{equation}\label{eq:first-averaged}
		\left\{
		\begin{aligned}
			d\bar X_t^\varepsilon
			={}&
			\bar b_1
			\left(
			\frac{t}{\varepsilon},
			\bar X_t^\varepsilon
			\right)\,dt
			+
			\sigma_1(\bar X_t^\varepsilon)\,dW_t^1
			+
			g_1(\bar X_{t-}^\varepsilon)\,dL_t^1,
			\\
			\bar X_0^\varepsilon
			={}&x
		\end{aligned}
		\right.
	\end{equation}where the time-dependent averaged drift $\bar b_1$ is defined by the probability measures $\{\mu_t^x\}_{t\in\mathbb{R}}$, i.e.,
	\[
	\bar b_1(t,x)
	:=
	\int_{\mathbb{R}^{d_2}}
	b_1(x,y)\,\mu_t^x(dy).
	\]

	Our first main Theorem~\ref{thm:strong-averaging-section2} is the strong averaging principle in the general case 
	and give the optimal convergence rate of $1-1/\alpha_2$, that is
	\begin{equation}
		\mathbb{E}
		\left[
		\sup_{0\leq t\leq T}
		|X_t^\varepsilon-\bar X_t^\varepsilon|^m
		\right]
		\leq
		C_{T,m,x,y}
		\bigl(
		1+|x|^m+|y|^m
		\bigr)
		\varepsilon^{m(1-1/\alpha_2)},
		\label{eq:strong-rate}
	\end{equation}
	where $ 1\leq m<\min\{\alpha_1,\alpha_2\}$, $T> 0$ and  $C_{T,m,x,y}>0$ is a constant independent of
	$\varepsilon\in(0,1]$.
	\smallskip
	
	The averaged process obtained in the first step still depends on
	$\varepsilon$. Indeed, this step averages out the oscillations of the
	fast state variable, while the temporal oscillations remain. We therefore
	perform a second averaging step to eliminate these temporal oscillations.
	Assuming that the coefficients of the fast equation are uniformly
	\emph{quasi-periodic} in time\footnote{Actually, time periodic case is an specific example in our study.}
	and the time-dependent averaged drift admits the representation
	\[
	\bar b_1(t,x)=B(\omega t,x),
	\]
	where $\omega\in\mathbb{R}^q$ satisfies the Diophantine
	condition with exponent $\nu$. We then define the averaged coefficient
	\[
	\widehat b_1(x)
	:=
	\int_{\mathbb{T}^q}
	B(\theta,x)\,d\theta,
	\]
	where $d\theta$ denotes the normalized Haar measure on
	$\mathbb{T}^q$. 
	
	The second averaged equation is then given by
	\begin{equation}\label{secondaveraging}
		\left\{
		\begin{aligned}
			d\widehat X_t
			={}&
			\widehat b_1(\widehat X_t)\,dt
			+
			\sigma_1(\widehat X_t)\,dW_t^1
			+
			g_1(\widehat X_{t-})\,dL_t^1,
			\\
			\widehat X_0
			={}&x,
		\end{aligned}
		\right.
	\end{equation}
	which is an autonomous system. 
	
	\smallskip
	
	Our second main Theorem~\ref{thm:quasiperiodic-averaging-section2} is the strong averaging principle in the quasi-periodic case and give the convergence rate of $\min\{1-1/\alpha_2,\,1/(\nu+2)\}$, that is
	\begin{equation}
		\mathbb E\left[\sup_{0\leqslant t\leqslant T}|X_t^\varepsilon-\widehat X_t|^m\right]
		\leqslant C_{T,m,x,y}\varepsilon^{m \min\{1-1/\alpha_2,\,1/(\nu+2)\}},
		\label{qpnew:main}
	\end{equation}
	where $ 1\leq m<\min\{\alpha_1,\alpha_2\}$, $T> 0$ and  $C_{T,m,x,y}>0$ is a constant independent of
	$\varepsilon\in(0,1]$.

	\smallskip
	
	We point out that our approach to get the first averaging principle employs  both coupling methods and regularity estimates, which is inspired from \cite{Wang2018,WangXiaXieYang2024}. 
	
	In the presence of multiplicative noise, the two particles may have different jump amplitudes. We construct a coupling based on reflected jumps that preserves the law of each marginal process and yields an exponential contraction estimate for the expected distance between the coupled particles. This establishes exponential contraction
	of the frozen dynamics with respect to its initial value. 
	
	We further 
	consider the time-inhomogeneous Poisson equation
	\[
	\partial_s\Phi(s,x,y)
	+
	\mathcal{L}_2(s,x)\Phi(s,x,y)
	+
	b_1(x,y)-\bar b_1(s,x)
	=
	0,
	\]
	where $\mathcal{L}_2(s,x)$ denotes the infinitesimal generator of the frozen
	fast process. Combining the regularity estimates for the Poisson equation
	with mollification arguments and moment estimates adapted to the stable
	noises, we prove that the slow component is strongly approximated by the
	equation with the time-dependent averaged drift. The resulting convergence
	order is
	$
	1-1/\alpha_2.
	$

	In the second averaging process, our novelty is to focus on systems with quasi-periodic time dependence. Existing studies on convergence rates for averaging principles in non-autonomous stochastic differential equations have primarily concerned with the periodic case. To the best of our knowledge, this work is the first to establish an explicit convergence rate in the quasi-periodic setting.
	
	\smallskip
	
	The paper is organized as follows. In Section~2, we
	introduce the notation and assumptions, describe the frozen fast dynamics and
	the associated evolution system of measures, and state the main strong
	averaging results. Section~3 is devoted to well-posedness and moment estimates
	for the slow-fast system and the frozen dynamics. In Section~4, we establish
	exponential contraction and regularity estimates for the frozen equation and
	study the associated non-autonomous Poisson equation. Section~5 contains the
	proof of the strong averaging principle and its convergence rate in general. In
	Section~6, we study further the averaging of the time-dependent averaged equation in the time quasi-periodic case.
	
	\section{Notation, Assumptions and Main Results}\label{sec: preliminary}
	
	In this section, we introduce the notations and assumptions used throughout
	the paper, describe the frozen fast dynamics and the associated evolution
	system of measures, and state the main strong averaging results. Besides the
	first averaging step with respect to the fast state variable, we also summarize
	the further averaging results for quasi-periodic temporal structures.
	
	\subsection{Notation}
	
	Throughout the paper, $C>0$ denotes a generic constant whose value may
	change from line to line. If the dependence of a constant on some parameters
	is relevant, we indicate it by subscripts.
	
	For $d\geq1$, let $\mathbb{R}^{d}$ be the $d$-dimensional Euclidean space
	equipped with the Euclidean norm $|\cdot|$ and inner product
	$\langle\cdot,\cdot\rangle$. For a matrix $A$, we denote by
	$\|A\|$ its operator norm and by $A^{*}$ its transpose. The unit sphere in
	$\mathbb{R}^{d}$ is denoted by
	\[
	\mathbb{S}^{d-1}
	:=
	\{z\in\mathbb{R}^{d}:|z|=1\}.
	\]
	
	For a bounded function $f$ defined on $\mathbb{R}^{d}$, we write
	\[
	\|f\|_{\infty}
	:=
	\sup_{x\in\mathbb{R}^{d}}|f(x)|.
	\]
	For $\gamma\in(0,1)$, the H\"older seminorm is defined by
	\[
	[f]_{C^\gamma}
	:=
	\sup_{x\neq y}
	\frac{|f(x)-f(y)|}{|x-y|^\gamma},
	\]
	and the corresponding H\"older norm is
	\[
	\|f\|_{C^\gamma}
	:=
	\|f\|_{\infty}+[f]_{C^\gamma}.
	\]
	The spaces $C_b^k$ and $C_b^{k+\gamma}$ are understood in the usual sense.
	
	Let $\{W_t^1\}_{t\geq0}$ and $\{W_t^2\}_{t\geq0}$ be independent
	$d_1$- and $d_2$-dimensional Brownian motions, respectively. Let
	$\{L_t^1\}_{t\geq0}$ and $\{L_t^2\}_{t\geq0}$ be independent isotropic
	$\alpha_1$-stable and $\alpha_2$-stable L\'evy processes in
	$\mathbb{R}^{d_1}$ and $\mathbb{R}^{d_2}$, respectively, where
	\[
	1<\alpha_1,\alpha_2<2.
	\]
	We assume that $W^1,W^2,L^1$ and $L^2$ are mutually independent.
	
	For $k=1,2$, the L\'evy measure of $L^k$ is given by
	\[
	\nu_k(dz)
	=
	c_{\alpha_k,d_k}
	\frac{dz}{|z|^{d_k+\alpha_k}},
	\qquad
	z\in\mathbb{R}^{d_k}\setminus\{0\},
	\]
	where $c_{\alpha_k,d_k}>0$ is a normalization constant. We denote by
	$N^k(dt,dz)$ the associated Poisson random measure and by
	\[
	\widetilde N^k(dt,dz)
	=
	N^k(dt,dz)-dt\,\nu_k(dz)
	\]
	its compensated Poisson random measure.
	
	For the fast stable process, it will also be convenient to use polar
	coordinates. Writing
	\[
	z=r\hat z,
	\qquad
	r=|z|>0,
	\qquad
	\hat z=\frac{z}{|z|}
	\in\mathbb{S}^{d_2-1},
	\]
	we may write
	\[
	\nu_2(dz)
	=
	\bar\nu_2(d\hat z)\frac{dr}{r^{1+\alpha_2}},
	\]
	where $\bar\nu_2$ is a symmetric finite Borel measure on
	$\mathbb{S}^{d_2-1}$.
	
	\subsection{Assumptions}
	
	We impose the following assumptions on the coefficients of
	\eqref{eq:system-section2}.
	
	\begin{hypo}[Uniform non-degeneracy]
		\label{h1}
		There exist constants $c_l,c_u>0$ such that, for every
		$t\in\mathbb{R}$, $x\in\mathbb{R}^{d_1}$,
		$y\in\mathbb{R}^{d_2}$,
		$\hat\omega_1\in\mathbb{S}^{d_1-1}$ and
		$\hat\omega_2\in\mathbb{S}^{d_2-1}$,
		\[
		c_l
		\leq
		|g_1(x)\hat\omega_1|
		\leq
		c_u,\qquad c_l
		\leq
		|\sigma_1(x)\hat\omega_1|
		\leq
		c_u
		\]
		and
		\[
		c_l
		\leq
		|g_2(t,x,y)\hat\omega_2|
		\leq
		c_u,
		\qquad
		c_l
		\leq
		|\sigma_2(t,x,y)\hat\omega_2|
		\leq
		c_u.
		\]
	\end{hypo}
	
	\begin{hypo}[Dissipativity of the fast dynamics]
		\label{h2}
		There exists a constant $\alpha>0$ such that, for every
		$t\in\mathbb{R}$, $x\in\mathbb{R}^{d_1}$ and
		$y_1,y_2\in\mathbb{R}^{d_2}$,
		\begin{align}
			&\left\langle
			b_2(t,x,y_1)-b_2(t,x,y_2),
			y_1-y_2
			\right\rangle
			\nonumber\\
			&\qquad
			+
			\|\sigma_2(t,x,y_1)-\sigma_2(t,x,y_2)\|^2
			\leq
			-\alpha|y_1-y_2|^2.
			\label{deepa1}
		\end{align}
		We further assume that
		\[
		\sup_{t\in\mathbb{R},\,x\in\mathbb{R}^{d_1}}
		|b_2(t,x,0)|
		<\infty.
		\]
		Consequently, there exists a constant $C>0$ such that
		\[
		\langle b_2(t,x,y),y\rangle
		\leq
		-\alpha|y|^2+C|y|.
		\]
	\end{hypo}
	
	\begin{hypo}[Global Lipschitz continuity]
		\label{h3}
		There exists a constant $C>0$ such that, for every
		$t_1,t_2\in\mathbb{R}$,
		$x_1,x_2\in\mathbb{R}^{d_1}$ and
		$y_1,y_2\in\mathbb{R}^{d_2}$,
		\[
		|b_1(x_1,y_1)-b_1(x_2,y_2)|
		\leq
		C\bigl(|x_1-x_2|+|y_1-y_2|\bigr),
		\]
		\[
		|b_2(t_1,x_1,y_1)-b_2(t_2,x_2,y_2)|
		\leq
		C\bigl(|t_1-t_2|+|x_1-x_2|+|y_1-y_2|\bigr),
		\]
		\[
		\|\sigma_1(x_1)-\sigma_1(x_2)\|
		\leq
		C|x_1-x_2|,
		\]
		\[
		\|g_1(x_1)-g_1(x_2)\|
		\leq
		C|x_1-x_2|,
		\]
		and
		\[
		\|\sigma_2(t,x_1,y_1)-\sigma_2(t,x_2,y_2)\|
		\leq
		C\bigl(|t_1-t_2|+|x_1-x_2|+|y_1-y_2|\bigr),
		\]
		\[
		\|g_2(t,x_1,y_1)-g_2(t,x_2,y_2)\|
		\leq
		C\bigl(|t_1-t_2|+|x_1-x_2|+|y_1-y_2|\bigr).
		\]
	\end{hypo}
	
	\begin{hypo}[Regularity of the fast coefficients]
		\label{h4}
		For every
		$(t,x)\in\mathbb{R}\times\mathbb{R}^{d_1}$,
		\[
		b_2(t,x,\cdot)
		\in
		C_b^2(\mathbb{R}^{d_2};\mathbb{R}^{d_2}),
		\]
		\[
		\sigma_2(t,x,\cdot)
		\in
		C_b^2
		(\mathbb{R}^{d_2};\mathbb{R}^{d_2\times d_2}),
		\]
		\[
		g_2(t,x,\cdot)
		\in
		C_b^1
		(\mathbb{R}^{d_2};\mathbb{R}^{d_2\times d_2}).
		\]
		and
		\begin{align*}
			\sup_{t,x}\big(\|b_2(t,x,\cdot)\|_{C_b^2} + \|\sigma_2(t,x,\cdot)\|_{C_b^2} + \|g_2(t,x,\cdot)\|_{C_b^1}\big) < \infty.
		\end{align*}
	\end{hypo}
	\begin{remark}
		Assumptions \ref{h1}--\ref{h4} are mainly imposed on the frozen fast
		dynamics. Assumption \ref{h2} provides the uniform dissipative mechanism
		responsible for the long-time contraction, while Assumption \ref{h1}
		provides the non-degeneracy needed for the regularity estimates.
	\end{remark}
	\begin{remark}
		In \ref{h2}, we impose a constant dissipativity coefficient
		$\alpha>0$ in order to keep the main argument transparent and to emphasize
		the effects of the multiplicative $\alpha$-stable noise and the quantitative
		averaging rate.
		
		For more general time-inhomogeneous systems, the dissipativity coefficient
		may depend explicitly on time. In such a case, the exponential contraction
		of the frozen dynamics is typically governed by an accumulated
		dissipativity factor of the form
		\[
		\exp\left\{
		-c\int_s^t \alpha(r)\,dr
		\right\},
		\]
		rather than by $e^{-c\alpha(t-s)}$. Consequently, the regularity estimates
		for the associated nonautonomous Poisson equation, and hence the
		quantitative averaging estimates, may depend on the temporal behavior of
		$\alpha(\cdot)$. In particular, quantities of the form
		\[
		\Lambda_\gamma(t)
		:=
		\int_t^\infty
		\exp\left\{
		-\gamma\int_t^r \alpha(u)\,du
		\right\}\,dr
		\]
		naturally arise in the convergence estimates.
		
		This phenomenon has been studied for time-inhomogeneous multi-scale SDEs
		in \cite{SunWangXie2024}, where the convergence rates are shown to depend
		explicitly on the time-dependent dissipativity profile. Since the main
		purpose of the present paper is to investigate the role of multiplicative
		stable noises and the corresponding strong averaging rate, we restrict
		ourselves to the uniform constant-dissipativity setting.
	\end{remark}
	\subsection{Frozen dynamics and the averaged equation}
	
	For each fixed $x\in\mathbb{R}^{d_1}$ and $s\in\mathbb{R}$, recall 
	the frozen fast equation
	\begin{equation}\label{eq:frozen-section2}
		\left\{
		\begin{aligned}
			dY_t^{s,x,y}
			={}&
			b_2(t,x,Y_t^{s,x,y})\,dt
			+\sigma_2(t,x,Y_t^{s,x,y})\,dW_t^2
			\\
			&+
			g_2(t,x,Y_{t-}^{s,x,y})\,dL_t^2,
			\qquad t\geq s,
			\\
			Y_s^{s,x,y}
			={}&y.
		\end{aligned}
		\right.
	\end{equation}
	Let $\{P_{s,t}^x\}_{s\leq t}$ denote its two-parameter transition
	family:
	\[
	P_{s,t}^xf(y)
	:=
	\mathbb{E}f(Y_t^{s,x,y}),
	\qquad
	s\leq t.
	\]
	
	A family of probability measures
	$\{\mu_t^x\}_{t\in\mathbb{R}}$ on $\mathbb{R}^{d_2}$ is called an
	evolution system of measures for $\{P_{s,t}^x\}_{s\leq t}$ if
	\[
	\int_{\mathbb{R}^{d_2}}
	P_{s,t}^xf(y)\,\mu_s^x(dy)
	=
	\int_{\mathbb{R}^{d_2}}
	f(y)\,\mu_t^x(dy)
	\]
	for every $s\leq t$ and every bounded measurable function $f$.
	
	Under Assumptions \ref{h1}--\ref{h4}, the frozen dynamics admits a
	unique evolution system of measures
	\[
	\{\mu_t^x\}_{t\in\mathbb{R}}.
	\]
	The construction and quantitative properties of this family will be
	established in Sections~3 and~4.
	
	We define the time-dependent averaged drift by
	\begin{equation}\label{eq:averaged-drift-section2}
		\bar b_1(t,x)
		:=
		\int_{\mathbb{R}^{d_2}}
		b_1(x,y)\,\mu_t^x(dy),
		\qquad
		(t,x)\in
		\mathbb{R}\times\mathbb{R}^{d_1}.
	\end{equation}
	The first averaged equation is then given by
	\begin{equation}\label{eq:first-averaged-section2}
		\left\{
		\begin{aligned}
			d\bar X_t^\varepsilon
			={}&
			\bar b_1
			\left(
			\frac{t}{\varepsilon},
			\bar X_t^\varepsilon
			\right)\,dt
			+
			\sigma_1(\bar X_t^\varepsilon)\,dW_t^1
			+
			g_1(\bar X_{t-}^\varepsilon)\,dL_t^1,
			\\
			\bar X_0^\varepsilon
			={}&x.
		\end{aligned}
		\right.
	\end{equation}
	
	As will be shown in Section~4, the averaged drift
	$\bar b_1(t,\cdot)$ is globally Lipschitz continuous uniformly in $t$.
	Hence \eqref{eq:first-averaged-section2} admits a unique strong
	solution.
	
	Notice that \eqref{eq:first-averaged-section2} is still
	$\varepsilon$-dependent because the averaged drift
	$\bar b_1(t,x)$ retains the explicit time dependence of the frozen
	fast dynamics. In Section~6, we further average this remaining
	fast-time dependence under quasi-periodic temporal structure.
	
	\subsection{Main results}
	Now we proceed to state the two main results of this paper.
	
	\begin{theorem}[Strong averaging principle]
		\label{thm:strong-averaging-section2}
		Suppose that Assumptions \ref{h1}--\ref{h4} hold. Let
		$T>0$, $(x,y)\in\mathbb{R}^{d_1}\times\mathbb{R}^{d_2}$ and
		\[
		1\leq m<\min\{\alpha_1,\alpha_2\}.
		\]
		Then there exists a constant $C_{T,m}>0$, independent of
		$\varepsilon\in(0,1]$, such that
		\begin{equation}
			\mathbb{E}
			\left[
			\sup_{0\leq t\leq T}
			|X_t^\varepsilon-\bar X_t^\varepsilon|^m
			\right]
			\leq
			C_{T,m}
			\bigl(
			1+|x|^m+|y|^m
			\bigr)
			\varepsilon^{m(1-1/\alpha_2)}.
			\label{eq:strong-rate-section2}
		\end{equation}
	\end{theorem}
	\begin{theorem}[Strong averaging under quasi-periodic temporal oscillations]
		\label{thm:quasiperiodic-averaging-section2}
		Suppose that Assumptions \ref{h1}--\ref{h4} hold and that
		\ref{h6} in Section~6 is satisfied. In particular,
		assume that the time-dependent averaged drift admits the representation
		\[
		\bar b_1(t,x)=B(\omega t,x),
		\]
		where $\omega\in\mathbb{R}^q$ satisfies the prescribed Diophantine
		condition.
		
		Define
		\[
		\widehat b_1(x)
		:=
		\int_{\mathbb{T}^q}
		B(\theta,x)\,d\theta,
		\]
		where $d\theta$ denotes the normalized Haar measure on
		$\mathbb{T}^q$, and let $\widehat X$ be the unique strong solution of
		\[
		\left\{
		\begin{aligned}
			d\widehat X_t
			={}&
			\widehat b_1(\widehat X_t)\,dt
			+
			\sigma_1(\widehat X_t)\,dW_t^1
			+
			g_1(\widehat X_{t-})\,dL_t^1,
			\\
			\widehat X_0
			={}&x.
		\end{aligned}
		\right.
		\]
		For every $T>0$ and $1\leq m<\alpha_1\wedge\alpha_2$,
		\begin{equation}
			\mathbb E\left[\sup_{0\leq t\leq T}|X_t^\varepsilon-\widehat X_t|^m\right]
			\leq C_{T,m,x,y}\varepsilon^{m\kappa}.
			\label{qpnew:main-rate-moment}
		\end{equation}
		where
		$\kappa=\min\{1-1/\alpha_2,\,1/(\nu+2)\}$.
	\end{theorem}
	
	\	\section{Well-posedness and Moment Estimates for the Frozen Equation}
	
	Recall that $L_t^i$, $i=1,2$, denote the isotropic $\alpha$-stable processes associated with $X_t^\varepsilon$ and $Y_t^\varepsilon$ respectively, the corresponding Poisson random measures are defined by \cite{Applebaum2009},
	\[
	N^i(t,A) = \sum_{s\leq t} 1_A(L_s^i - L_{s-}^i),\quad \forall A\in\mathcal{B}(\mathbb{R}^{d_i}),
	\]
	then compensated Poisson measures will be
	\[
	\tilde{N}^i(t,A) = N^i(t,A) - t\nu_i(A),
	\]
	where $\nu_i(dz) = \frac{c_{\alpha_i,d_i}}{|z|^{d_i+\alpha_i}}dz$ is symmetric L\'evy measure, $c_{\alpha_i,d_i}>0$ is constant. By L\'evy-It\^o decomposition and symmetry of $\nu_i(dz)$, we have
	\[
	L_t^i = \int_{|z|\leq 1} z\tilde{N}^i(t,dz) + \int_{|z|>1} zN^i(t,dz),
	\]
	so \eqref{eq:system-section2} with initial data $X_0^\varepsilon = x\in\mathbb{R}^{d_1}$, $Y_0^\varepsilon = y\in\mathbb{R}^{d_2}$ can be rewritten in the form of Poisson processes
	\begin{equation}\label{main1}
		\begin{split}
			dX_t^\varepsilon = &b_1(X_t^\varepsilon,Y_t^\varepsilon)dt +\sigma_1(X_t^\varepsilon)dW_t^1+ \int_{|z|\leq 1} g_1(X_{t-}^\varepsilon)z\tilde{N}^1(dt,dz) \\
			&+\int_{|z|>1} g_1(X_{t-}^\varepsilon)zN^1(dt,dz), \\
			dY_t^\varepsilon =& \frac{1}{\varepsilon}b_2(\frac{t}{\varepsilon},X_t^\varepsilon,Y_t^\varepsilon)dt +\frac{1}{\sqrt{\varepsilon}}\sigma_2(\frac{t}{\varepsilon},X_t^\varepsilon,Y_t^\varepsilon)dW_t^2 \\
			&+\frac{1}{\varepsilon^{\frac{1}{\alpha_2}}}\left( \int_{|z|\leq 1} g_2(\frac{t}{\varepsilon},X_{t-}^\varepsilon,Y_{t-}^\varepsilon)z\tilde{N}^2(dt,dz) + \int_{|z|>1} g_2(\frac{t}{\varepsilon},X_{t-}^\varepsilon,Y_{t-}^\varepsilon)zN^2(dt,dz) \right).
		\end{split}
	\end{equation}
	
	\begin{theorem}[Well-posedness of \eqref{eq:system-section2}]
		Let \ref{h1}-\ref{h4} hold, $\forall \varepsilon>0$, given any initial data $x\in\mathbb{R}^{d_1}$, $y\in\mathbb{R}^{d_2}$, there exists unique strong solution $(X_t^\varepsilon,Y_t^\varepsilon)$ to \eqref{eq:system-section2}.
	\end{theorem}
	\begin{proof}
		Under  \ref{h3} on $b_1, b_2, \sigma_1, \sigma_2,g_1 $ and $g_2$, the well-posedness of \eqref{main1} can be established by the same procedures outlined in [\citenum{Applebaum2009} ,Theorem 6.2.3, Theorem 6.2.9, Theorem 6.2.11], which lead to well-posedness of \eqref{eq:system-section2}.
	\end{proof}
	\begin{theorem}\label{theorem3.2}
		Let \ref{h1}-\ref{h4} hold. For any solution $(X_t^\varepsilon,Y_t^\varepsilon)$ to (\ref{eq:system-section2}), $\forall p\in [1,\alpha_1\land \alpha_2 )$, $t\geq 0$, $\exists  C>0$,  s.t.,
		\begin{align}
			\sup_{\varepsilon\in(0,1)} \mathbb{E}\big(\sup_{t\in [0,T]} |X_t^\varepsilon|^{p}\big) &\leq C_{p,T}(1+|x|^{p}+|y|^{p}), \label{well1}\\
			\sup_{\varepsilon\in(0,1)} \sup_{t\geq 0} \mathbb{E}|Y_t^\varepsilon|^{p} &\leq C_{p,T}(1+|y|^{p}).\label{well2}
		\end{align}
	\end{theorem}
	
	\begin{proof}
		Our proof refers to \cite{LiuRocknerSunXie2020}, [\citenum{SunXieXie2022},Lemma A.1], especially [\citenum{Yin2025} ,Theorem 3.2]. We observe that
		\begin{align*}
			X_t^\varepsilon &= x + \int_0^t b_1(X_s^\varepsilon,Y_s^\varepsilon)ds + \int_0^t\sigma_1(X_s^\varepsilon)dWs\\
			&+ \int_0^t g_1(X_{s-}^\varepsilon)\left( \int_{|z|\leq 1} z\tilde{N}^1(ds,dz) + \int_{|z|>1} zN^1(ds,dz) \right),
		\end{align*}
		Since $b_1$ is globally Lipschitz continuous with respect to $(x,y)$ and $\sigma_1$ and $g_1$ are globally Lipschitz continuous with respect to $x$, by Burkholder-Davis-Gundy's inequality (see, e.g., \cite{PeszatZabczyk2007}), we have for any $1\leq p<\alpha_1$,
		
		\begin{align*}
			\mathbb{E}\Bigl(\sup_{t\in[0,T]} |L_t^1|^p\Bigr)
			&\leq C_p\mathbb{E}\biggl[\int_{|x|\leq 1} |x|^2 N^1(T,dx)\biggr]^{p/2}
			+ C_{p,T}\int_{|x|>1} |x|^p \nu_1(dx)
			+ C_{p,T}\biggl[\int_{|x|>1} |x|\nu_1(dx)\biggr]^p \\
			&\leq C_{p,T}\biggl\{\biggl[\int_{|x|\leq 1} |x|^2\nu_1(dx)\biggr]^{p/2}
			+ \int_{|x|>1}|x|^p\nu_1(dx)
			+ \biggl[\int_{|x|>1}|x|\nu_1(dx)\biggr]^p\biggr\} \\
			&\leq C_{p,T}.
		\end{align*}
		
		By the Lipschitz continuity of $\sigma_1$ and the upper boundness of $g_1$, it is easy to see
		\begin{align*}
			\mathbb{E}\Bigl(\sup_{t\in[0,T]} |X_t^\epsilon|^p\Bigr)
			&\leq C_p|x|^p + C_{p,T}\int_0^T \mathbb{E}|X_s^\epsilon|^p ds
			+ C_{p,T}\int_0^T \mathbb{E}|Y_s^\epsilon|^p ds\\
			&+ C_{p,T}\mathbb{E}\left(\int_{0}^{T}\int_{|z|\leq 1}g_1(X_{s-}^\varepsilon)  z{N}^1(ds,dz)\right)^p + C_{p,T}\mathbb{E}\left(\int_{0}^{T}\int_{|z|>1} g_1(X_{s-}^\varepsilon)  zN^1(ds,dz)\right)^p  \\
			&\leq C_{p,T}\bigl(|x|^p \bigr) + C_{p,T}\int_0^T \mathbb{E}|X_s^\epsilon|^p ds
			+ C_{p,T}\int_0^T \mathbb{E}|Y_s^\epsilon|^p ds\\
			&+ C_{p,T}\mathbb{E}\left(\int_{0}^{T}\int_{|z|\leq 1}|g_1(X_{s-}^\varepsilon)  z|^2{N}^1(ds,dz)\right)^{\frac{p}{2}} + C_{p,T}\mathbb{E}\int_{0}^{T}\int_{|z|>1} |g_1(X_{s-}^\varepsilon)  z|^pv_1(dz)ds  \\
			&\leq C_{p,T}\bigl(|x|^p + 1\bigr) + C_{p,T}\int_0^T \mathbb{E}|X_s^\epsilon|^p ds+ C_{p,T}\int_0^T \mathbb{E}|Y_s^\epsilon|^p ds.
		\end{align*}
		
		Thus, Gronwall's inequality yields
		\begin{align}
			\mathbb{E}\Bigl(\sup_{t\in[0,T]} |X_t^\epsilon|^p\Bigr)
			\leq C_{p,T}\bigl(|x|^p + 1\bigr) + C_{p,T}\int_0^T \mathbb{E}|Y_s^\epsilon|^p ds.
		\end{align}
		Next, we estimate $Y_t^\varepsilon$. 
		From \eqref{main1}, we also deduce that
		\begin{align*}
			Y_t^\varepsilon = &y_0 + \int_0^t \frac{1}{\varepsilon}b_2(\frac{s}{\varepsilon},X_s^\varepsilon,Y_s^\varepsilon)ds + \int_0^t \frac{1}{\sqrt{\varepsilon}}\sigma_2(\frac{s}{\varepsilon},X_s^\varepsilon,Y_s^\varepsilon)dWs\\
			&+\int_0^t \frac{1}{\varepsilon^{\frac{1}{\alpha_2}}}g_2(\frac{s}{\varepsilon},X_{s-}^\varepsilon,Y_{s-}^\varepsilon)\left( \int_{|z|\leq \varepsilon^{\frac{1}{\alpha_2}}} z\tilde{N}^2(ds,dz) + \int_{|z|>\varepsilon^{\frac{1}{\alpha_2}}} zN^2(ds,dz) \right),
		\end{align*}
		we define
		\[
		U(y) = (|y|^2+1)^{\frac{p}{2}},
		\]
		we can see that $U(y)>0$, and
		\begin{align}
			|DU(y)| &= \left| \frac{p y}{(|y|^2+1)^{1-\frac{p}{2}}} \right| \leq C_{p} |y|^{p-1}, \label{DUU}\\
			|D^2U(y)| &= \left| \left( \frac{p I_{d_2\times d_2}}{(|y|^2+1)^{1-\frac{p}{2}}} - \frac{p (p-2)y\otimes y}{(|y|^2+1)^{2-\frac{p}{2}}} \right) \right| \leq \frac{C_{p} }{(|y|^2+1)^{1-\frac{p}{2}}} \leq C_{p} .\label{D2UU}
		\end{align}
		Applying It\^o's formula and taking expectation on both sides, with $\mathbb{E}(dW_t)=0$ and $\mathbb{E}(\tilde{N}^2(ds,dz))=0$, we derive,
		\begin{align*}
			\frac{d\mathbb{E}U(Y_t^\varepsilon)}{dt} &= \mathbb{E}\frac{1}{\varepsilon}\langle b_2(\frac{t}{\varepsilon},X_t^\varepsilon,Y_t^\varepsilon), DU(Y_t^\varepsilon)\rangle \\
			&\quad+\frac{1}{2\varepsilon}\mathbb{E} |\sigma_2(\frac{t}{\varepsilon},X_t^\varepsilon,Y_t^\varepsilon)\cdot \sigma_2(\frac{t}{\varepsilon},X_t^\varepsilon,Y_t^\varepsilon)^\top|\cdot|D^2U(Y_t^\varepsilon)|\\
			&\quad + \mathbb{E}\int_{|z|\leq \varepsilon^{\frac{1}{\alpha_2}}} \left( U(Y_t^\varepsilon + \varepsilon^{-\frac{1}{\alpha_2}}g_2\cdot z) - U(Y_t^\varepsilon) - \langle DU(Y_t^\varepsilon), \varepsilon^{-\frac{1}{\alpha_2}}g_2\cdot z\rangle \right)\nu_2(dz) \\
			&\quad + \mathbb{E}\int_{|z|>\varepsilon^{\frac{1}{\alpha_2}}} \left( U(Y_t^\varepsilon + \varepsilon^{-\frac{1}{\alpha_2}}g_2\cdot z) - U(Y_t^\varepsilon) \right)\nu_2(dz) \\
			&= I_1 + I_2 + I_3+I_4,
		\end{align*}
		from the upper boundness of $g_2(t,x,y)$, the proof of \eqref{well2} is quite similar to [\citenum{SunXieXie2022}, Lemma A.1] and [\citenum{Yin2025}, Theorem 3.2], for $I_1$, by  Assumption \ref{h2},
		\begin{align*}
			I_1 &= \mathbb{E}\frac{1}{\varepsilon}\langle b_2(\frac{t}{\varepsilon},X_t^\varepsilon,Y_t^\varepsilon), DU(Y_t^\varepsilon)\rangle \\
			&\leq \frac{1}{\varepsilon}\mathbb{E}\frac{\langle b_2(\frac{t}{\varepsilon},X_t^\varepsilon,Y_t^\varepsilon) - b_2(\frac{t}{\varepsilon},X_t^\varepsilon,0), p Y_t^\varepsilon\rangle + \langle b_2(\frac{t}{\varepsilon},X_t^\varepsilon,0), p Y_t^\varepsilon\rangle}{(|Y_t^\varepsilon|^2 + 1)^{1-\frac{p}{2}}} \\
			&\leq \frac{C}{\varepsilon}\mathbb{E}\frac{C|Y_t^\varepsilon| - \alpha|Y_t^\varepsilon|^2}{(|Y_t^\varepsilon|^2 + 1)^{1-\frac{p}{2}}} \\
			&\leq  \frac{C}{\varepsilon} - \frac{\alpha\mathbb{E}\bigl[U(Y_t^\varepsilon)\bigr]}{2\varepsilon},
		\end{align*}
		in addition, taking $y = \varepsilon^{-\frac{1}{\alpha_2}}z$, we obtain
		\[
		\nu_2(dz) = \frac{c}{|z|^{d_2+\alpha_2}}dz = \frac{c}{|\varepsilon^{\frac{1}{\alpha_2}}y|^{d_2+\alpha_2}}(\varepsilon^{\frac{1}{\alpha_2}})^{d_2} dy = \frac{1}{\varepsilon}\frac{c}{|y|^{d_2+\alpha_2}}dy = \frac{1}{\varepsilon}\nu_2(dy),
		\]
		thus for $I_2$, according to the upper boundness of $\sigma_2$ and \eqref{D2UU}, we have \begin{align}
			I_2 \leq \frac{ C_{{p}}}{\varepsilon},
		\end{align}
		for  $I_3$ ,by the upper boundness of $g_2(t,x,y)$ and young inequality, we have
		\begin{align*}
			I_3 &= \frac{1}{\varepsilon}\mathbb{E}\int_{|y|>1} \left( U(Y_t^\varepsilon + g_2\cdot y) - U(Y_t^\varepsilon) \right)\nu_2(dy) \\
			&\leq \frac{C_{p}}{\varepsilon}\mathbb{E}\int_{|y|>1} \left( |Y_t^\varepsilon|^{p-1}\cdot|g_2y| + |g_2y|^{p} \right)\nu_2(dy)\\ &\leq \frac{C_{p}}{\varepsilon}\mathbb{E}\int_{|y|>1} \left( |Y_t^\varepsilon|^{p-1}\cdot|g_2y| + |g_2y|^{p} \right)\nu_2(dy)\\& \leq \frac{C_{p}}{\varepsilon} + \frac{\delta \mathbb{E}U(Y_t^\varepsilon)}{\varepsilon},
		\end{align*}
		then
		\[
		I_4 = \frac{1}{\varepsilon}\mathbb{E}\int_{|y|\leq 1} \left( U(Y_t^\varepsilon + g_2\cdot y) - U(Y_t^\varepsilon) - \langle DU(Y_t^\varepsilon), g_2\cdot y\rangle \right)\nu_2(dy) \leq \frac{C_{p}}{\varepsilon},
		\]
		combining above estimates and taking $\delta=\frac{\alpha}{4}$, we derive
		\[
		\frac{d\mathbb{E}U(Y_t^\varepsilon)}{dt} \leq \frac{C}{\varepsilon} - \frac{-\alpha\mathbb{E}U(Y_t^\varepsilon)}{4\varepsilon},
		\]
		so that by Gronwall's inequality we have
		\[
		\mathbb{E}U(Y_t^\varepsilon) \leq e^{-C\frac{\alpha t}{4\varepsilon}}(|y|^2 + 1)^{\frac{p}{2}} + \frac{C}{\varepsilon}\int_0^t e^{-\frac{\alpha(t-s)}{4\varepsilon}}ds,
		\]
		which means
		\[
		\sup_{\varepsilon\in(0,1)} \sup_{t\geq 0}\mathbb{E}(|Y_t^\varepsilon|^{p}) \leq C_{p}(1 + |y|^{p}),
		\]
		and so the proof is complete.
	\end{proof}
	\begin{proposition}
		Let \ref{h1}-\ref{h4} hold, we have for $m\in[1,\alpha_1 \land \alpha_2)$, for $T\geq 1$,
		\begin{align}
			\sup_{t\geq s}\mathbb{E}|Y_t^{s,x,y}|^m &\leq C_m(1+|y|^m), \label{frozen1}\\
			\mathbb{E}\left( \sup_{t\in[s,s+T]}|Y_t^{s,x,y}|^m \right) &\leq C_m\left(T^{\frac{m}{\alpha_2}}+|y|^m\right). \label{frozen2}
		\end{align}
	\end{proposition}
	\begin{proof}
		The proof follows from [\citenum{SunXieXie2022}, Lemma A.1], [\citenum{Yin2025}, Theorem 4.1], (\ref{frozen1}) can be derived from (\ref{well2}) and we shall prove (\ref{frozen2}). We define
		\[
		U_T(y) = (|y|^2+T^{\frac{2}{\alpha_2}})^{\frac{m}{2}},
		\]
		so that similar to (\ref{DUU}) and (\ref{D2UU}),
		\begin{align}
			|DU_T(y)| &= \left| \frac{my}{(|y|^2+T^{\frac{2}{\alpha_2}})^{1-\frac{m}{2}}} \right| \leq C_m|y|^{m-1}, \\
			|D^2U_T(y)| &= \left| \frac{mI_{d_2\times d_2}}{(|y|^2+T^{\frac{2}{\alpha_2}})^{1-\frac{m}{2}}} - \frac{m(m-2)y\otimes y}{(|y|^2+T^{\frac{2}{\alpha_2}})^{2-\frac{m}{2}}} \right| \leq C_mT^{\frac{2}{\alpha_2}(\frac{m}{2}-1)},\label{hhhhh}
		\end{align}
		by It\^o's formula,
		\begin{align*}
			U_T(Y_t^{s,x,y}) 
			&=  U_T(y)+\int_s^t \langle b_2(r,x,Y_r^{s,x,y}), DU_T(Y_r^{s,x,y}) \rangle dr \\
			&\quad+\frac{1}{2}\int_s^t Tr[\sigma_2(r,x,Y_r^{s,x,y}) \sigma_2(r,x,Y_r^{s,x,y})^\top D^2U_T(Y_r^{s,x,y})]dr\\
			&\quad+ \int_s^t \int_{|z|\leq T^{\frac{1}{\alpha_2}}} \left( U_T(Y_r^{s,x,y}+g_2z) - U_T(Y_r^{s,x,y}) \right)\tilde{N}^2(dr,dz) \\
			&\quad + \int_s^t \int_{|z|\leq T^{\frac{1}{\alpha_2}}} \left( U_T(Y_r^{s,x,y}+g_2z) - U_T(Y_r^{s,x,y}) - \langle DU_T(Y_r^{s,x,y}), g_2z \rangle \right)\nu_2(dz)dr \\
			&\quad + \int_s^t \int_{|z|>T^{\frac{1}{\alpha_2}}} \left( U_T(Y_r^{s,x,y}+g_2z) - U_T(Y_r^{s,x,y}) \right)N^2(dr,dz) + U_T(y) \\
			& \quad+ \int_{s}^{t} \left\langle D U_T \big(Y_r^{s,x,y}\big), \sigma_2\big(r, x, Y_r^{s,x,y}\big) dW_r^2 \right\rangle\\
			&= U_T(y)+ \hat{I}_1 + \hat{I}_2 + \hat{I}_3 + \hat{I}_4 +\hat{I}_5+\hat{I}_6,
		\end{align*}
		so by Assumption \ref{h2} of $b_2(t,x,y)$ and $T\geq 1$, obviously we have for $\hat{I}_1$,
		\[
		\mathbb{E}\left( \sup_{t\in[s,s+T]} \hat{I}_1(t) \right) \leq \int_s^{s+T }\frac{C_m}{(|Y_r^{s,x,y}|^2+T^{\frac{2}{\alpha_2}})^{1-\frac{m}{2}}} dr \leq C_mT^{\frac{m}{\alpha_2}-\frac{2}{\alpha_2}+1} ,
		\]
		for $\hat{I}_2$,by Assumption \ref{h3} the boundness of $\sigma_2$
		and \eqref{hhhhh}, we have that
		\begin{align}
			\mathbb{E}\left( \sup_{t\in[s,s+T]} |\hat{I}_2(t)| \right) \leq C_mT^{\frac{m}{\alpha_2}-\frac{2}{\alpha_2}+1},
		\end{align}
		for $\hat{I}_3$, from Assumption \ref{h3} of $g_2(t,x,y)$, by Burkholder-Davies-Gundy's inequality and (\ref{frozen1}),
		\begin{align*}
			&\mathbb E\left(
			\sup_{t\in[s,s+T]}|\hat I_3(t)|
			\right)
			\\
			&\leq C_m\mathbb E\left[
			\left(
			\int_s^{s+T}\int_{|z|\leq T^{1/\alpha_2}}
			\left|
			U_T\left(
			Y_{r-}^{s,x,y}+g_2(r,x,Y_{r-}^{s,x,y})z
			\right)
			-U_T\left(Y_{r-}^{s,x,y}\right)
			\right|^2
			N^2(dr,dz)
			\right)^{1/2}
			\right]
			\\
			&\leq C_m\mathbb E\left[
			\left(
			\int_s^{s+T}\int_{|z|\leq T^{1/\alpha_2}}
			\left(
			|Y_{r-}^{s,x,y}|^{2m-2}|z|^2+|z|^{2m}
			\right)
			N^2(dr,dz)
			\right)^{1/2}
			\right]
			\\
			&\leq C_m\mathbb E\left[
			\left(
			\sup_{r\in[s,s+T]}|Y_r^{s,x,y}|^{m-1}
			\right)
			\left(
			\int_s^{s+T}\int_{|z|\leq T^{1/\alpha_2}}
			|z|^2\,N^2(dr,dz)
			\right)^{1/2}
			\right]
			\\
			&\quad
			+ C_m\mathbb E\left[
			\left(
			\int_s^{s+T}\int_{|z|\leq T^{1/\alpha_2}}
			|z|^{2m}\,N^2(dr,dz)
			\right)^{1/2}
			\right]
			\\
			&\leq \frac14\mathbb E\left(
			\sup_{r\in[s,s+T]}|Y_r^{s,x,y}|^m
			\right)
			+ C_m\mathbb E\left[
			\left(
			\int_s^{s+T}\int_{|z|\leq T^{1/\alpha_2}}
			|z|^2\,N^2(dr,dz)
			\right)^{m/2}
			\right]
			\\
			&\quad
			+ C_m\mathbb E\left[
			\left(
			\int_s^{s+T}\int_{|z|\leq T^{1/\alpha_2}}
			|z|^{2m}\,N^2(dr,dz)
			\right)^{1/2}
			\right]
			\\
			&\leq \frac14\mathbb E\left(
			\sup_{r\in[s,s+T]}|Y_r^{s,x,y}|^m
			\right)
			+ C_m\left(
			T\int_{|z|\leq T^{1/\alpha_2}}
			|z|^2\,\nu_2(dz)
			\right)^{m/2}
			\\
			&\quad
			+ C_m\left(
			T\int_{|z|\leq T^{1/\alpha_2}}
			|z|^{2m}\,\nu_2(dz)
			\right)^{1/2}
			\\
			&\leq \frac14\mathbb E\left(
			\sup_{r\in[s,s+T]}|Y_r^{s,x,y}|^m
			\right)
			+ C_m T^{m/\alpha_2}.
		\end{align*}
		Here we used Jensen's inequality and the identities
		\[
		\int_{|z|\leq T^{1/\alpha_2}}
		|z|^2\,\nu_2(dz)
		=C T^{(2-\alpha_2)/\alpha_2},
		\qquad
		\int_{|z|\leq T^{1/\alpha_2}}
		|z|^{2m}\,\nu_2(dz)
		=C_m T^{(2m-\alpha_2)/\alpha_2}.
		\]
		Both integrals are finite since $2m\geq 2>\alpha_2$.
		For $\hat{I}_4$,
		\[
		\mathbb{E}\left( \sup_{t\in[s,s+T]} |\hat{I}_4(t)| \right) \leq C_mT^{\frac{2}{\alpha_2}(\frac{m}{2}-1)} \int_s^{s+T} \int_{|z|\leq T^{\frac{1}{\alpha_2}}} |g_2\cdot z|^2\nu_2(dz)dr \leq C_mT^{\frac{m}{\alpha_2}},
		\]
		For $\hat I_5$, the total variation bound, the compensation
		formula, and the boundedness of $g_2$ yield
		\begin{align*}
			&\mathbb E\left(
			\sup_{t\in[s,s+T]}|\hat I_5(t)|
			\right)
			\\
			&\leq \mathbb E\left[
			\int_s^{s+T}\int_{|z|>T^{1/\alpha_2}}
			\left|
			U_T\left(
			Y_{r-}^{s,x,y}+g_2(r,x,Y_{r-}^{s,x,y})z
			\right)
			-U_T\left(Y_{r-}^{s,x,y}\right)
			\right|
			N^2(dr,dz)
			\right]
			\\
			&\leq C_m
			\int_s^{s+T}\int_{|z|>T^{1/\alpha_2}}
			\left(
			\mathbb E|Y_r^{s,x,y}|^{m-1}|z|+|z|^m
			\right)
			\nu_2(dz)\,dr
			\\
			&\leq C_m\mathbb E\left(
			\sup_{r\in[s,s+T]}|Y_r^{s,x,y}|^{m-1}
			\right)
			\left(
			T\int_{|z|>T^{1/\alpha_2}}|z|\,\nu_2(dz)
			\right)
			\\
			&\quad
			+ C_m T\int_{|z|>T^{1/\alpha_2}}
			|z|^m\,\nu_2(dz)
			\\
			&\leq \frac14\mathbb E\left(
			\sup_{r\in[s,s+T]}|Y_r^{s,x,y}|^m
			\right)
			+ C_m\left(
			T\int_{|z|>T^{1/\alpha_2}}
			|z|\,\nu_2(dz)
			\right)^m
			\\
			&\quad
			+ C_m T\int_{|z|>T^{1/\alpha_2}}
			|z|^m\,\nu_2(dz)
			\\
			&\leq \frac14\mathbb E\left(
			\sup_{r\in[s,s+T]}|Y_r^{s,x,y}|^m
			\right)
			+ C_m T^{m/\alpha_2}.
		\end{align*}
		Indeed,
		\[
		\int_{|z|>T^{1/\alpha_2}}
		|z|\,\nu_2(dz)
		=C T^{(1-\alpha_2)/\alpha_2},
		\qquad
		\int_{|z|>T^{1/\alpha_2}}
		|z|^m\,\nu_2(dz)
		=C_m T^{(m-\alpha_2)/\alpha_2},
		\]
		which are finite because $1\leqslant m<\alpha_2$.

		The applications of Young's inequality above are for $m>1$.
		When $m=1$, the factors involving the zeroth power are
		interpreted as one, and the same final bounds follow directly.
		All estimates are first carried out after localization.
		combining above estimates we derive (\ref{frozen2}).
	\end{proof}
	\begin{lemma}\label{lemma3.4}
		$\forall t\in[0,T]$, $T\geqslant 1$, we have
		\[
		\mathbb{E}\left( \sup_{t\in[0,T]} |Y_t^\varepsilon|^m \right) \leq C_{T,m}\left( \varepsilon^{-\frac{m}{\alpha_2}} + |y|^m \right).
		\]
	\end{lemma}
	\begin{proof}
		Denote $
		\widetilde{Y}_r^\varepsilon := Y_{\varepsilon r}^\varepsilon, \quad
		\widetilde{W}_r^2 := \varepsilon^{-1/2}W_{\varepsilon r}^2,\quad
		\widetilde{L}_r^2 := \varepsilon^{-1/\alpha_2}L_{\varepsilon r}^2
		$
		, so that
		\begin{align*}
			\tilde{Y}_t^\varepsilon &= y +\frac{1}{\varepsilon}\int_0^{t\varepsilon} b_2(\frac{s}{\varepsilon},X_s^\varepsilon,{Y}_s^\varepsilon)ds+ \frac{1}{\sqrt{\varepsilon}}\int_0^{t\varepsilon} \sigma_2(\frac{s}{\varepsilon},X_s^\varepsilon,{Y}_s^\varepsilon)dW_s^2 +\frac{1}{\varepsilon^{\frac{1}{\alpha_2}}}\int_{0}^{t\varepsilon} g_2(\frac{s}{\varepsilon},X_{s-}^\varepsilon,{Y}_{s-}^\varepsilon)dL_{s}^2 \\
			&= y + \int_0^t b_2(s,X_{s\varepsilon}^\varepsilon,\tilde{Y}_{s}^\varepsilon)ds +\int_0^{t} \sigma_2(s,X_{s\varepsilon}^\varepsilon,\tilde{Y}_s^\varepsilon)d\widetilde{W}_s^2+\int_{0}^{t} g_2(s,X_{(s\varepsilon)-}^\varepsilon,\tilde{Y}_{s-}^\varepsilon)d\widetilde{L}_s^2.
		\end{align*}
		With respect to the filtration
		$(\mathcal{F}_{\varepsilon r})_{r\geq 0}$,
		$\widetilde W^2$ and $\widetilde L^2$ are, respectively,
		a Brownian motion and an isotropic $\alpha_2$-stable
		L\'evy process with the same laws as the original noises.
		By \ref{h1} and \ref{h2}, the boundedness of $\sigma_2$ and $g_2$,
		the dissipativity condition, and the bound on $b_2(t,x,0)$
		are uniform in the slow variable $x$. Then similar to the proof of (\ref{frozen2}), we have
		\[
		\mathbb{E}\left( \sup_{t\in[0,T]} |\tilde{Y}_t^\varepsilon|^m \right) \leq C_m\left( T^{\frac{m}{\alpha_2}} + |y|^m \right),
		\]
		from \ref{h2}, and (\ref{frozen2}), for any $T\geq 1$,
		\[
		\mathbb{E}\left( \sup_{t\in[0,T]} |Y_t^\varepsilon|^m \right) = \mathbb{E}\left( \sup_{t\in[0,\frac{T}{\varepsilon}]} |\tilde{Y}_t^\varepsilon|^m \right) \leq C_m\left( \left(\frac{T}{\varepsilon}\right)^{\frac{m}{\alpha_2}} + |y|^m \right) \leq C_{T,m}\left( \varepsilon^{-\frac{m}{\alpha_2}} + |y|^m \right).\qedhere
		\]
	\end{proof}
	\begin{remark}
		Lemma \ref{lemma3.4} is only used, but very important, for strong convergence estimates, see Theorem \ref{theorem5.4} for details.
	\end{remark}
	\section{Lyapunov Estimates and Regularity of the Poisson Equation}
	Let $\mathcal{L}_2(t,x, y)$ be the infinitesimal generator of the frozen process $Y_t^{x,y}$ defined by
	
	\begin{align}\label{frozen}
		dY^{x,y}_t = b_2(t,x, Y^{x,y}_t)dt + \sigma_2(t,x, Y^{x,y}_t)dW_t^2 + g_2(t,x, Y^{x,y}_t)dL_t^2,
	\end{align}
	
	where $L_t^2$ is a $d_2$-dimensional isotropic $\alpha_2$-stable process with L\'evy measure 
	\begin{align*}
		\nu_2(dz) = c_{\alpha_2, d_2} |z|^{-d_2 - \alpha_2} dz.
	\end{align*}
	
	For a sufficiently smooth function $\Phi: \mathbb{R}^{d_1} \times \mathbb{R}^{d_2} \to \mathbb{R}$, the generator is given by
	
	\begin{align}\label{4.2}
		\begin{aligned}
			\mathcal{L}_2(t,x, y)\Phi(t,x, y) :=& \langle b_2(t,x, y), \nabla_y \Phi(x, y) \rangle \\
			&+ \frac{1}{2} \sum_{i,j=1}^{d_2} \left( \sigma_2(t,x, y) \sigma_2(t,x, y)^\top \right)_{ij} \frac{\partial^2 \Phi}{\partial y_i \partial y_j}(t,x, y) \\
			&+ \int_{\mathbb{R}^{d_2}} \Big[ \Phi(t,x, y + g_2(t,x, y) z) - \Phi(t,x, y) \\
			&\qquad - \mathbf{1}_{\{|z| \leq 1\}} \langle \nabla_y \Phi(t,x, y), g_2(t,x, y) z \rangle \Big] \nu_2(dz).
		\end{aligned}
	\end{align}
	
	Since the coefficients of the fast equation depend on time, the frozen equation does not possess an invariant measure. In this case, the measure associated with the solution of the equation is time-dependent, which we denote by $\mu^{x}_t$. Below we establish the exponential dependence of solutions to the frozen equation on their initial values. Due to the presence of multiplicative $\alpha$-stable noise, the classical It{\^o}'s
	formula cannot be directly employed for our computations. According to \cite{YangYin2026,LiuWang2019,Wang2018}, we can normalize the multiplicative $\alpha$-stable noise. For the infinitesimal generator $\mathcal{L}_2$ defined in (\ref{4.2}), we have an equivalent variant form by above analysis
	\begin{align}\label{42}
		\begin{aligned}
			\mathcal{L}_2(t,x, y)\Phi(t,x, y) :=& \langle b_2(t,x, y), \nabla_y \Phi(x, y) \rangle \\
			&+ \frac{1}{2} \sum_{i,j=1}^{d_2} \left( \sigma_2(t,x, y) \sigma_2(t,x, y)^\top \right)_{ij} \frac{\partial^2 \Phi}{\partial y_i \partial y_j}(t,x, y) \\
			&+ \int_{\mathbb{R}^{d_2}} \Big[ \Phi(t,x, y + z) - \Phi(t,x, y) \\
			&\qquad - \mathbf{1}_{\{|z| \leq 1\}} \langle \nabla_y \Phi(t,x, y),  z \rangle \Big]H(t,y,\hat{z}) \nu_2(dz).
		\end{aligned}
	\end{align}
	where $H(t,y,\hat{z})=\big(|det(g_2(t,x,y))|\cdot|g_2^{-1}(t,x,y)\hat{z}|^{\alpha_2+d_2}\big)^{-1}$.
	\subsection{Lyapunov Estimate}
	Next, under Assumption \ref{h1} we apply coupling method for nonlocal operators to obtain exponential contractivity with respect to $L^p$-Wasserstein distance, then the exponential dependence of solutions to the frozen equation on the initial data follows.
	
	The idea is based on \cite{LiuWang2019,Wang2018}. We introduce a Markov coupling operator for $\mathcal{L}_2$ in terms of \eqref{4.2}, $\forall y_1,y_2,z\in\mathbb{R}^{d_2}$, let
	\[
	\varphi_{y_1,y_2}(z)=
	\begin{cases}
		z-\frac{2\langle y_1-y_2,z\rangle}{|y_1-y_2|^2}(y_1-y_2), & y_1\neq y_2;\\
		-z, & y_1=y_2.
	\end{cases}
	\]
	then $\varphi_{y_1,y_2}(z):\mathbb{R}^{d_2}\to\mathbb{R}^{d_2}$ enjoys the following properties:
	\begin{itemize}
		\item[(A1)] $\varphi_{y_1,y_2}(z)=\varphi_{y_2,y_1}(z)$ and $\varphi_{y_1,y_2}^2(z)=z$, i.e., $\varphi_{y_1,y_2}^{-1}=\varphi_{y_1,y_2}(z)$;
		\item[(A2)] $|\varphi_{y_1,y_2}(z)|=|z|$;
		\item[(A3)] $(z-\varphi_{y_1,y_2}(z))\parallel (y_1-y_2)$ and $(z+\varphi_{y_1,y_2}(z))\perp(y_1-y_2)$.
	\end{itemize}
	For convenient we denote $\varphi_{y_1,y_2}(z)$ by $\varphi(z)$. $\forall y_1,y_2,z\in\mathbb{R}^{d_2}$, let
	\[
	\widetilde{H}(t,y_1,y_2,z)=H(t,y_1,z)\wedge H(t,y_2,z)\wedge H(t,y_1,\varphi(z))\wedge H(t,y_2,\varphi(z)),
	\]
	and $\forall F\in C_b^2(\mathbb{R}^{2d_2})$,
	\[
	\nabla_{y_1}F(y_1,y_2):=\left(\frac{\partial F(y_1,y_2)}{\partial y_1^i}\right)_{1\leq i\leq d_2},\quad
	\nabla_{y_2}F(y_1,y_2):=\left(\frac{\partial F(y_1,y_2)}{\partial y_2^i}\right)_{1\leq i\leq d_2},
	\]
	For $a\in(0,\frac12)$, we can define for $z\in\mathbb{R}^{d_2}\setminus\{0\}$
	\begin{align*}
		\widetilde{\mathcal{L}} F(y_1,y_2)
		&=\frac12\Bigg[
		\int_{\{|z|\leq a|y_1-y_2|\}}
		\Big(
		F(y_1+z,y_2+\varphi(z))-F(y_1,y_2)
		-\langle\nabla_{y_1}F(y_1,y_2),z\rangle I_{\{|z|\leq1\}}
		\\
		&\qquad -\langle\nabla_{y_2}F(y_1,y_2),\varphi(z)\rangle I_{\{|z|\leq1\}}
		\Big)
		\widetilde{H}(t,y_1,y_2,z)\nu_2(dz)
		\\
		&\quad +\int_{\{|z|\leq a|y_1-y_2|\}}
		\Big(
		F(y_1+\varphi(z),y_2+z)-F(y_1,y_2)
		-\langle\nabla_{y_2}F(y_1,y_2),z\rangle I_{\{|z|\leq1\}}
		\\
		&\qquad -\langle\nabla_{y_1}F(y_1,y_2),\varphi(z)\rangle I_{\{|z|\leq1\}}
		\Big)
		\widetilde{H}(t,y_1,y_2,z)\nu_2(dz)\Bigg]
		\\
		&\quad +\int_{\{|z|\leq a|y_1-y_2|\}}
		\Big(
		F(y_1+z,y_2+z)-F(y_1,y_2)
		-\langle\nabla_{y_1}F(y_1,y_2),z\rangle I_{\{|z|\leq1\}}
		\\
		&\qquad -\langle\nabla_{y_2}F(y_1,y_2),z\rangle I_{\{|z|\leq1\}}
		\Big)
		\times (H(t,y_1,z)\wedge H(t,y_2,z)-\widetilde{H}(t,y_1,y_2,z))\nu_2(dz)
		\\
		&\quad +\int_{\{|z|>a|y_1-y_2|\}}
		\Big(
		F(y_1+z,y_2+z)-F(y_1,y_2)
		-\langle\nabla_{y_1}F(y_1,y_2),z\rangle I_{\{|z|\leq1\}}
		\\
		&\qquad -\langle\nabla_{y_2}F(y_1,y_2),z\rangle I_{\{|z|\leq1\}}
		\Big)
		\times (H(t,y_1,z)\wedge H(t,y_2,z))\nu_2(dz)
		\\
		&+\int_{\mathbb R^{d_2}\setminus\{0\}}
		\Bigl[
		F(y_1+z,y_2)-F(y_1,y_2)
		-\bigl\langle\nabla_{y_1}F(y_1,y_2),z\bigr\rangle
		I_{\{|z|\le1\}}
		\Bigr]
		\\
		&\qquad\times
		\Bigl[
		H(t,y_1,z)
		-\bigl(H(t,y_1,z)\wedge H(t,y_2,z)\bigr)
		\Bigr]\,\nu_2(dz)
		\\[4pt]
		&+\int_{\mathbb R^{d_2}\setminus\{0\}}
		\Bigl[
		F(y_1,y_2+z)-F(y_1,y_2)
		-\bigl\langle\nabla_{y_2}F(y_1,y_2),z\bigr\rangle
		I_{\{|z|\le1\}}
		\Bigr]
		\\
		&\qquad\times
		\Bigl[
		H(t,y_2,z)
		-\bigl(H(t,y_1,z)\wedge H(t,y_2,z)\bigr)
		\Bigr]\,\nu_2(dz)\\
		&\quad \quad+\frac12\operatorname{Tr}\Biggl[
		\begin{pmatrix}
			\sigma_2(t,x,y_1)^\top &
			\sigma_2(t,x,y_2)^\top
		\end{pmatrix}
		D^2F(y_1,y_2)
		\begin{pmatrix}
			\sigma_2(t,x,y_1)\\
			\sigma_2(t,x,y_2)
		\end{pmatrix}
		\Biggr]\\
		&\quad \quad +\langle\nabla_{y_1}F(y_1,y_2),b_2(t,x,y_1)\rangle
		+\langle\nabla_{y_2}F(y_1,y_2),b_2(t,x,y_2)\rangle.
	\end{align*}
	From above definition we have the following lemma. The proof is analogous to that in [\citenum{LiuWang2019}, Lemma 2.2] and \cite{YangYin2026}   and is omitted.
	
	\begin{lemma}\label{lem:4.1}
		$\forall z\in\mathbb{R}^{d_2}\setminus\{0\}$, $\widetilde{\mathcal{L}}$ is coupling operator of $\mathcal{L}_2$, i.e., $\forall F(y_1)\in C_b^2(\mathbb{R}^{d_2})$
		\begin{equation}\label{eq4.4}
			\widetilde{\mathcal{L}}F(y_1,y_2)=\mathcal{L}_2F(y_1)+\mathcal{L}_2F(y_2).
		\end{equation}
		where $F(y_1,y_2)=F(y_1)+F(y_2)$.
	\end{lemma}

	Then, we need to construct a coupling process associated with $\widetilde{\mathcal{L}}$, this method is based on [\citenum{Wang2018}, Section 2.2]. let
	For any $y_1\ne y_2$ and
	$A\in\mathcal B(\mathbb R^{2d_2}\setminus\{0\})$, define
	\begin{align*}
		&\mu(t,y_1,y_2,A)
		\\
		&=\frac12
		\int_{\{(z,\varphi(z))\in A,\ |z|\leq a|y_1-y_2|\}}
		\widetilde H(t,y_1,y_2,z)\,\nu_2(dz)
		\\
		&\quad+\frac12
		\int_{\{(\varphi(z),z)\in A,\ |z|\leq a|y_1-y_2|\}}
		\widetilde H(t,y_1,y_2,z)\,\nu_2(dz)
		\\
		&\quad+
		\int_{\{(z,z)\in A,\ |z|\leq a|y_1-y_2|\}}
		\Bigl[
		H(t,y_1,z)\wedge H(t,y_2,z)
		-\widetilde H(t,y_1,y_2,z)
		\Bigr]\,\nu_2(dz)
		\\
		&\quad+
		\int_{\{(z,z)\in A,\ |z|>a|y_1-y_2|\}}
		\bigl[
		H(t,y_1,z)\wedge H(t,y_2,z)
		\bigr]\,\nu_2(dz)
		\\
		&\quad+
		\int_{\{(z,0)\in A\}}
		\Bigl[
		H(t,y_1,z)
		-\bigl(H(t,y_1,z)\wedge H(t,y_2,z)\bigr)
		\Bigr]\,\nu_2(dz)
		\\
		&\quad+
		\int_{\{(0,z)\in A\}}
		\Bigl[
		H(t,y_2,z)
		-\bigl(H(t,y_1,z)\wedge H(t,y_2,z)\bigr)
		\Bigr]\,\nu_2(dz).
	\end{align*}
	Here $a=\frac14$, and all integrals are over nonzero jump vectors.
	
	On the diagonal, define
	\[
	\mu(t,y,y,A)
	=
	\int_{\mathbb R^{d_2}\setminus\{0\}}
	\mathbf 1_A(z,z)\,H(t,y,z)\,\nu_2(dz).
	\]
	then $\forall F\in C_b^2(\mathbb{R}^{2d_2})$,
	For any $F\in C_b^2(\mathbb R^{2d_2})$, we have
	\begin{align*}
		\widetilde{\mathcal L}F(y_1,y_2)
		={}&
		\int_{\mathbb R^{2d_2}\setminus\{0\}}
		\Bigl[
		F(y_1+u_1,y_2+u_2)-F(y_1,y_2)
		\\
		&\qquad
		-\Bigl\langle
		\bigl(\nabla_{y_1}F(y_1,y_2),
		\nabla_{y_2}F(y_1,y_2)\bigr),
		(u_1,u_2)
		\Bigr\rangle
		I_{\{|u_1|\le1,\ |u_2|\le1\}}
		\Bigr]
		\\
		&\qquad\qquad\times
		\mu(t,y_1,y_2,du_1,du_2)
		\\
		&+\langle\nabla_{y_1}F(y_1,y_2),b_2(t,x,y_1)\rangle
		\\
		&+\langle\nabla_{y_2}F(y_1,y_2),b_2(t,x,y_2)\rangle
		\\
		&+\frac12\operatorname{Tr}\Biggl[
		\begin{pmatrix}
			\sigma_2(t,x,y_1)^\top &
			\sigma_2(t,x,y_2)^\top
		\end{pmatrix}
		D^2F(y_1,y_2)
		\begin{pmatrix}
			\sigma_2(t,x,y_1)\\
			\sigma_2(t,x,y_2)
		\end{pmatrix}
		\Biggr].
	\end{align*}
	Moreover, for any $h\in C_b(\mathbb R^{2d_2})$,
	writing $u=(u_1,u_2)$, we have
	\begin{align*}
		&\int_{\mathbb R^{2d_2}\setminus\{0\}}
		h(u)\frac{|u|^2}{1+|u|^2}
		\,\mu(t,y_1,y_2,du)
		\\
		&=
		\int_{\{|z|\leq a|y_1-y_2|\}}
		h(z,\varphi(z))
		\frac{|z|^2}{1+2|z|^2}
		\widetilde H(t,y_1,y_2,z)\,\nu_2(dz)
		\\
		&\quad+
		\int_{\{|z|\leq a|y_1-y_2|\}}
		h(\varphi(z),z)
		\frac{|z|^2}{1+2|z|^2}
		\widetilde H(t,y_1,y_2,z)\,\nu_2(dz)
		\\
		&\quad+
		2\int_{\{|z|\leq a|y_1-y_2|\}}
		h(z,z)\frac{|z|^2}{1+2|z|^2}
		\\
		&\qquad\qquad\times
		\Bigl[
		H(t,y_1,z)\wedge H(t,y_2,z)
		-\widetilde H(t,y_1,y_2,z)
		\Bigr]\,\nu_2(dz)
		\\
		&\quad+
		2\int_{\{|z|>a|y_1-y_2|\}}
		h(z,z)\frac{|z|^2}{1+2|z|^2}
		\\
		&\qquad\qquad\times
		\bigl[H(t,y_1,z)\wedge H(t,y_2,z)\bigr]
		\,\nu_2(dz)
		\\
		&\quad+
		\int_{\mathbb R^{d_2}\setminus\{0\}}
		h(z,0)\frac{|z|^2}{1+|z|^2}
		\\
		&\qquad\qquad\times
		\Bigl[
		H(t,y_1,z)
		-\bigl(H(t,y_1,z)\wedge H(t,y_2,z)\bigr)
		\Bigr]\,\nu_2(dz)
		\\
		&\quad+
		\int_{\mathbb R^{d_2}\setminus\{0\}}
		h(0,z)\frac{|z|^2}{1+|z|^2}
		\\
		&\qquad\qquad\times
		\Bigl[
		H(t,y_2,z)
		-\bigl(H(t,y_1,z)\wedge H(t,y_2,z)\bigr)
		\Bigr]\,\nu_2(dz),
	\end{align*}
	which means that $(y_1,y_2)\mapsto\int_{\mathbb{R}^{2d_2}}h(u)\frac{|u|^2}{1+|u|^2}\mu(y_1,y_2,du)$ is continuous on $\mathbb{R}^{2d_2}$, we also have the fact that $f(y)$ is continuous on $\mathbb{R}^{d_2}$.
	
	Then similar to analysis in \cite{LiuWang2019} and [\citenum{Wang2018}, Section 2.2], there exists a martingale solution with respect to $\widetilde{\mathcal{L}}$, i.e., there exist a probability space $(\widetilde{\Omega},\widetilde{\mathscr{F}},(\widetilde{\mathscr{F}}_t)_{t\geq0},\widetilde{\mathbb{P}})$ and $\mathbb{R}^{2d_2}$-valued process $(\widetilde{Y}_t)_{t\geq0}$ which is $(\widetilde{\mathscr{F}}_t)_{t\geq0}$-progressively measurable, and $\forall f\in C_b^2(\mathbb{R}^{2d_2})$,
	\[
	f(\widetilde{Y}_t)-f(\widetilde{Y}_s)-\int_s^{t\wedge\tau}\widetilde{\mathcal{L}}f(\widetilde{Y}_u)du,\quad t\geq s,
	\]
	is $(\widetilde{\mathscr{F}}_t)_{t\geq s}$-local martingale, where $\tau$ is the explosion time of $\widetilde{Y}_t$,
	\[
	\tau=\lim_{n\to\infty}\inf\{t\geq s:|\widetilde{Y}_t|\geq n\},
	\]
	here $\widetilde{Y}_t=(Y_t',Y_t'')$, then $Y_t'$ and $Y_t''$ are two stochastic processes on $\mathbb{R}^{d_2}$. Since $\widetilde{\mathcal{L}}$ is coupling operator of $\mathcal{L}_2$, the infinitesimal generators of processes $Y_t'$ and $Y_t''$ are $\mathcal{L}_2$, then both processes are solutions to the martingale problem of $\mathcal{L}_2$. From \ref{h3} and \ref{h4}, we can derive the well-posedness of the pathwise unique strong solution to \eqref{frozen}, so the well-posedness of weak solution follows, and the well-posedness of martingale problem to $(\mathcal{L}_2,C_b^2)$. We conclude that $Y_t'$ and $Y_t''$ are well-posed, and $\tau=\infty$, $\widetilde{\mathcal{L}}$ generates a non-explosive process $\widetilde{Y}_t$. Let $T$ be the coupling time of $Y_t'$ and $Y_t''$, i.e.
	\[
	T=\inf\{t\geq0:Y_t'=Y_t''\},
	\]
	so $T$ is an $(\widetilde{\mathscr{F}}_t)_{t\geq0}$-stopping time. Define a new process $(Z_t')_{t\geq0}$ as follows
	\[
	Z_t'=
	\begin{cases}
		Y_t'', & t<T;\\
		Y_t', & t\geq T,
	\end{cases}
	\]
	from [\citenum{Wang2018}, Section 2.2], we conclude that $(Y_t',Z_t')_{t\geq0}$ is also non-explosive coupling process of $Y_t'$ s.t. $Z_t'=Y_t'$ $\forall t\geq T$, the generator of $(Y_t',Z_t')_{t\geq0}$ before the coupling time $T$ is the coupling operator $\widetilde{\mathcal{L}}$ defined in \eqref{eq4.4}, and we know that $\forall y_1,y_2\in\mathbb{R}^{d_2}$, $f\in C_b(\mathbb{R}^{d_2})$,
	\[
	P_tf(y_1)=\mathbb{E}^{y_1}f(Y_t')=\widetilde{\mathbb{E}}^{y_1,y_2}f(Y_t'),\quad P_tf(y_2)=\mathbb{E}^{y_2}f(Z_t')=\widetilde{\mathbb{E}}^{(y_1,y_2)}f(Z_t').
	\]
	
	Then we discuss the exponential contractivity of $Y_t^{x,y}$ with respect to $L^1$-Wasserstein distance. 
	Fix $x\in\mathbb R^{d_1}$ and, for notational simplicity, write
	\[
	b_t(y)=b_2(t,x,y),\qquad
	\sigma_t(y)=\sigma_2(t,x,y),\qquad
	H_t(y,z)=H(t,x,y,z).
	\]
	Under \ref{h1}--\ref{h3}, there exist constants
	$L_\sigma,M_\sigma,h_0,h_1,L_H$, independent of $t$ and $x$, such that
	\begin{align}
		\langle b_t(y_1)-b_t(y_2),y_1-y_2\rangle
		&\leq-\alpha|y_1-y_2|^2,\label{eq:lyap-drift-bound}\\
		\|\sigma_t(y_1)-\sigma_t(y_2)\|_{\mathrm{HS}}
		&\leq L_\sigma|y_1-y_2|,
		\qquad \|\sigma_t(y)\|_{\mathrm{HS}}\leq M_\sigma,
		\label{eq:lyap-sigma-bound}\\
		0<h_0\leq H_t(y,z)&\leq h_1,
		\qquad H_t(y,-z)=H_t(y,z),\label{eq:lyap-H-bound}\\
		|H_t(y_1,z)-H_t(y_2,z)|&\leq L_H|y_1-y_2|.
		\label{eq:lyap-H-Lipschitz}
	\end{align}
	Here $\|\cdot\|_{\mathrm{HS}}$ denotes the Hilbert--Schmidt norm.
	Indeed, \eqref{eq:lyap-drift-bound} follows from \ref{h2} by discarding its
	nonnegative diffusion term. The operator-norm bounds in \ref{h1} and \ref{h3}
	imply \eqref{eq:lyap-sigma-bound} after multiplication by $\sqrt{d_2}$.
	\ref{h1} gives evenness and uniform positive
	upper and lower bounds for $H_t$. Its uniform Lipschitz bound follows from
	\[
	G_1^{-1}-G_2^{-1}=G_1^{-1}(G_2-G_1)G_2^{-1},
	\qquad G_i=g_2(t,x,y_i),
	\]
	the uniform bounds for $G_i^{-1}$, and the Lipschitz continuity of
	the determinant on bounded sets of matrices.
	
	We use the following reflection--synchronous coupling. For $y_1\ne y_2$,
	set
	\[
	D=y_1-y_2,\qquad r=|D|,\qquad e=\frac{D}{r},\qquad a=\frac14,
	\qquad \varphi_e(z)=z-2\langle e,z\rangle e.
	\]
	Writing $H_i(z)=H_t(y_i,z)$, define
	\begin{align*}
		m(z)&=H_1(z)\wedge H_2(z),\\
		\widetilde H(z)&=H_1(z)\wedge H_2(z)
		\wedge H_1(\varphi_e(z))\wedge H_2(\varphi_e(z)),\\
		R(z)&=\widetilde H(z)\mathbf1_{\{|z|\le ar\}}.
	\end{align*}
	For Borel $A\subset\mathbb R^{2d_2}\setminus\{0\}$, let
	\begin{align}
		\mu_{t,x}(y_1,y_2,A)
		={}&\int_{\mathbb R^{d_2}\setminus\{0\}}
		\Big[\tfrac12R(z)\mathbf1_A(z,\varphi_e(z))
		+\tfrac12R(z)\mathbf1_A(\varphi_e(z),z)\notag\\
		&\qquad +(m(z)-R(z))\mathbf1_A(z,z)\notag\\
		&\qquad +(H_1(z)-m(z))\mathbf1_A(z,0)
		+(H_2(z)-m(z))\mathbf1_A(0,z)\Big]\nu_2(\mathrm dz).
		\label{eq:lyap-coupling-measure}
	\end{align}
	On the diagonal we retain only common jumps, namely
	$\mu_{t,x}(y,y,A)=\int\mathbf1_A(z,z)H_t(y,z)\nu_2(\mathrm dz)$.
	The Brownian parts are coupled synchronously. Put
	\[
	\Sigma_t(y_1,y_2)=
	\begin{pmatrix}\sigma_t(y_1)\\\sigma_t(y_2)\end{pmatrix}
	\]
	and define
	\begin{align*}
		\Delta_{z_1,z_2}F(y_1,y_2)
		={}&F(y_1+z_1,y_2+z_2)-F(y_1,y_2)\\
		&-\mathbf1_{\{|z_1|\leq1\}}\langle\nabla_{y_1}F(y_1,y_2),z_1\rangle\\
		&-\mathbf1_{\{|z_2|\leq1\}}\langle\nabla_{y_2}F(y_1,y_2),z_2\rangle.
	\end{align*}
	The resulting coupling operator is
	\begin{align}
		\widetilde{\mathcal L}_{t,x}F(y_1,y_2)
		={}&\langle b_t(y_1),\nabla_{y_1}F\rangle
		+\langle b_t(y_2),\nabla_{y_2}F\rangle\notag\\
		&+\frac12\operatorname{tr}
		\big(\Sigma_t\Sigma_t^\top D^2F\big)\notag\\
		&+\int_{\mathbb R^{2d_2}\setminus\{0\}}
		\Delta_{z_1,z_2}F(y_1,y_2)
		\mu_{t,x}(y_1,y_2,\mathrm dz_1,\mathrm dz_2).
		\label{eq:lyap-coupling-generator}
	\end{align}
	All densities in \eqref{eq:lyap-coupling-measure} are nonnegative.
	Moreover, $R(\varphi_e(z))=R(z)$ and orthogonal reflection preserves
	$\nu_2$. The two half-reflection terms therefore contribute total density
	$R$ to either marginal. The marginal densities are consequently
	$R+(m-R)+(H_i-m)=H_i$, so \eqref{eq:lyap-coupling-generator}
	is a coupling operator for the normalized frozen generator.
	
	\begin{lemma}[Lyapunov estimate]\label{lem:lyapunov-corrected}
		Suppose that \ref{h1}--\ref{h4} hold. Let
		$\widetilde{\mathcal L}_{t,x}$ be given by
		\eqref{eq:lyap-coupling-generator}, with the isotropic stable measure
		$\nu_2$ above. Set $\omega_{d_2}=|\mathbb S^{d_2-1}|$ and
		\begin{align*}
			K_\nu&=c_{\alpha_2,d_2}\omega_{d_2}
			\left(\frac{2^{\alpha_2-2}}{2-\alpha_2}
			+\frac{2^{\alpha_2}}{\alpha_2-1}\right),\\
			A_\sigma&=\frac12L_\sigma^2,
			\qquad B_\sigma=2M_\sigma^2,
			\qquad A_J=K_\nu L_H,
			\qquad B_J=K_\nu h_1,\\
			D_{\mathrm{ref}}&=
			\frac{2h_0c_{\alpha_2,d_2}\omega_{d_2}a^{2-\alpha_2}}
			{d_2(2-\alpha_2)}.
		\end{align*}
		Choose, in the following order,
		\begin{align}
			R_0&=\max\left\{1,
			\left(\frac{4B_\sigma}{\alpha}\right)^{1/2},
			\left(\frac{4B_J}{\alpha}\right)^{1/\alpha_2}\right\},\notag\\
			K&=A_\sigma R_0^{\alpha_2-1}+A_J,\notag\\
			c&=1+\frac{2K}{D_{\mathrm{ref}}},
			\qquad \ell=\frac12e^{-cR_0},
			\qquad \lambda=\frac{\alpha\ell}{2}>0.
			\label{eq:lyap-parameter-choice}
		\end{align}
		Define
		\begin{equation}\label{eq:lyap-distance-function}
			\psi(r)=\ell r+\frac{1-\ell}{c}(1-e^{-cr}),\qquad r\geq0.
		\end{equation}
		Then $\psi$ is increasing and concave, and
		\begin{equation}\label{eq:lyap-linear-comparison}
			\ell r\leq\psi(r)\leq r,\qquad r\geq0.
		\end{equation}
		Furthermore, for $F(y_1,y_2)=\psi(|y_1-y_2|)$ and every $y_1\ne y_2$,
		\begin{equation}\label{eq:lyap-final-estimate}
			\widetilde{\mathcal L}_{t,x}F(y_1,y_2)
			\leq-\frac\alpha2|y_1-y_2|\psi'(|y_1-y_2|)
			\leq-\lambda\psi(|y_1-y_2|).
		\end{equation}
		The constants $c,\ell,\lambda$ are independent of $t,x,y_1,y_2$.
	\end{lemma}
	
	\begin{proof}
		Fix $t,x$ and $y_1\ne y_2$, and retain $D=re$ as above.
		Direct differentiation yields
		\begin{equation}\label{eq:lyap-derivatives}
			\psi'(r)=\ell+(1-\ell)e^{-cr},\qquad
			\psi''(r)=-c(1-\ell)e^{-cr}.
		\end{equation}
		Thus $\ell\leq\psi'(r)\leq1$ and $\psi''(r)<0$.
		Together with $\psi(0)=0$, these inequalities imply
		\eqref{eq:lyap-linear-comparison}.
		The generator applied to $F$ is well defined off the diagonal:
		$F$ is smooth locally there, giving a quadratic small-jump remainder,
		and its linear growth is compatible with
		$\int_{|z|>1}|z|\nu_2(\mathrm dz)<\infty$.
		
		Decompose the generator as
		\[
		\widetilde{\mathcal L}_{t,x}F
		=I_b+I_\sigma+I_{\mathrm{ref}}+I_{\mathrm{rem}}.
		\]
		The common-displacement jumps contribute zero, since
		$F(y_1+z,y_2+z)=F(y_1,y_2)$ and
		$\nabla_{y_1}F+\nabla_{y_2}F=0$.
		
		\medskip\noindent\textit{Drift and Brownian contributions.}
		By \eqref{eq:lyap-drift-bound},
		\begin{equation}\label{eq:lyap-Ib}
			I_b=\psi'(r)\frac{\langle b_t(y_1)-b_t(y_2),D\rangle}{r}
			\leq-\alpha r\psi'(r).
		\end{equation}
		Writing $\Delta\sigma=\sigma_t(y_1)-\sigma_t(y_2)$, synchronous
		Brownian coupling gives
		\begin{align}
			I_\sigma
			={}&\frac12\psi''(r)|\Delta\sigma^\top e|^2
			+\frac{\psi'(r)}{2r}
			\left(\|\Delta\sigma\|_{\mathrm{HS}}^2
			-|\Delta\sigma^\top e|^2\right)\notag\\
			\leq{}&\min\left\{A_\sigma r,\frac{B_\sigma}{r}\right\}\psi'(r).
			\label{eq:lyap-Isigma}
		\end{align}
		Here we used concavity and the two bounds
		$\|\Delta\sigma\|_{\mathrm{HS}}\leq L_\sigma r$ and
		$\|\Delta\sigma\|_{\mathrm{HS}}\leq2M_\sigma$.
		
		\medskip\noindent\textit{Unilateral jump contributions.}
		Let $a_i(z)=H_i(z)-m(z)$. These densities are even and satisfy
		$a_1(z)+a_2(z)=|H_1(z)-H_2(z)|$.
		Changing $z$ to $-z$ in the second unilateral integral gives
		\begin{align}
			I_{\mathrm{rem}}
			=\int_{\mathbb R^{d_2}\setminus\{0\}}
			\big[\psi(|D+z|)-\psi(r)-\psi'(r)\langle e,z\rangle\big]
			|H_1(z)-H_2(z)|\nu_2(\mathrm dz).
			\label{eq:lyap-rem-identity}
		\end{align}
		In this identity the linear compensation has also been included for
		$|z|>1$. This is legitimate because the density is even and the
		large-jump first moment is finite. Small-jump differences and their
		linear compensation are kept together throughout.
		
		Put $N(D,z)=|D+z|-r-\langle e,z\rangle$.
		Convexity of the Euclidean norm gives $N(D,z)\geq0$.
		Taylor's formula, using $|D+sz|\geq r/2$ for $s\in[0,1]$
		and $|z|\leq r/2$, gives $N(D,z)\leq |z|^2/r$ on this set.
		The triangle inequality also gives $N(D,z)\leq 2|z|$ for all $z$.
		Consequently,
		\begin{align}
			\int N(D,z)\nu_2(\mathrm dz)
			&\leq\frac1r\int_{|z|\le r/2}|z|^2\nu_2(\mathrm dz)
			+2\int_{|z|>r/2}|z|\nu_2(\mathrm dz)\notag\\
			&=K_\nu r^{1-\alpha_2}.
			\label{eq:lyap-norm-remainder}
		\end{align}
		By concavity of $\psi$,
		\[
		\psi(|D+z|)-\psi(r)-\psi'(r)\langle e,z\rangle
		\leq\psi'(r)N(D,z).
		\]
		Combining \eqref{eq:lyap-H-bound}, \eqref{eq:lyap-H-Lipschitz},
		\eqref{eq:lyap-rem-identity}, and \eqref{eq:lyap-norm-remainder}, we obtain
		\begin{equation}\label{eq:lyap-Irem}
			I_{\mathrm{rem}}
			\leq K_\nu\min\{L_Hr,h_1\}r^{1-\alpha_2}\psi'(r)
			=\min\{A_Jr^{2-\alpha_2},B_Jr^{1-\alpha_2}\}\psi'(r).
		\end{equation}
		
		\medskip\noindent\textit{Reflection contribution.}
		For $|z|\leq ar$ with $a=1/4$, we have $2|\langle e,z\rangle|\leq r/2$.
		The distances after the two reflected jumps are therefore
		$r+2\langle e,z\rangle$ and $r-2\langle e,z\rangle$.
		The two linear compensation terms cancel, and hence
		\begin{align}
			I_{\mathrm{ref}}=\frac12\int_{|z|\leq ar}
			\big[\psi(r+2\langle e,z\rangle)
			+\psi(r-2\langle e,z\rangle)-2\psi(r)\big]
			\widetilde H(z)\nu_2(\mathrm dz).
			\label{eq:lyap-ref-identity}
		\end{align}
		For $|s|<r$, the particular choice of $\psi$ gives the exact identity
		\[
		\psi(r+s)+\psi(r-s)-2\psi(r)
		=-\frac{2(1-\ell)e^{-cr}}{c}\big(\cosh(cs)-1\big).
		\]
		Using $\cosh q-1\geq q^2/2$, $\widetilde H(z)\geq h_0$, and isotropy
		of $\nu_2$, it follows that
		\begin{equation}
			\begin{split}
				I_{\mathrm{ref}}
				&\leq-2c(1-\ell)e^{-cr}h_0
				\int_{|z|\leq ar}\langle e,z\rangle^2\nu_2(\mathrm dz)\\
				&=-\frac{2h_0c_{\alpha_2,d_2}\omega_{d_2}a^{2-\alpha_2}}
				{d_2(2-\alpha_2)}c(1-\ell)e^{-cr}r^{2-\alpha_2}\\
				&=-D_{\mathrm{ref}}c(1-\ell)e^{-cr}r^{2-\alpha_2}\leq0.
				\label{eq:lyap-Iref}
			\end{split}
		\end{equation}
		Only the isotropic reference measure $\nu_2$ is used here;
		rotational invariance of $\widetilde H(z)\nu_2(\mathrm dz)$ is not required.
		
		\medskip\noindent\textit{Small distances: $0<r\leq R_0$.}
		Since $\ell=\tfrac12e^{-cR_0}\leq\tfrac12$,
		\[
		(1-\ell)e^{-cr}\geq\tfrac12e^{-cR_0}=\ell,
		\qquad
		(1-\ell)e^{-cr}\geq\tfrac12\psi'(r).
		\]
		Moreover, $r\leq R_0^{\alpha_2-1}r^{2-\alpha_2}$.
		Combining \eqref{eq:lyap-Ib}, \eqref{eq:lyap-Isigma},
		\eqref{eq:lyap-Irem}, and \eqref{eq:lyap-Iref}, we deduce
		\begin{align*}
			\widetilde{\mathcal L}_{t,x}F
			&\leq-\alpha r\psi'(r)
			+\left(A_\sigma R_0^{\alpha_2-1}+A_J
			-\frac12D_{\mathrm{ref}}c\right)
			r^{2-\alpha_2}\psi'(r)\\
			&=-\alpha r\psi'(r)
			+\left(K-\frac12D_{\mathrm{ref}}c\right)
			r^{2-\alpha_2}\psi'(r)\\
			&\leq-\alpha r\psi'(r),
		\end{align*}
		where the last inequality follows from
		\[
		K-\frac12D_{\mathrm{ref}}c=-\frac12D_{\mathrm{ref}}<0.
		\]
		
		\medskip\noindent\textit{Large distances: $r\geq R_0$.}
		The choice of $R_0$ ensures that
		\[
		\frac{B_\sigma}{r}\leq\frac\alpha4r,
		\qquad B_Jr^{1-\alpha_2}\leq\frac\alpha4r.
		\]
		Using \eqref{eq:lyap-Iref} and the large-distance bounds in
		\eqref{eq:lyap-Isigma} and \eqref{eq:lyap-Irem}, we obtain
		\[
		\widetilde{\mathcal L}_{t,x}F
		\leq\left(-\alpha r+\frac{B_\sigma}{r}
		+B_Jr^{1-\alpha_2}\right)\psi'(r)
		\leq-\frac\alpha2r\psi'(r).
		\]
		Combining both distance regimes gives the first inequality in
		\eqref{eq:lyap-final-estimate}. Finally,
		$r\psi'(r)\geq\ell r\geq\ell\psi(r)$ yields
		\[
		\widetilde{\mathcal L}_{t,x}F
		\leq-\frac\alpha2r\psi'(r)
		\leq-\frac{\alpha\ell}{2}\psi(r)=-\lambda\psi(r).
		\]
		The uniformity of the constants follows from
		\eqref{eq:lyap-drift-bound}--\eqref{eq:lyap-H-Lipschitz}.
	\end{proof}
	
	\begin{lemma}\label{lemma4.3} Suppose \ref{h1}-\ref{h4} hold, for any fixed $x\in\mathbb{R}^{d_1}$, $\forall t\geq 0$, $y_1,y_2\in\mathbb{R}^{d_2}$, we have,
		\begin{align}\label{ex}
			\mathcal W_1\!\left(
			P_{s,t}^{x,*}\delta_{y_1},P_{s,t}^{x,*}\delta_{y_2}
			\right) \leq C\cdot e^{-\lambda(t-s)}|y_1 - y_2|.
		\end{align}
		Here $P_{s,t}^{x,*}$ denotes the action of the transition
		family on probability measures. More generally, for every
		$\nu_1,\nu_2\in\mathcal P_1(\mathbb R^{d_2})$,
		\begin{equation}
			\mathcal W_1\!\left(
			P_{s,t}^{x,*}\nu_1,P_{s,t}^{x,*}\nu_2
			\right)
			\leq \ell^{-1}e^{-\lambda(t-s)}
			\mathcal W_1(\nu_1,\nu_2),
			\label{eq:frozen-measure-w1-revised}
		\end{equation}
		where $\mathcal P_1(\mathbb R^{d_2})$ is the space of
		probability measures with finite first moment.
		Consequently, for every globally Lipschitz function $f$,
		\begin{equation}
			\operatorname{Lip}(P_{s,t}^{x}f)
			\leq \ell^{-1}e^{-\lambda(t-s)}
			\operatorname{Lip}(f).
			\label{eq:frozen-semigroup-lip-revised}
		\end{equation}
		All constants are independent of $s,t$ and $x$.
	\end{lemma}
	\begin{proof}
		We apply the coupling process $(Y_t', Z_t')$ constructed in Section 4.1. Taking the similar procedure in [\citenum{Wang2018}, p.1012, Proof of Theorem 1.2], we have
		\[
		\mathbb{E}\psi(|Y_t^{s,x,y_1} - Y_t^{s,x,y_2}|) \leq \psi(|y_1 - y_2|)e^{-\lambda(t-s)},
		\]
		Using the global comparison $\ell r\leq\psi(r)\leq r$, we conclude that
		\begin{align}\label{4.7}
			\mathcal W_1\!\left(
			P_{s,t}^{x,*}\delta_{y_1},P_{s,t}^{x,*}\delta_{y_2}
			\right)&\leq \mathbb{E}|Y_t^{s,x,y_1} - Y_t^{s,x,y_2}| \leq \frac{1}{l}\mathbb{E}\psi(|Y_t^{s,x,y_1} - Y_t^{s,x,y_2}|)\notag\\& \leq\frac{1}{l} e^{-\lambda(t-s)}\psi(|y_1 - y_2|)
			\leq \frac{1}{l} e^{-\lambda(t-s)}|y_1 - y_2|, 
		\end{align}
		the second inequality follows from definition of $\psi$ in Lemma \ref{lem:lyapunov-corrected}.
		For every globally Lipschitz $f$,
		\[
		|P_{s,t}^xf(y_1)-P_{s,t}^xf(y_2)|
		\leq \operatorname{Lip}(f)\,
		\mathcal W_1\!\left(
		P_{s,t}^{x,*}\delta_{y_1},
		P_{s,t}^{x,*}\delta_{y_2}\right),
		\]
		which proves \eqref{eq:frozen-semigroup-lip-revised}.
		Applying the Kantorovich--Rubinstein duality and this
		Lipschitz estimate gives
		\begin{align*}
			&\mathcal W_1(P_{s,t}^{x,*}\nu_1,P_{s,t}^{x,*}\nu_2)=
			\sup_{\operatorname{Lip}(f)\leq1}
			\left|\int P_{s,t}^xf\,d\nu_1-\int P_{s,t}^xf\,d\nu_2\right|\leq
			\ell^{-1}e^{-\lambda(t-s)}\mathcal W_1(\nu_1,\nu_2).
		\end{align*}
		This proves \eqref{eq:frozen-measure-w1-revised}.
	\end{proof}
	In order to prove the existence and the uniqueness of evolution system of measures for the inhomogeneous SDE (\ref{frozen}), we primarily adhere to the corresponding concept in \cite{DaPratoRockner2006}. For this aim, we need the following lemma.
	\begin{lemma}[Evolution system of measures and pullback convergence]
		\label{lem:frozen-evolution-revised}
		Under the assumptions of Lemma~\ref{lemma4.3},
		for each fixed $x\in\mathbb R^{d_1}$ there exists an
		evolution system of probability measures
		$\{\mu_t^x\}_{t\in\mathbb R}\subset
		\mathcal P_1(\mathbb R^{d_2})$ satisfying
		\begin{equation}
			P_{s,t}^{x,*}\mu_s^x=\mu_t^x,\qquad s\leq t.
			\label{eq:evolution-property-revised}
		\end{equation}
		For every $p\in[1,\alpha_2)$,
		\begin{equation}
			\sup_{x\in\mathbb R^{d_1}}\sup_{t\in\mathbb R}
			\int_{\mathbb R^{d_2}}|z|^p\,\mu_t^x(dz)\leq C_p.
			\label{eq:evolution-moments-revised}
		\end{equation}
		Moreover, for every $y\in\mathbb R^{d_2}$ and $s\leq t$,
		\begin{equation}
			\mathcal W_1(P_{s,t}^{x,*}\delta_y,\mu_t^x)
			\leq C(1+|y|)e^{-\lambda(t-s)}.
			\label{eq:pullback-w1-revised}
		\end{equation}
		In particular,
		\[
		\mu_t^x
		=\lim_{s\to-\infty}P_{s,t}^{x,*}\delta_y
		\quad\text{in }\mathcal W_1,
		\]
		and the limit is independent of the fixed initial point $y$.
		For each fixed $x$, this evolution system is unique among
		all evolution systems $\{\nu_t^x\}_{t\in\mathbb R}$ satisfying
		\[
		\sup_{t\in\mathbb R}\int_{\mathbb R^{d_2}}|z|\,\nu_t^x(dz)
		<\infty.
		\]
		Consequently, for every globally Lipschitz function $f$,
		\begin{equation}
			\left|P_{s,t}^xf(y)-\int f(z)\,\mu_t^x(dz)\right|
			\leq C\operatorname{Lip}(f)(1+|y|)e^{-\lambda(t-s)}.
			\label{eq:pullback-observable-revised}
		\end{equation}
		The constants above are independent of $s,t$ and $x$.
	\end{lemma}
	
	\begin{proof}
		According to \eqref{frozen1}, we have
		\begin{equation}
			\sup_{x\in\mathbb R^{d_1}}\sup_{u\le v}
			\int |z|^p\,P_{u,v}^{x,*}\delta_0(dz)\leq M_p.
			\label{eq:frozen-zero-moment-revised}
		\end{equation}
		
		Fix $t\in\mathbb R$ and $x\in\mathbb R^{d_1}$, and put
		$\nu_s=P_{s,t}^{x,*}\delta_0$ for $s<t$.
		For $r<s<t$, the evolution property of the transition family
		and \eqref{eq:frozen-measure-w1-revised} give
		\begin{align*}
			\mathcal W_1(\nu_r,\nu_s)
			&=\mathcal W_1\!\left(
			P_{s,t}^{x,*}(P_{r,s}^{x,*}\delta_0),
			P_{s,t}^{x,*}\delta_0\right)\\
			&\leq\ell^{-1}e^{-\lambda(t-s)}
			\mathcal W_1(P_{r,s}^{x,*}\delta_0,\delta_0)\\
			&=\ell^{-1}e^{-\lambda(t-s)}
			\int |z|\,P_{r,s}^{x,*}\delta_0(dz)\\
			&\leq Ce^{-\lambda(t-s)}.
		\end{align*}
		Thus $\{\nu_s\}_{s<t}$ is a Cauchy family in
		$\mathcal W_1$ as $s\to-\infty$.
		Completeness of
		$(\mathcal P_1(\mathbb R^{d_2}),\mathcal W_1)$ gives a limit
		$\mu_t^x$. By lower semicontinuity and
		\eqref{eq:frozen-zero-moment-revised}, this family satisfies
		\eqref{eq:evolution-moments-revised}.
		
		For $u<t$, the map $P_{u,t}^{x,*}$ is Lipschitz continuous
		in $\mathcal W_1$ by
		\eqref{eq:frozen-measure-w1-revised}. Therefore,
		\begin{align*}
			P_{u,t}^{x,*}\mu_u^x
			=\lim_{s\to-\infty}
			P_{u,t}^{x,*}(P_{s,u}^{x,*}\delta_0)=\lim_{s\to-\infty}P_{s,t}^{x,*}\delta_0
			=\mu_t^x.
		\end{align*}
		This proves \eqref{eq:evolution-property-revised}.
		
		Using this evolution property and
		\eqref{eq:frozen-measure-w1-revised}, we obtain
		\begin{align*}
			\mathcal W_1(P_{s,t}^{x,*}\delta_y,\mu_t^x)
			&=\mathcal W_1(
			P_{s,t}^{x,*}\delta_y,P_{s,t}^{x,*}\mu_s^x)\\
			&\leq\ell^{-1}e^{-\lambda(t-s)}
			\mathcal W_1(\delta_y,\mu_s^x)\\
			&\leq\ell^{-1}e^{-\lambda(t-s)}
			\left(|y|+\int|z|\,\mu_s^x(dz)\right)\\
			&\leq C(1+|y|)e^{-\lambda(t-s)}.
		\end{align*}
		This proves \eqref{eq:pullback-w1-revised} and the
		independence of the pullback limit from $y$.
		
		Finally, let $\{\widehat\mu_t^x\}_{t\in\mathbb R}$ be another
		evolution system with
		$\widehat M_x=\sup_t\int|z|\,\widehat\mu_t^x(dz)<\infty$.
		For every $s<t$,
		\begin{align*}
			\mathcal W_1(\mu_t^x,\widehat\mu_t^x)
			\leq\ell^{-1}e^{-\lambda(t-s)}
			\mathcal W_1(\mu_s^x,\widehat\mu_s^x)\leq\ell^{-1}e^{-\lambda(t-s)}
			(M_1+\widehat M_x).
		\end{align*}
		Letting $s\to-\infty$ proves uniqueness in the stated class.
		Estimate \eqref{eq:pullback-observable-revised} follows
		from \eqref{eq:pullback-w1-revised} and
		the Kantorovich--Rubinstein inequality.
	\end{proof}
	\subsection{Regularity of the Poisson Equation}
	Inspired by \cite{RocknerXie2021, SunXieXie2022, Yin2025}, we consider the nonlocal Poisson equation as follows,
	\begin{align}\label{possionequation}
		\partial_s \Phi(s,x,y)
		+\mathcal{L}_2(s,x,y) \Phi(s,x,y) + b_1(x,y) - \bar{b}_1(s,x) = 0,
	\end{align}
	here we define 
	\begin{align}\label{avco}
		\bar{b}_1(s, x) = \int_{\mathbb{R}^m} b_1(x, y) \, \mu_s^x(dy),
	\end{align}
	and we have the following regularity estimates.
	
	\begin{theorem}\label{throrem4.5}
		Let \ref{h1}-\ref{h4} hold, $\forall x \in \mathbb{R}^{d_1}$, $y \in \mathbb{R}^{d_2}$, we define
		\[
		\Phi(s,x,y) = \int_s^\infty \big[ \mathbb{E} b_1(x,Y_r^{s,x,y}) - \bar{b}_1(r,x) \big] dr, 
		\]
		then $\Phi(s,x,y)$ is a solution of (\ref{possionequation}). And there exists $C > 0$ s.t.,
		\begin{align}\label{y}
			\sup_{x \in \mathbb{R}^{d_1}} |\Phi(s,x,y)| \leq C_T(1+|y|). 
		\end{align}
		Additionally, we have the Lipschitz estimate,
		\begin{align}\label{lipschitz estimate}
			\operatorname{Lip}_{y}\Phi(s,x,\cdot)
			\leq
			C.
		\end{align}
	\end{theorem}
	\begin{proof}
		Our proof refers to [\citenum{SunXieXie2022}, Proposition 3.3] and [\citenum{Yin2025}, Theorem 5.2]. We first establish the integral representation and the
		growth and Lipschitz estimates for $\Phi$.
		The classical solvability of the non-autonomous Poisson
		equation will be justified after Corollary~4.11, once the
		required spatial regularity has been established., we can deduce that $\Phi(s,x,y)$ is a solution of (\ref{possionequation}). 
		Let 
		\begin{align*}
			H(r,x,y)=b_1(x,y)-\bar{b}_1(r,x),
		\end{align*}
		Obviously, we have 
		\begin{align*}
			\mathbb{E}H(s,x,\eta_s^x)=\int_{\mathbb{R}^{d_2}} H(s,x,y)\,\mu_s^x(dy)=0,\quad s\in\mathbb{R}, x\in\mathbb{R}^{d_1},
		\end{align*}
		for $r \geqslant s$,
		\[
		\mathbb{E}H(r,x,Y_r^{s,x,y}) = \mathbb{E}H(r,x,Y_r^{s,x,y}) - \mathbb{E}H(r,x,\eta_r^x).
		\]
		Note that, by Lipschitz condition \ref{h3}, we have
		\[
		|H(r,x,y_1)-H(r,x,y_2)| \leqslant C|y_1-y_2|.
		\]
		Then, according to (\ref{eq:pullback-w1-revised}),
		\begin{align*}
			|\Phi(s,x,y)| &\leqslant \int_{s}^{+\infty} \big|\mathbb{E}H(r,x,Y_r^{s,x,y}) - \mathbb{E}H(r,x,\eta_r^x)\big|\,dr \\
			&\leqslant C(1+|y|)\int_{s}^{+\infty} e^{-\lambda(r-s)} dr = C(1+|y|),
		\end{align*}
		which proves the first inequality \eqref{y}.
		
		Next we deal with gradient estimate (\ref{lipschitz estimate}). 
		\[
		\begin{aligned}
			\left|
			\Phi(s,x,y_1)-\Phi(s,x,y_2)
			\right|
			&\leq
			\int_s^\infty
			\left|
			P_{s,r}^{x} b_1(x,\cdot)(y_1)
			-
			P_{s,r}^{x} b_1(x,\cdot)(y_2)
			\right|\,dr
			\\
			&\leq
			C|y_1-y_2|
			\int_s^\infty
			e^{-\lambda(r-s)}.
		\end{aligned}
		\]
		So, we have
		\[
		\operatorname{Lip}_{y}\Phi(s,x,\cdot)
		\leq
		C. \qedhere
		\]
	\end{proof}
	We define the averaging coefficient $\bar{b}(s, x)$ as follows,
	\begin{align}
		\bar{b}(s, x) = \int_{\mathbb{R}^m} b(x, y) \, \mu_s^x(dy),
	\end{align}
	
	Next we prove the Lipschitz property of the averaged coefficient with respect to $x$. Before proving the main result, we state several lemmas as preliminaries.
	
	For each fixed $x\in\mathbb R^{d_1}$, let
	$\{P_{s,t}^x\}_{s\le t}$ denote the evolution family associated with
	the frozen equation. We decompose its generator as
	\[
	\mathcal L_t^x=\mathcal A_t^x+\mathcal J_t^x,
	\]
	where
	\[
	\mathcal A_t^x\varphi(y)
	=
	\langle b_2(t,x,y),\nabla\varphi(y)\rangle
	+
	\frac12
	\operatorname{Tr}\big(
	a_2(t,x,y)D^2\varphi(y)
	\big),
	\]
	with
	\[
	a_2(t,x,y)
	=
	\sigma_2(t,x,y)\sigma_2(t,x,y)^*,
	\]
	and $\mathcal J_t^x$ denotes the nonlocal part generated by the
	multiplicative $\alpha_2$-stable noise.
	
	Let $\{Q_{s,t}^x\}_{s\le t}$ be the evolution family associated with
	the reference diffusion operator $\mathcal A_t^x$.
	\begin{lemma}[Estimates for the reference diffusion]\cite{MenozziPesceZhang2021}
		\label{lem:reference-diffusion}
		Under Assumption \ref{h1}-\ref{h4}, there exists a constant
		$C>0$, independent of $s\in\mathbb R$ and
		$x\in\mathbb R^{d_1}$, such that, for $0<t-s\leq1$,
		\begin{align}\label{41}
			\|\nabla Q_{s,t}^x h\|_\infty
			\leq
			C(t-s)^{-1/2}\|h\|_\infty,
		\end{align}
		and
		\begin{align}\label{42}
			\|D^2 Q_{s,t}^x h\|_\infty
			\leq
			C(t-s)^{-1}\|h\|_\infty.
		\end{align}
		Moreover, if $h$ is Lipschitz continuous, then
		\begin{align}\label{43}
			\|\nabla Q_{s,t}^x h\|_\infty
			\leq
			C\operatorname{Lip}(h),
		\end{align}
		and
		\begin{align}\label{44}
			\|D^2 Q_{s,t}^x h\|_\infty
			\leq
			C(t-s)^{-1/2}\operatorname{Lip}(h).
		\end{align}
		More generally, for every $\theta\in(0,1)$,
		\begin{align}\label{45}
			[\nabla Q_{s,t}^x h]_{C^\theta}
			\leq
			C(t-s)^{-(1+\theta)/2}\|h\|_\infty,
		\end{align}
		and, for Lipschitz $h$,
		\begin{align}\label{46}
			[\nabla Q_{s,t}^x h]_{C^\theta}
			\leq
			C(t-s)^{-\theta/2}\operatorname{Lip}(h).
		\end{align}
	\end{lemma}
	\begin{proof}
		The estimates (\ref{41})-(\ref{42}) follow from the first and second order
		transition density estimates for uniformly non-degenerate
		time-inhomogeneous diffusions with linearly growing drift.
		
		We briefly explain why the estimates are uniform in the frozen
		parameter $x$. By Assumption~\ref{h4}, the
		ellipticity constants, the linear-growth constant of $b_2$, and the
		spatial regularity bounds of $b_2$ and $\sigma_2$ are all uniform in
		$(t,x)$. Hence the constants in the corresponding transition-density
		estimates can be chosen independently of $x$ and of the initial time.
		
		Let $q^x(s,y;t,z)$ denote the transition density of $Q_{s,t}^x$.
		The cancellation identities
		\[
		\int_{\mathbb R^{d_2}}
		\nabla_y q^x(s,y;t,z)\,dz=0,
		\qquad
		\int_{\mathbb R^{d_2}}
		D_y^2 q^x(s,y;t,z)\,dz=0
		\]
		and the Gaussian-type derivative bounds yield, for Lipschitz $h$,
		\[
		\|\nabla Q_{s,t}^x h\|_\infty
		\leq C\operatorname{Lip}(h),
		\]
		and
		\[
		\|D^2Q_{s,t}^x h\|_\infty
		\leq
		C(t-s)^{-1/2}\operatorname{Lip}(h).
		\]
		Finally, (\ref{45})-(\ref{46}) follow by interpolation between the gradient
		and Hessian estimates.
	\end{proof}
	Fix
	\begin{align}
		\alpha_2-1<\beta<1,
		\qquad
		\gamma:=1+\beta-\alpha_2.
	\end{align}
	Then $\gamma\in(0,2-\alpha_2)$.
	\begin{lemma}[H\"older estimate for the jump operator]
		\label{lem:jump-holder}
		There exists a constant $C>0$, uniform in $(t,x)$, such that
		\begin{align}\label{423}
			[\mathcal J_t^x u]_{C^\gamma}
			\leq
			C\left(
			\|\nabla u\|_\infty
			+
			[\nabla u]_{C^\beta}
			\right)
		\end{align}
		for every $u\in C_b^{1+\beta}(\mathbb R^{d_2})$.
		Moreover,
		\begin{align}\label{424}
			\|\mathcal J_t^x u\|_\infty
			\leq
			C\left(
			\|\nabla u\|_\infty
			+
			[\nabla u]_{C^\beta}
			\right).
		\end{align}
	\end{lemma}
	\begin{proof}
		By the uniform invertibility of $g_2$, the nonlocal operator can be
		written in stable-like form
		\[
		\mathcal J_t^x u(y)
		=
		\int_{\mathbb R^{d_2}}
		\delta u(y;z)
		\frac{\kappa^x(t,y,z)}
		{|z|^{d_2+\alpha_2}}\,dz,
		\]
		where
		\[
		\delta u(y;z)
		=
		u(y+z)-u(y)
		-\mathbf 1_{\{|z|\leq1\}}
		\langle\nabla u(y),z\rangle.
		\]
		The assumptions on $g_2$ imply that $\kappa^x$ is uniformly bounded
		and uniformly H\"older continuous in $y$; in particular,
		\begin{align}
			|\kappa^x(t,y_1,z)-\kappa^x(t,y_2,z)|
			\leq C|y_1-y_2|^\gamma.
		\end{align}
		
		For $|z|\leq1$, the $C^{1+\beta}$ regularity of $u$ gives
		\begin{align}
			|\delta u(y;z)|
			\leq
			C[\nabla u]_{C^\beta}|z|^{1+\beta}.
		\end{align}
		Since
		\[
		1+\beta-\alpha_2=\gamma>0,
		\]
		the small-jump integral is finite.
		
		Let $h=|y_1-y_2|\leq1$. Splitting the jump integral into
		\[
		|z|\leq h,\qquad h<|z|\leq1,\qquad |z|>1,
		\]
		we obtain
		\[
		\int_{|z|\leq h}
		|\delta u(y_1;z)-\delta u(y_2;z)|
		\frac{dz}{|z|^{d_2+\alpha_2}}
		\leq
		C[\nabla u]_{C^\beta}h^{1+\beta-\alpha_2}\leq
		C[\nabla u]_{C^\beta}h^\gamma.
		\]
		On $h<|z|\leq1$,
		\[
		|\delta u(y_1;z)-\delta u(y_2;z)|
		\leq
		C[\nabla u]_{C^\beta}h^\beta|z|,
		\]
		and therefore
		\begin{align}
			\begin{aligned}
				&\int_{h<|z|\leq1}
				|\delta u(y_1;z)-\delta u(y_2;z)|
				\frac{dz}{|z|^{d_2+\alpha_2}}
				\\
				&\qquad\leq
				C[\nabla u]_{C^\beta}
				h^\beta
				\int_h^1\rho^{-\alpha_2}\,d\rho
				\\
				&\qquad\leq
				C[\nabla u]_{C^\beta}
				h^{1+\beta-\alpha_2}
				=
				C[\nabla u]_{C^\beta}h^\gamma.
			\end{aligned}
		\end{align}
		
		Finally, for $|z|>1$,
		\[
		|\delta u(y_1;z)-\delta u(y_2;z)|
		\leq
		2\|\nabla u\|_\infty h,
		\]
		and since the large-jump part of the stable L\'evy measure has finite
		mass,
		\begin{align}
			\int_{|z|>1}
			|\delta u(y_1;z)-\delta u(y_2;z)|
			\frac{dz}{|z|^{d_2+\alpha_2}}
			\leq
			C\|\nabla u\|_\infty h^\gamma.
		\end{align}
		It remains to estimate the contribution arising from the spatial
		dependence of the kernel. More precisely, writing
		\[
		\begin{aligned}
			\mathcal J_t^x u(y_1)-\mathcal J_t^x u(y_2)
			={}&
			\int_{\mathbb R^{d_2}}
			\bigl[\delta u(y_1;z)-\delta u(y_2;z)\bigr]
			\frac{\kappa^x(t,y_1,z)}
			{|z|^{d_2+\alpha_2}}\,dz
			\\
			&+
			\int_{\mathbb R^{d_2}}
			\delta u(y_2;z)
			\frac{\kappa^x(t,y_1,z)-\kappa^x(t,y_2,z)}
			{|z|^{d_2+\alpha_2}}\,dz
			\\
			=:{}& I_1+I_2,
		\end{aligned}
		\]
		the preceding estimates give
		\[
		|I_1|
		\leq
		C |y_1-y_2|^\gamma
		\left(
		\|\nabla u\|_\infty+[\nabla u]_{C^\beta}
		\right).
		\]
		
		For $I_2$, by the uniform $C^\gamma$ regularity of
		$\kappa^x$ in the $y$-variable,
		\[
		|\kappa^x(t,y_1,z)-\kappa^x(t,y_2,z)|
		\leq C|y_1-y_2|^\gamma.
		\]
		Hence, setting $h=|y_1-y_2|\leq1$,
		\[
		\begin{aligned}
			|I_2|
			&\leq
			C h^\gamma
			\int_{\mathbb R^{d_2}}
			|\delta u(y_2;z)|
			\frac{dz}{|z|^{d_2+\alpha_2}}
			\\
			&\leq
			C h^\gamma [\nabla u]_{C^\beta}
			\int_{|z|\leq1}
			\frac{|z|^{1+\beta}}
			{|z|^{d_2+\alpha_2}}\,dz
			\\
			&\quad+
			C h^\gamma \|\nabla u\|_\infty
			\int_{|z|>1}
			\frac{|z|}
			{|z|^{d_2+\alpha_2}}\,dz.
		\end{aligned}
		\]
		Since $1+\beta-\alpha_2=\gamma>0$ and $\alpha_2>1$,
		both integrals are finite. Therefore,
		\[
		|I_2|
		\leq
		C h^\gamma
		\left(
		[\nabla u]_{C^\beta}
		+\|\nabla u\|_\infty
		\right).
		\]
		Combining the estimates for $I_1$ and $I_2$ yields
		\[
		[\mathcal J_t^x u]_{C^\gamma}
		\leq
		C\left(
		\|\nabla u\|_\infty
		+[\nabla u]_{C^\beta}
		\right).
		\]
		The estimate (\ref{424}) follows in the same way from (\ref{423}).
	\end{proof}
	\begin{lemma}[Fractional gradient regularization]
		\label{lem:fractional-gradient}
		For every $\beta\in(\alpha_2-1,1)$, there exists $C>0$ such that
		\begin{align}\label{431}
			[\nabla P_{s,t}^x f]_{C^\beta}
			\leq
			C(t-s)^{-\beta/2}\operatorname{Lip}(f),
			\qquad
			0<t-s\leq1,
		\end{align}
		uniformly in $s\in\mathbb R$ and $x\in\mathbb R^{d_1}$.
	\end{lemma}
	\begin{proof}
		We first prove the estimate on a sufficiently short time interval.
		The variation-of-constants formula gives
		\begin{align}\label{432}
			P_{s,t}^xf
			=
			Q_{s,t}^xf
			+
			\int_s^t
			Q_{s,r}^x\mathcal J_r^xP_{r,t}^xf\,dr.
		\end{align}
		Using Lemma~\ref{lem:reference-diffusion}, we obtain
		\[
		\begin{aligned}
			[\nabla P_{s,t}^xf]_{C^\beta}
			\leq{}&
			C(t-s)^{-\beta/2}\operatorname{Lip}(f)
			\\
			&+
			C\int_s^t
			(r-s)^{-(1+\beta)/2}
			\|\mathcal J_r^xP_{r,t}^xf\|_\infty\,dr.
		\end{aligned}
		\]
		
		By the fixed-parameter contraction estimate,
		\[
		\|\nabla P_{r,t}^xf\|_\infty
		\leq C\operatorname{Lip}(f).
		\]
		Hence Lemma~\ref{lem:jump-holder} yields
		\[
		\|\mathcal J_r^xP_{r,t}^xf\|_\infty
		\leq
		C\left(
		\operatorname{Lip}(f)
		+
		[\nabla P_{r,t}^xf]_{C^\beta}
		\right).
		\]
		
		Set
		\[
		F(\tau)
		:=
		\sup_{\substack{t-s=\tau\\x\in\mathbb R^{d_1}}}
		\frac{
			[\nabla P_{s,t}^xf]_{C^\beta}
		}{
			\operatorname{Lip}(f)
		}.
		\]
		Then
		\begin{align}\label{430}
			F(\tau)
			\leq
			C\tau^{-\beta/2}
			+
			C\int_0^\tau
			u^{-(1+\beta)/2}
			\bigl(1+F(\tau-u)\bigr)\,du.
		\end{align}
		
		For $T\leq1$, define
		\[
		M_T:=\sup_{0<\tau\leq T}\tau^{\beta/2}F(\tau).
		\]
		Multiplying (\ref{430}) by $\tau^{\beta/2}$ and using the beta-function
		identity gives
		\[
		M_T
		\leq
		C+CT^{1/2}
		+
		CT^{(1-\beta)/2}M_T.
		\]
		Choose $T_0>0$ sufficiently small so that
		\[
		CT_0^{(1-\beta)/2}\leq\frac12.
		\]
		Then
		\[
		M_{T_0}\leq C,
		\]
		which proves
		\[
		[\nabla P_{s,t}^xf]_{C^\beta}
		\leq
		C(t-s)^{-\beta/2}\operatorname{Lip}(f),
		\qquad
		0<t-s\leq T_0.
		\]
		
		For $T_0<t-s\leq1$, use the evolution property
		\[
		P_{s,t}^x
		=
		P_{s,s+T_0/2}^x
		P_{s+T_0/2,t}^x
		\]
		together with the fixed-parameter Lipschitz contraction estimate.
		This extends (\ref{431}) to the whole interval $0<t-s\leq1$.
	\end{proof}
	\begin{proposition}[Short-time Hessian estimate]
		\label{prop:hessian-frozen}
		Under Assumptions~\ref{h1}-
		\ref{h4}, there exists a constant $C>0$,
		independent of $s\in\mathbb R$ and
		$x\in\mathbb R^{d_1}$, such that
		\begin{align}\label{hessian2}
			\|D_y^2P_{s,t}^xf\|_\infty
			\leq
			C(t-s)^{-1/2}\operatorname{Lip}(f),
			\qquad
			0<t-s\leq1.
		\end{align}
	\end{proposition}
	\begin{proof}
		Fix
		\[
		\beta\in(\alpha_2-1,1),
		\qquad
		\gamma=1+\beta-\alpha_2.
		\]
		By Lemmas~\ref{lem:jump-holder} and
		\ref{lem:fractional-gradient},
		\[
		[\mathcal J_r^xP_{r,t}^xf]_{C^\gamma}
		\leq
		C(t-r)^{-\beta/2}\operatorname{Lip}(f).
		\]
		
		Differentiating the variation-of-constants formula (\ref{432}) twice
		with respect to the initial spatial variable gives
		\[
		D_y^2P_{s,t}^xf
		=
		D_y^2Q_{s,t}^xf
		+
		\int_s^t
		D_y^2Q_{s,r}^x
		\mathcal J_r^xP_{r,t}^xf\,dr.
		\]
		
		The first term satisfies
		\begin{align}\label{435}
			\|D_y^2Q_{s,t}^xf\|_\infty
			\leq
			C(t-s)^{-1/2}\operatorname{Lip}(f)
		\end{align}
		by Lemma~\ref{lem:reference-diffusion}.
		
		For the second term, the cancellation property of the transition
		density of $Q^x$ yields
		\[
		\|D_y^2Q_{s,r}^xh\|_\infty
		\leq
		C(r-s)^{-1+\gamma/2}[h]_{C^\gamma}.
		\]
		Consequently, by (\ref{423}),
		\begin{align}\label{434}
			\begin{aligned}
				&\left\|
				\int_s^t
				D_y^2Q_{s,r}^x
				\mathcal J_r^xP_{r,t}^xf\,dr
				\right\|_\infty
				\\
				&\quad\leq
				C\operatorname{Lip}(f)
				\int_s^t
				(r-s)^{-1+\gamma/2}
				(t-r)^{-\beta/2}\,dr.
			\end{aligned}
		\end{align}
		
		Since $\gamma>0$ and $\beta<1$, the integral is finite. By the
		beta-function identity,
		\[
		\begin{aligned}
			&\int_s^t
			(r-s)^{-1+\gamma/2}
			(t-r)^{-\beta/2}\,dr
			\\
			&\quad=
			B\left(
			\frac{\gamma}{2},
			1-\frac{\beta}{2}
			\right)
			(t-s)^{(\gamma-\beta)/2}.
		\end{aligned}
		\]
		Since
		\[
		\gamma-\beta=1-\alpha_2,
		\]
		we obtain
		\[
		\left\|
		\int_s^t
		D_y^2Q_{s,r}^x
		\mathcal J_r^xP_{r,t}^xf\,dr
		\right\|_\infty
		\leq
		C(t-s)^{-(\alpha_2-1)/2}.
		\]
		As $\alpha_2\in(1,2)$,
		\[
		\frac{\alpha_2-1}{2}<\frac12.
		\]
		Hence, for $0<t-s\leq1$,
		\[
		(t-s)^{-(\alpha_2-1)/2}
		\leq
		(t-s)^{-1/2}.
		\]
		Combining (\ref{435})-(\ref{434}) proves
		\[
		\|D_y^2P_{s,t}^xf\|_\infty
		\leq
		C(t-s)^{-1/2}\operatorname{Lip}(f). \qedhere
		\]
	\end{proof}
	\begin{corollary}
		\label{cor:gradient-Phi}
		Under Assumptions \ref{h1}-\ref{h4}, for every
		$(s,x)\in\mathbb{R}\times\mathbb{R}^{d_1}$,
		the mapping
		\[
		y\longmapsto \Phi(s,x,y)
		\]
		is continuously differentiable. Moreover,
		\[
		\sup_{s\in\mathbb R,\,
			x\in\mathbb R^{d_1},\,
			y\in\mathbb R^{d_2}}
		|\nabla_y\Phi(s,x,y)|
		<\infty.
		\]
	\end{corollary}
	\begin{proof}
		Recall that the solution of the non-autonomous Poisson equation admits
		the representation
		\[
		\Phi(s,x,y)
		=
		\int_s^\infty
		\left[
		P_{s,r}^{x}b_1(x,\cdot)(y)
		-
		\bar b_1(r,x)
		\right]\,dr .
		\]
		Since $\bar b_1(r,x)$ is independent of $y$, it does not contribute
		to the derivative with respect to $y$.
		
		By the regularization estimates established above, for every $r>s$,
		the mapping
		\[
		y\longmapsto P_{s,r}^{x}b_1(x,\cdot)(y)
		\]
		is continuously differentiable. On the other hand, the
		fixed-parameter Lipschitz contraction estimate gives
		\[
		\begin{aligned}
			\left\|
			\nabla_y P_{s,r}^{x}b_1(x,\cdot)
			\right\|_{\infty}
			&\leq
			\operatorname{Lip}
			\left(
			P_{s,r}^{x}b_1(x,\cdot)
			\right)
			\leq
			C
			e^{-\lambda(r-s)}
			\operatorname{Lip}_y
			\left(
			b_1(x,\cdot)
			\right).
		\end{aligned}
		\]
		Since $b_1$ is uniformly Lipschitz continuous in the fast variable,
		we obtain
		\[
		\left\|
		\nabla_y P_{s,r}^{x}b_1(x,\cdot)
		\right\|_{\infty}
		\leq
		C
		e^{-\lambda(r-s)}.
		\]
		Hence, by differentiation under the
		integral sign,
		\[
		\nabla_y\Phi(s,x,y)
		=
		\int_s^\infty
		\nabla_y
		P_{s,r}^{x}b_1(x,\cdot)(y)\,dr .
		\]
		Consequently,
		\[
		\begin{aligned}
			|\nabla_y\Phi(s,x,y)|
			\leq
			\int_s^\infty
			\left|
			\nabla_yP_{s,r}^{x}b_1(x,\cdot)(y)
			\right|\,dr
			\leq
			C
			\int_s^\infty
			e^{-\lambda(r-s)}dr=
			C.
		\end{aligned}
		\]
		Taking the supremum over $y$ yields
		\[
		\|\nabla_y\Phi(s,x,\cdot)\|_{\infty}
		\leq
		C.\qedhere
		\]
	\end{proof}
	\begin{corollary}
		\label{cor:hessian-Phi}
		Under Assumptions \ref{h1}-\ref{h4}, the solution $\Phi$ of the
		non-autonomous Poisson equation satisfies
		\[
		\Phi(s,x,\cdot)\in C^2(\mathbb R^{d_2}),
		\]
		and
		\[
		\sup_{s\in\mathbb R,\;x\in\mathbb R^{d_1}}
		\|D_y^2\Phi(s,x,\cdot)\|_\infty
		\leq C.
		\]
	\end{corollary}
	\begin{proof}
		Recall that
		\[
		\Phi(s,x,y)
		=
		\int_s^\infty
		\left[
		P_{s,r}^x b_1(x,\cdot)(y)
		-
		\bar b_1(r,x)
		\right]\,dr .
		\]
		Since $\bar b_1(r,x)$ is independent of $y$, formally,
		\[
		D_y^2\Phi(s,x,y)
		=
		\int_s^\infty
		D_y^2P_{s,r}^x b_1(x,\cdot)(y)\,dr .
		\]
		We justify this representation by showing that the integral on the
		right-hand side is absolutely convergent.
		
		For $0<r-s\leq1$, Proposition~\ref{prop:hessian-frozen}
		gives
		\[
		\left\|
		D_y^2P_{s,r}^x b_1(x,\cdot)
		\right\|_\infty
		\leq
		C(r-s)^{-1/2}
		\operatorname{Lip}_y(b_1).
		\]
		Hence,
		\[
		\int_s^{s+1}
		\left\|
		D_y^2P_{s,r}^x b_1(x,\cdot)
		\right\|_\infty\,dr
		\leq
		C\operatorname{Lip}_y(b_1)
		\int_0^1 \tau^{-1/2}\,d\tau
		\leq C.
		\]
		
		For $r-s>1$, by the evolution property,
		\[
		P_{s,r}^x
		=
		P_{s,s+1}^xP_{s+1,r}^x.
		\]
		Applying the unit-time Hessian regularization estimate to
		$P_{s+1,r}^x b_1(x,\cdot)$ yields
		\[
		\begin{aligned}
			\left\|
			D_y^2P_{s,r}^x b_1(x,\cdot)
			\right\|_\infty
			\leq
			C\operatorname{Lip}
			\left(
			P_{s+1,r}^x b_1(x,\cdot)
			\right)
			\leq
			C
			\exp\left\{
			-\int_{s+1}^{r}\lambda\,dq
			\right\}
			\operatorname{Lip}_y(b_1),
		\end{aligned}
		\]
		where the second inequality follows from the fixed-parameter
		Lipschitz contraction estimate. Therefore,
		\[
		\begin{aligned}
			\int_{s+1}^{\infty}
			\left\|
			D_y^2P_{s,r}^x b_1(x,\cdot)
			\right\|_\infty\,dr
			&\leq
			C
			\int_{s+1}^{\infty}
			\exp\left\{
			-\int_{s+1}^{r}\lambda\,dq
			\right\}\,dr\leq C,
		\end{aligned}
		\]
		where the last inequality follows from the backward resolvent
		condition.
		
		Consequently, differentiation under the integral sign is justified,
		and
		\[
		D_y^2\Phi(s,x,y)
		=
		\int_s^\infty
		D_y^2P_{s,r}^x b_1(x,\cdot)(y)\,dr.
		\]
		It follows that
		\[
		\sup_{s,x}
		\|D_y^2\Phi(s,x,\cdot)\|_\infty
		\leq C.
		\]
		The proof is completed.
	\end{proof}
	Next we consider the Lipschitz property of $\bar{b}_1(t,\cdot)$. If the Lipschitz property of $\bar{b}_1(t,\cdot)$ is established, then that of $\Phi(s,\cdot,y)$ follows immediately.
	\begin{lemma}\label{lem:Lipschitz-evolution-measure}
		Assume that \ref{h1}-\ref{h4} hold. 
		Then there exists a constant $C>0$, independent of
		$t\in\mathbb{R}$ and $x_1,x_2\in\mathbb{R}^{d_1}$, such that
		\begin{equation}\label{eq:measure-x-lipschitz}
			W_1(\mu_t^{x_1},\mu_t^{x_2})
			\leq
			C|x_1-x_2|.
		\end{equation}
		Consequently,
		\begin{equation}\label{eq:averaged-drift-lipschitz}
			|\bar b_1(t,x_1)-\bar b_1(t,x_2)|
			\leq
			C|x_1-x_2|.
		\end{equation}
	\end{lemma}
	
	\begin{proof}
		Fix $t\in\mathbb{R}$ and $s<t$. Let
		$f:\mathbb{R}^{d_2}\to\mathbb{R}$ be such that
		\[
		\operatorname{Lip}(f)\leq 1.
		\]
		For $r\in[s,t]$, define
		\[
		u(r,y):=P_{r,t}^{x_2}f(y).
		\]
		Then $u$ solves the backward Kolmogorov equation
		\[
		\partial_r u(r,y)
		+
		\mathcal{L}_r^{x_2}u(r,y)
		=0,
		\qquad
		u(t,y)=f(y).
		\]
		Therefore, by the variation-of-constants formula,
		\begin{equation}\label{eq:duhamel-x}
			P_{s,t}^{x_1}f(y)-P_{s,t}^{x_2}f(y)
			=
			\int_s^t
			P_{s,r}^{x_1}
			\left[
			(\mathcal{L}_r^{x_1}-\mathcal{L}_r^{x_2})
			P_{r,t}^{x_2}f
			\right](y)\,dr.
		\end{equation}
		
		Set
		\[
		\delta:=|x_1-x_2|.
		\]
		The frozen generator is given by
		\begin{align}\label{eq:frozen-generator-x}
			\mathcal{L}_r^x u(y)
			={}&
			\langle b_2(r,x,y),\nabla u(y)\rangle
			+
			\frac12
			\operatorname{Tr}
			\left[
			a_2(r,x,y)D^2u(y)
			\right]
			\nonumber\\
			&+
			\int_{\mathbb{R}^{d_2}}
			\Big(
			u(y+g_2(r,x,y)z)-u(y)
			-\mathbf{1}_{\{|z|\leq1\}}
			\langle \nabla u(y),g_2(r,x,y)z\rangle
			\Big)
			\nu_2(dz),
		\end{align}
		where
		\[
		a_2(r,x,y)
		=
		\sigma_2(r,x,y)\sigma_2(r,x,y)^{\ast}.
		\]
		
		By the Lipschitz assumption on $b_2$,
		\[
		|b_2(r,x_1,y)-b_2(r,x_2,y)|
		\leq C\delta,
		\]
		and hence
		\begin{equation}\label{eq:drift-difference}
			|I_b|
			\leq
			C\delta
			\|\nabla u(r,\cdot)\|_{\infty}.
		\end{equation}
		
		Since $\sigma_2$ is uniformly bounded and globally Lipschitz in $x$,
		\[
		|a_2(r,x_1,y)-a_2(r,x_2,y)|
		\leq C\delta,
		\]
		so that
		\begin{equation}\label{eq:brownian-difference}
			|I_{\sigma}|
			\leq
			C\delta
			\|D^2u(r,\cdot)\|_{\infty}.
		\end{equation}
		
		It remains to estimate the jump part. Set
		\[
		G_i:=g_2(r,x_i,y),
		\qquad i=1,2.
		\]
		By the Lipschitz continuity of $g_2$,
		\[
		|G_1-G_2|
		\leq C\delta.
		\]
		For $|z|\leq1$, Taylor's formula yields
		\begin{align}
			\Big|
			u(y+G_1z)-u(y+G_2z)
			-
			\langle\nabla u(y),(G_1-G_2)z\rangle
			\Big|
			\nonumber
			\leq
			C\delta
			\|D^2u(r,\cdot)\|_{\infty}
			|z|^2.
			\label{eq:small-jump-generator-difference}
		\end{align}
		Since $\alpha_2<2$,
		\[
		\int_{|z|\leq1}|z|^2\nu_2(dz)<\infty.
		\]
		For $|z|>1$, we have
		\begin{equation}\label{eq:large-jump-generator-difference}
			|u(y+G_1z)-u(y+G_2z)|
			\leq
			C\delta
			\|\nabla u(r,\cdot)\|_{\infty}
			|z|.
		\end{equation}
		Since $\alpha_2>1$,
		\[
		\int_{|z|>1}|z|\nu_2(dz)<\infty.
		\]
		Consequently,
		\begin{equation}\label{eq:jump-difference}
			|I_g|
			\leq
			C\delta
			\left(
			\|\nabla u(r,\cdot)\|_{\infty}
			+
			\|D^2u(r,\cdot)\|_{\infty}
			\right).
		\end{equation}
		
		Combining
		\eqref{eq:drift-difference},
		\eqref{eq:brownian-difference},
		and \eqref{eq:jump-difference}, we obtain
		\begin{align}
			\left\|
			(\mathcal{L}_r^{x_1}-\mathcal{L}_r^{x_2})
			P_{r,t}^{x_2}f
			\right\|_{\infty}
			\nonumber
			\leq
			C\delta
			\left[
			\|\nabla_yP_{r,t}^{x_2}f\|_{\infty}
			+
			\|D_y^2P_{r,t}^{x_2}f\|_{\infty}
			\right].
			\label{eq:generator-difference-estimate}
		\end{align}
		
		We next split the integral in \eqref{eq:duhamel-x} into the short-time and long-time parts.
		
		For $r\in[(t-1)\vee s,t)$, by
		\eqref{eq:frozen-measure-w1-revised} and \eqref{hessian2},
		\[
		\|\nabla_yP_{r,t}^{x_2}f\|_{\infty}
		+
		\|D_y^2P_{r,t}^{x_2}f\|_{\infty}
		\leq
		C\left(1+(t-r)^{-\frac12}\right).
		\]
		Thus,
		\begin{align}\label{eq:short-time-part}
			\int_{(t-1)\vee s}^{t}
			\left\|
			(\mathcal{L}_r^{x_1}-\mathcal{L}_r^{x_2})
			P_{r,t}^{x_2}f
			\right\|_{\infty}
			\,dr
			\leq
			C\delta
			\int_0^1
			\left(1+\tau^{-\frac12}\right)\,d\tau
			\leq
			C\delta.
		\end{align}
		
		For $r\leq t-1$, by the semigroup property,
		\[
		P_{r,t}^{x_2}
		=
		P_{r,r+1}^{x_2}P_{r+1,t}^{x_2}.
		\]
		Applying the unit-time regularization estimates to
		$P_{r+1,t}^{x_2}f$, we obtain
		\begin{equation}		
			\|\nabla_yP_{r,t}^{x_2}f\|_{\infty}
			+
			\|D_y^2P_{r,t}^{x_2}f\|_{\infty}	
			\leq
			C\operatorname{Lip}
			\left(P_{r+1,t}^{x_2}f\right).
			\label{eq:long-regularization}
		\end{equation}
		Using the contraction estimate \eqref{eq:frozen-measure-w1-revised},
		\begin{align*}
			\operatorname{Lip}\big(P_{r+1,t}^{x_2} f\big) \leq C e^{-\lambda(t-r-1)} \operatorname{Lip}(f).
		\end{align*}
		Hence
		\begin{align}
			\int_s^{t-1}
			\left\|
			(\mathcal{L}_r^{x_1}-\mathcal{L}_r^{x_2})
			P_{r,t}^{x_2}f
			\right\|_{\infty}
			\,dr
			\leq
			C\delta
			\int_s^{t-1}
			e^{-\lambda(t-r-1)}
			dr\leq C\delta.
			\label{eq:long-time-part-1}
		\end{align}
		Combining
		\eqref{eq:duhamel-x},
		\eqref{eq:short-time-part},
		and \eqref{eq:long-time-part-1}, we obtain
		\begin{equation}\label{eq:semigroup-x-stability}
			|P_{s,t}^{x_1}f(y)-P_{s,t}^{x_2}f(y)|
			\leq
			C|x_1-x_2|,
		\end{equation}
		where $C$ is independent of $s,t,y,x_1,x_2$.
		
		Letting $s\to-\infty$, and using the pullback convergence
		\[
		P_{s,t}^{x}f(y)
		\longrightarrow
		\int_{\mathbb{R}^{d_2}}f(z)\,\mu_t^x(dz),
		\]
		we obtain
		\[
		\left|
		\int_{\mathbb{R}^{d_2}} f(y)\,\mu_t^{x_1}(dy)
		-
		\int_{\mathbb{R}^{d_2}} f(y)\,\mu_t^{x_2}(dy)
		\right|
		\leq
		C|x_1-x_2|.
		\]
		Taking the supremum over all $f$ satisfying
		$\operatorname{Lip}(f)\leq1$, the Kantorovich--Rubinstein duality yields
		\[
		W_1(\mu_t^{x_1},\mu_t^{x_2})
		\leq
		C|x_1-x_2|.
		\]
		
		Finally, by the definition
		\[
		\bar b_1(t,x)
		=
		\int_{\mathbb{R}^{d_2}}
		b_1(x,y)\,\mu_t^x(dy),
		\]
		we have
		\begin{align*}
			&
			|\bar b_1(t,x_1)-\bar b_1(t,x_2)|
			\\
			&\leq
			\left|
			\int_{\mathbb{R}^{d_2}}
			\big(
			b_1(x_1,y)-b_1(x_2,y)
			\big)
			\mu_t^{x_1}(dy)
			\right|
			\\
			&\quad
			+
			\left|
			\int_{\mathbb{R}^{d_2}}
			b_1(x_2,y)\,\mu_t^{x_1}(dy)
			-
			\int_{\mathbb{R}^{d_2}}
			b_1(x_2,y)\,\mu_t^{x_2}(dy)
			\right|
			\\
			&\leq
			C|x_1-x_2|
			+
			C W_1(\mu_t^{x_1},\mu_t^{x_2})
			\\
			&\leq
			C|x_1-x_2|.
		\end{align*}
		The proof is complete.
	\end{proof}
	\section{The First-Step Averaging Principle}
	We begin with introducing mollification of functions which will be used to tackle the difficulties related to the slow component. 
	Let $\rho\in C_c^\infty(\mathbb R^{d_1})$ be a standard mollifier
	satisfying
	\[
	\rho\geq0,\qquad
	\operatorname{supp}\rho\subset\overline{B_1(0)},\qquad
	\int_{\mathbb R^{d_1}}\rho(z)\,dz=1.
	\]
	For $n\geq1$, set
	\[
	\rho_n(z)=n^{d_1}\rho(nz),\qquad
	g_n(x,y)=\int_{\mathbb R^{d_1}}g(x-z,y)\rho_n(z)\,dz.
	\]
	We impose no pointwise domination of $D^\beta\rho$ by $\rho$.
	For each multi-index $\beta$ with $|\beta|=j\ge1$,
	compact support and smoothness imply
	\[
	A_\beta:=\int_{\mathbb R^{d_1}}|D^\beta\rho(w)|\,dw<\infty,
	\qquad
	B_\beta:=\int_{\mathbb R^{d_1}}|w|\,|D^\beta\rho(w)|\,dw<\infty.
	\]
	By scaling and cancellation,
	\begin{align}
		\int_{\mathbb R^{d_1}}|D^\beta\rho_n(z)|\,dz
		&=n^j A_\beta, \label{eq:mollifier-l1-revised}\\
		\int_{\mathbb R^{d_1}}|z|\,|D^\beta\rho_n(z)|\,dz
		&=n^{j-1}B_\beta, \label{eq:mollifier-weighted-revised}\\
		\int_{\mathbb R^{d_1}}D^\beta\rho_n(z)\,dz
		&=0. \label{eq:mollifier-cancellation-revised}
	\end{align}
	
	\begin{lemma}[Mollification in the slow variable]
		\label{yang}
		Suppose that $g:\mathbb R^{d_1}\times\mathbb R^{d_2}\to\mathbb R$
		satisfies
		\[
		|g(x_1,y)-g(x_2,y)|
		\leq L(1+|y|)|x_1-x_2|,
		\qquad x_1,x_2\in\mathbb R^{d_1},\quad y\in\mathbb R^{d_2}.
		\]
		Then $g_n(\cdot,y)\in C^\infty(\mathbb R^{d_1})$, and
		\begin{align}
			\|g_n(\cdot,y)-g(\cdot,y)\|_\infty
			&\leq Cn^{-1}(1+|y|), \label{eq:mollification-error-revised}\\
			\|\nabla_xg_n(\cdot,y)\|_\infty
			&\leq C(1+|y|), \label{2222}\\
			\|D_x^2g_n(\cdot,y)\|_\infty
			&\leq Cn(1+|y|). \label{3333}
		\end{align}
		Assume in addition that $1<\alpha_1<2$,
		$\nu_1(dz)=c_{\alpha_1,d_1}|z|^{-d_1-\alpha_1}\,dz$, and
		$\sup_x\|g_1(x)\|<\infty$. Define the slow jump operator by
		\[
		\mathcal J_1 f(x)
		=\int_{\mathbb R^{d_1}\setminus\{0\}}
		\left[
		f(x+g_1(x)z)-f(x)
		-\mathbf 1_{\{|z|\le1\}}\langle\nabla f(x),g_1(x)z\rangle
		\right]\nu_1(dz).
		\]
		Then
		\begin{equation}
			\|\mathcal J_1 g_n(\cdot,y)\|_\infty
			\leq Cn^{\alpha_1-1}(1+|y|).
			\label{4444}
		\end{equation}
		Here $C$ is independent of $n,x,y$.
		The estimates extend componentwise to vector-valued functions.
	\end{lemma}
	
	\begin{proof}
		The assumed Lipschitz estimate implies at most linear growth in $x$
		for fixed $y$, so convolution with the compactly supported kernel
		is well defined and smooth in $x$. First,
		\begin{align*}
			|g_n(x,y)-g(x,y)|
			&\leq L(1+|y|)\int_{\mathbb R^{d_1}}|z|\rho_n(z)\,dz\\
			&=Ln^{-1}(1+|y|)\int_{\mathbb R^{d_1}}|w|\rho(w)\,dw.
		\end{align*}
		This proves \eqref{eq:mollification-error-revised}.
		
		For a multi-index $\beta$ with $|\beta|=j\geq1$, differentiate the
		kernel in
		$g_n(x,y)=\int g(w,y)\rho_n(x-w)\,dw$.
		By \eqref{eq:mollifier-cancellation-revised},
		\[
		D_x^\beta g_n(x,y)
		=\int_{\mathbb R^{d_1}}
		[g(x-z,y)-g(x,y)]D^\beta\rho_n(z)\,dz.
		\]
		Consequently,
		\begin{align*}
			|D_x^\beta g_n(x,y)|
			&\leq L(1+|y|)
			\int_{\mathbb R^{d_1}}|z|\,|D^\beta\rho_n(z)|\,dz\\
			&=L B_\beta n^{j-1}(1+|y|).
		\end{align*}
		Taking $j=1$ and $j=2$ proves
		\eqref{2222} and
		\eqref{3333}.
		
		To prove \eqref{4444}, write
		\[
		\Delta_n(x,y;z)
		=g_n(x+g_1(x)z,y)-g_n(x,y)
		-\mathbf 1_{\{|z|\leq1\}}
		\langle\nabla_xg_n(x,y),g_1(x)z\rangle.
		\]
		For $|z|\leq n^{-1}$, Taylor's formula and
		\eqref{3333} yield
		\[
		|\Delta_n(x,y;z)|\leq Cn(1+|y|)|z|^2.
		\]
		For $|z|>n^{-1}$, the global gradient estimate
		\eqref{2222} bounds both the increment
		and the compensation term, when present, by
		\[
		|\Delta_n(x,y;z)|\leq C(1+|y|)|z|.
		\]
		It follows that
		\begin{align*}
			|\mathcal J_1g_n(\cdot,y)(x)|
			&\leq C(1+|y|)
			\left[
			n\int_{|z|\leq n^{-1}}|z|^2\nu_1(dz)
			+\int_{|z|>n^{-1}}|z|\nu_1(dz)
			\right]\\
			&\leq C(1+|y|)
			\left[
			n\int_0^{n^{-1}}r^{1-\alpha_1}\,dr
			+\int_{n^{-1}}^\infty r^{-\alpha_1}\,dr
			\right]\\
			&\leq Cn^{\alpha_1-1}(1+|y|).
		\end{align*}
		The integrals are finite because $1<\alpha_1<2$.
	\end{proof}
	\begin{theorem}\label{theorem5.4}
		Let \ref{h1}-\ref{h4} hold, for $m \in [1, \alpha_1\land\alpha_2 )$, then we have
		\[
		\mathbb{E}\left( \sup_{t \in [0,T]} \left| \int_0^t \big(b_1(X_s^\varepsilon,Y_s^\varepsilon) - \bar{b}_1(\frac{s}{\varepsilon},X_s^\varepsilon)\big) ds \right|^m \right)  \leq C_{T,m,x,y} \left( \varepsilon^{m(1 - \frac{1}{\alpha_2})} + \varepsilon^{\frac{m}{2}}\right). 
		\]
	\end{theorem}
	\begin{proof}
		From \eqref{possionequation} we know that there exist $\Phi(\cdot,\cdot,\cdot) $ such that
		\begin{align}
			\partial_s \Phi(s,x,y)
			+\mathcal{L}_2(s,x,y) \Phi(s,x,y) + b_1(x,y) - \bar{b}_1(s,x) = 0,
		\end{align} 
		For $n\geq1$, we mollify the corrector $\Phi$ only with respect to the
		slow variable $x$. Let $\rho\in C_c^\infty(\mathbb R^{d_1})$ be a
		standard mollifier defined as lemma \ref{yang}.
		Define
		\[
		\Phi_n(s,x,y)
		:=
		\int_{\mathbb R^{d_1}}
		\Phi(s,x-z,y)\rho_n(z)\,dz.
		\]
		Then $\Phi_n$ is smooth in the $x$-variable. Notice that, since the
		frozen generator $\mathcal L_2^x$ depends on $x$, the mollified
		function $\Phi_n$ does not in general solve the original Poisson
		equation. We therefore keep track of the corresponding commutator
		error below.
		
		Set
		\[
		F(s,x,y):=b_1(x,y)-\bar b_1(s,x).
		\]
		Then the non-autonomous Poisson equation can be written as
		\[
		(\partial_s+\mathcal L_2^x)\Phi(s,x,y)
		=
		-F(s,x,y).
		\]
		
		For $n\geq1$, define the mollifications in the slow variable $x$ by
		\[
		\Phi_n(s,x,y)
		:=
		\int_{\mathbb R^{d_1}}
		\Phi(s,x-z,y)\rho_n(z)\,dz
		\]
		and
		\[
		F_n(s,x,y)
		:=
		\int_{\mathbb R^{d_1}}
		F(s,x-z,y)\rho_n(z)\,dz,
		\]
		where
		\[
		\rho_n(z)=n^{d_1}\rho(nz),
		\qquad
		\rho\in C_c^\infty(\mathbb R^{d_1}),
		\qquad
		\rho\geq0,
		\qquad
		\int_{\mathbb R^{d_1}}\rho(z)\,dz=1.
		\]
		
		Since the mollification acts only on the $x$-variable, we have
		\[
		\partial_s\Phi_n(s,x,y)
		=
		\int_{\mathbb R^{d_1}}
		\partial_s\Phi(s,x-z,y)\rho_n(z)\,dz.
		\]
		However, since the frozen generator $\mathcal L_2^x$ depends on the
		parameter $x$, it does not in general commute with the
		$x$-mollification. We therefore introduce the commutator	
		\begin{equation*}
			R_n(s,x,y)
			:={}
			\mathcal L_2^x\Phi_n(s,x,y)
			-
			\int_{\mathbb R^{d_1}}
			\mathcal L_2^{x-z}\Phi(s,x-z,y)\rho_n(z)\,dz.
		\end{equation*}
		Using the Poisson equation, we obtain
		\[
		\begin{aligned}
			(\partial_s+\mathcal L_2^x)\Phi_n(s,x,y)
			&=
			\int_{\mathbb R^{d_1}}
			\bigl(
			\partial_s+\mathcal L_2^{x-z}
			\bigr)
			\Phi(s,x-z,y)\rho_n(z)\,dz
			+
			R_n(s,x,y)
			\\
			&=
			-F_n(s,x,y)+R_n(s,x,y).
		\end{aligned}
		\]
		Consequently,
		\[
		F(s,x,y)
		=
		F(s,x,y)-F_n(s,x,y)
		-
		(\partial_s+\mathcal L_2^x)\Phi_n(s,x,y)
		+
		R_n(s,x,y).
		\]
		
		We next estimate the two error terms. Since $b_1$ and
		$\bar b_1$ are uniformly Lipschitz continuous with respect to the
		slow variable $x$, we have,
		\[
		|F(s,x_1,y)-F(s,x_2,y)|
		\leq
		C|x_1-x_2|.
		\]
		Hence,
		\begin{align}\label{nnnnn}
			\begin{aligned}
				|F(s,x,y)-F_n(s,x,y)|
				&\leq
				\int_{\mathbb R^{d_1}}
				|F(s,x,y)-F(s,x-z,y)|\rho_n(z)\,dz
				\\
				&\leq
				C
				\int_{\mathbb R^{d_1}}
				|z|\rho_n(z)\,dz
				\\
				&\leq
				Cn^{-1}.
			\end{aligned}
		\end{align}
		
		We now estimate the commutator $R_n$. Write the frozen generator as
		\[
		\mathcal L_2^x
		=
		\mathcal B_2^x+\mathcal A_2^x+\mathcal J_2^x,
		\]
		where
		\[
		\mathcal B_2^x u(y)
		=
		\langle b_2(s,x,y),\nabla_yu(y)\rangle,
		\]
		\[
		\mathcal A_2^x u(y)
		=
		\frac12
		a_2(s,x,y):D_y^2u(y),
		\qquad
		a_2=\sigma_2\sigma_2^*,
		\]
		and $\mathcal J_2^x$ denotes the nonlocal part associated with
		$g_2(s,x,y)$ and the $\alpha_2$-stable L\'evy measure.
		
		Accordingly,
		\[
		R_n=R_n^{(b)}+R_n^{(a)}+R_n^{(J)}.
		\]
		For the drift part,
		\begin{align*}
			\begin{aligned}
				R_n^{(b)}(s,x,y)
				={}&
				\int_{\mathbb R^{d_1}}
				\Big\langle
				b_2(s,x,y)-b_2(s,x-z,y),
				\\
				&\hspace{32mm}
				\nabla_y\Phi(s,x-z,y)
				\Big\rangle
				\rho_n(z)\,dz .
			\end{aligned}
		\end{align*}
		By the Lipschitz continuity of $b_2$ with respect to $x$ and the
		uniform gradient estimate for the corrector,
		\[
		\sup_{s,x,y}|\nabla_y\Phi(s,x,y)|\leq C,
		\]
		we obtain
		\begin{align}\label{driftpart}
			|R_n^{(b)}(s,x,y)|
			\leq
			C\int_{\mathbb R^{d_1}}|z|\rho_n(z)\,dz
			\leq
			Cn^{-1}.
		\end{align}
		
		For the Brownian part,
		\[
		\begin{aligned}
			R_n^{(a)}(s,x,y)
			=
			\frac12
			\int_{\mathbb R^{d_1}}
			\bigl[
			a_2(s,x,y)-a_2(s,x-z,y)
			\bigr]
			:
			D_y^2\Phi(s,x-z,y)
			\rho_n(z)\,dz.
		\end{aligned}
		\]
		By the Lipschitz continuity of $\sigma_2$ with respect to $x$,
		the same property holds for $a_2=\sigma_2\sigma_2^*$.
		Together with
		\[
		\sup_{s,x,y}|D_y^2\Phi(s,x,y)|\leq C,
		\]
		this gives
		\begin{align}\label{brownpart}
			|R_n^{(a)}(s,x,y)|
			\leq
			C\int_{\mathbb R^{d_1}}|z|\rho_n(z)\,dz
			\leq
			Cn^{-1}.
		\end{align}
		
		For the nonlocal part, using the stable-like representation,
		we write
		\[
		\mathcal J_2^x u(y)
		=
		\int_{\mathbb R^{d_2}}
		\delta u(y;\xi)
		\frac{\kappa^x(s,y,\xi)}
		{|\xi|^{d_2+\alpha_2}}\,d\xi,
		\]
		where
		\[
		\delta u(y;\xi)
		=
		u(y+\xi)-u(y)
		-
		\mathbf 1_{\{|\xi|\leq1\}}
		\langle\nabla_yu(y),\xi\rangle .
		\]
		Hence,
		\[
		\begin{aligned}
			R_n^{(J)}(s,x,y)
			={}&
			\int_{\mathbb R^{d_1}}
			\int_{\mathbb R^{d_2}}
			\delta\Phi(s,x-z,y;\xi)
			\\
			&\quad\times
			\frac{
				\kappa^x(s,y,\xi)
				-
				\kappa^{x-z}(s,y,\xi)
			}
			{|\xi|^{d_2+\alpha_2}}
			\,d\xi\,\rho_n(z)\,dz.
		\end{aligned}
		\]
		By the Lipschitz continuity of $g_2$ with respect to $x$ and its
		uniform non-degeneracy, the corresponding stable-like kernel satisfies
		\[
		\left|
		\kappa^x(s,y,\xi)
		-
		\kappa^{x-z}(s,y,\xi)
		\right|
		\leq
		C|z|.
		\]
		Moreover, for $|\xi|\leq1$, the boundedness of $D_y^2\Phi$ yields
		\[
		|\delta\Phi(s,x-z,y;\xi)|
		\leq
		C|\xi|^2,
		\]
		whereas for $|\xi|>1$, the boundedness of $\nabla_y\Phi$ gives
		\[
		|\delta\Phi(s,x-z,y;\xi)|
		\leq
		C|\xi|.
		\]
		Since $\alpha_2\in(1,2)$,
		\[
		\int_{|\xi|\leq1}
		\frac{|\xi|^2}
		{|\xi|^{d_2+\alpha_2}}\,d\xi
		<\infty
		\qquad \text{and}\qquad
		\int_{|\xi|>1}
		\frac{|\xi|}
		{|\xi|^{d_2+\alpha_2}}\,d\xi
		<\infty.
		\]
		Consequently,
		\begin{align}\label{nonlocalpart}
			\begin{aligned}
				|R_n^{(J)}(s,x,y)|
				&\leq
				C
				\int_{\mathbb R^{d_1}}
				|z|\rho_n(z)\,dz
				\\
				&\leq
				Cn^{-1}.
			\end{aligned}
		\end{align}
		
		Combining \eqref{driftpart}, \eqref{brownpart} and \eqref{nonlocalpart}, we conclude that
		\[
		|R_n(s,x,y)|
		\leq
		Cn^{-1}.
		\]
		Together with \eqref{nnnnn}, this yields
		\[
		|F(s,x,y)-F_n(s,x,y)|
		+
		|R_n(s,x,y)|
		\leq
		Cn^{-1}(1+|y|).
		\]
		
		Finally, we justify the application of It\^o's formula. By the
		regularity results established after Corollaries \ref{cor:gradient-Phi} and \ref{cor:hessian-Phi}, the corrector
		$\Phi$ possesses the required first- and second-order regularity in
		the fast variable $y$, and
		\[
		\sup_{s,x,y}
		\left(
		|\nabla_y\Phi(s,x,y)|
		+
		|D_y^2\Phi(s,x,y)|
		\right)
		<\infty.
		\]
		Moreover, the Poisson equation
		\[
		\partial_s\Phi
		=
		-\mathcal L_2^x\Phi-F
		\]
		provides the required regularity in the time variable. Since
		$\Phi_n$ is obtained by mollification only with respect to the slow
		variable $x$, it is smooth in $x$, while the above regularity in
		$(s,y)$ is preserved under the mollification. Therefore,
		$\Phi_n$ has sufficient regularity to apply the It\^o formula.
		We have
		\begin{align}\label{514}
			\begin{aligned}
				&b_1(X_s^\varepsilon,Y_s^\varepsilon)
				-
				\bar b_1\left(
				\frac{s}{\varepsilon},
				X_s^\varepsilon
				\right)
				\\
				={}&
				(F-F_n)\left(
				\frac{s}{\varepsilon},
				X_s^\varepsilon,
				Y_s^\varepsilon
				\right)
				\\
				&-
				(\partial_s+\mathcal L_2^{X_s^\varepsilon})
				\Phi_n\left(
				\frac{s}{\varepsilon},
				X_s^\varepsilon,
				Y_s^\varepsilon
				\right)
				\\
				&+
				R_n\left(
				\frac{s}{\varepsilon},
				X_s^\varepsilon,
				Y_s^\varepsilon
				\right).
			\end{aligned}
		\end{align}
		The second term on the right-hand side of \eqref{514} can now be
		handled directly by applying It\^o's formula to
		\[
		\Phi_n\left(
		\frac{t}{\varepsilon},
		X_t^\varepsilon,
		Y_t^\varepsilon
		\right).
		\]
		We have that,
		\begin{align*}
			\Phi_n(\frac{t}{\varepsilon},X_t^\varepsilon,Y_t^\varepsilon) &=\Phi_n(0,x,y) + \frac{1}{\varepsilon} \int_0^t (\partial_s\Phi_n(\frac{s}{\varepsilon},X_s^\varepsilon,Y_s^\varepsilon)+\mathcal{L}_2(\frac{s}{\varepsilon},X_s^\varepsilon,Y_s^\varepsilon)\Phi_n(\frac{s}{\varepsilon},X_s^\varepsilon,Y_s^\varepsilon)) ds \\
			&  \quad+  \int_0^t \mathcal{L}_1(\frac{s}{\varepsilon},X_s^\varepsilon,Y_s^\varepsilon)\Phi_n(\frac{s}{\varepsilon},X_s^\varepsilon,Y_s^\varepsilon) ds \\
			&\quad+\int_{0}^{t}\langle \partial_x \Phi_n(\frac{s}{\varepsilon},X_s^\varepsilon,Y_s^\varepsilon), \sigma_1(X_s^\varepsilon)\rangle dW^1_ s\\
			&\quad+\frac{1}{\sqrt{\varepsilon}}\int_{0}^{t}\langle \partial_y \Phi_n(\frac{s}{\varepsilon},X_s^\varepsilon,Y_s^\varepsilon), \sigma_2 (\frac{s}{\varepsilon},X_s^\varepsilon,Y_s^\varepsilon)\rangle dW^2_ s\\
			&\quad+ M_{n}^{1,\varepsilon} + M_{n}^{2,\varepsilon},
		\end{align*}
		where $\mathcal{L}_1(t,x, y)$is defined as
		\begin{equation}
			\begin{aligned}
				\mathcal{L}_1(t,x, y)\Phi(t,x, y) :=& \langle b_1(x, y), \nabla_x \Phi(t,x, y) \rangle \\
				&+ \frac{1}{2} \sum_{i,j=1}^{d_1} \left( \sigma_1(x) \sigma_1(x)^\top \right)_{ij} \frac{\partial^2 \Phi}{\partial x_i \partial x_j}(t,x, y) \\
				&+ \int_{\mathbb{R}^{d_1}} \Big[ \Phi(t,x+ g_1(x)z, y  ) - \Phi(t,x, y) \\
				&\qquad - \mathbf{1}_{\{|z| \leq 1\}} \langle \nabla_x \Phi(t,x, y), g_1(x) z \rangle \Big] \nu_1(dz).
			\end{aligned}
		\end{equation}
		and $M_{n,t}^{1,\varepsilon}$, $M_{n,t}^{2,\varepsilon}$ are two $\mathcal{F}_t$ martingales defined as
		\begin{align*}
			M_{n}^{1,\varepsilon} &= \int_0^t \int_{\mathbb{R}^{d_1}} \big( \Phi_n(\frac{s}{\varepsilon},X_{s-}^\varepsilon + g_1 \cdot x,Y_{s-}^\varepsilon) - \Phi_n(\frac{s}{\varepsilon},X_{s-}^\varepsilon,Y_{s-}^\varepsilon) \big) \tilde{N}^1(ds,dx), \\
			M_{n}^{2,\varepsilon} &= \int_0^t \int_{\mathbb{R}^{d_2}} \big( \Phi_n(\frac{s}{\varepsilon},X_{s-}^\varepsilon,Y_{s-}^\varepsilon + \varepsilon^{-\frac{1}{\alpha_2}} g_2 \cdot y) - \Phi_n(\frac{s}{\varepsilon},X_{s-}^\varepsilon,Y_{s-}^\varepsilon) \big) \tilde{N}^2(ds,dy),
		\end{align*}
		where $\tilde{N}^1$, $\tilde{N}^2$ are compensated Poisson measures.
		
		Above calculations lead us to
		\begin{align*}
			&\int_0^t (\partial_s\Phi_n(\frac{s}{\varepsilon},X_s^\varepsilon,Y_s^\varepsilon)+\mathcal{L}_2(\frac{s}{\varepsilon},X_s^\varepsilon,Y_s^\varepsilon)\Phi_n(\frac{s}{\varepsilon},X_s^\varepsilon,Y_s^\varepsilon)) ds\\
			&=-\varepsilon\Bigg[\Phi_n(0,x,y) - \Phi_n(\frac{t}{\varepsilon},X_t^\varepsilon,Y_t^\varepsilon)  +  \int_0^t \mathcal{L}_1(\frac{s}{\varepsilon},X_s^\varepsilon,Y_s^\varepsilon)\Phi_n(\frac{s}{\varepsilon},X_s^\varepsilon,Y_s^\varepsilon) ds \\
			&\quad+\int_{0}^{t}\langle \partial_x \Phi_n(\frac{s}{\varepsilon},X_s^\varepsilon,Y_s^\varepsilon), \sigma_1(X_s^\varepsilon)\rangle dW^1_ s+\frac{1}{\sqrt{\varepsilon}}\int_{0}^{t}\langle \partial_y \Phi_n(\frac{s}{\varepsilon},X_s^\varepsilon,Y_s^\varepsilon), \sigma_2 (\frac{s}{\varepsilon},X_s^\varepsilon,Y_s^\varepsilon)\rangle dW^2_ s\\
			&\quad+ M_{n}^{1,\varepsilon} + M_{n}^{2,\varepsilon} \Bigg],
		\end{align*}
		taking expectation from both sides above , we have
		\begin{align}\label{111}
			\begin{aligned}
				&\mathbb{E}\left( \sup_{t \in [0,T]} \left| \int_0^t (\partial_s\Phi_n(\frac{s}{\varepsilon},X_s^\varepsilon,Y_s^\varepsilon)+\mathcal{L}_2(\frac{s}{\varepsilon},X_s^\varepsilon,Y_s^\varepsilon)\Phi_n(\frac{s}{\varepsilon},X_s^\varepsilon,Y_s^\varepsilon)) ds \right|^m \right) \\
				&\leq C_{T,m} \cdot \varepsilon^m \Bigg[ \mathbb{E}\left( \sup_{t \in [0,T]} \left| \Phi_n(0,x,y) - \Phi_n(\frac{t}{\varepsilon},X_t^\varepsilon,Y_t^\varepsilon) \right|^m \right) \\
				&\quad+ \mathbb{E}\left( \int_0^T \left| \mathcal{L}_1(\frac{s}{\varepsilon},X_s^\varepsilon,Y_s^\varepsilon)\Phi_n(\frac{s}{\varepsilon},X_s^\varepsilon,Y_s^\varepsilon)   \right|^m ds \right) \\
				&\quad+ \mathbb{E}\left( \sup_{t \in [0,T]} \left|\int_{0}^{t}\langle \partial_x \Phi_n(\frac{s}{\varepsilon},X_s^\varepsilon,Y_s^\varepsilon), \sigma_1(X_s^\varepsilon)\rangle dW^1_ s \right|^m \right) \\
				&\quad+ \mathbb{E}\left( \sup_{t \in [0,T]} \left| \frac{1}{\sqrt{\varepsilon}}\int_{0}^{t}\langle \partial_y \Phi_n(\frac{s}{\varepsilon},X_s^\varepsilon,Y_s^\varepsilon), \sigma_2 (\frac{s}{\varepsilon},X_s^\varepsilon,Y_s^\varepsilon)\rangle dW^2_ s \right|^m \right) \\
				&\quad+ \mathbb{E}\left( \sup_{t \in [0,T]} \left| M_{n,t}^{1,\varepsilon} \right|^m \right) + \mathbb{E}\left( \sup_{t \in [0,T]} \left| M_{n}^{2,\varepsilon} \right|^m \right) \\
				&= C_{T,m} \cdot \varepsilon^m (I_1 + I_2 + I_3 + I_4 + I_5 + I_6 ), 
			\end{aligned}
		\end{align}
		we will estimate the above terms respectively.
		
		For $I_1$, by Theorem \ref{throrem4.5}, Lemma \ref{lemma3.4} and Theorem \ref{theorem3.2},
		\begin{align}\label{333}
			\begin{aligned}
				I_1 &= \mathbb{E}\left( \sup_{t \in [0,T]} \left| \Phi_n(0,x,y) - \Phi_n(\frac{t}{\varepsilon},X_t^\varepsilon,Y_t^\varepsilon)  \right|^m \right)\\
				&\leq C_{T,m}\left (1+ \mathbb{E}( \sup_{t \in [0,T]} |X_t^\varepsilon|^m)+\mathbb{E}+( \sup_{t \in [0,T]} |Y_t^\varepsilon|^m )\right ) +C (1+|x|^p+|y|^p)\\
				&\leq C_{T,m} (1+|x|^m+|y|^m) \varepsilon^{-\frac{m}{\alpha_2}}+C (1+|x|^m+|y|^m). 
			\end{aligned}
		\end{align}
		For $I_2$, since we have Assumption \ref{h3} and theorem \ref{theorem3.2}, by (\ref{2222}) ,(\ref{3333})  and (\ref{4444})  in Lemma \ref{yang}, moment estimates (\ref{well2}),
		\begin{align}
			\begin{aligned}
				I_2 &= \mathbb{E}\left( \int_0^T \left| \mathcal{L}_1(\frac{s}{\varepsilon},X_s^\varepsilon,Y_s^\varepsilon)\Phi_n(\frac{s}{\varepsilon},X_s^\varepsilon,Y_s^\varepsilon) \right|^m ds \right) \\
				&\leq C_{T,m} \mathbb{E}\left( \int_0^T \left| \big(b_1(X_s^\varepsilon,Y_s^\varepsilon), \nabla_x \Phi_n(\frac{s}{\varepsilon},X_s^\varepsilon,Y_s^\varepsilon)\big) \right|^m ds \right) \\
				&\quad+ C_{T,m} \mathbb{E}\left( \int_0^T \left| \sigma_1(X_s^\varepsilon)D^2(\Phi_n(\frac{s}{\varepsilon},X_s^\varepsilon,Y_s^\varepsilon))\sigma_1(X_s^\varepsilon)^\top \right|^m ds \right) \\
				&\quad + C_{T,m} \mathbb{E}\left( \int_0^T \left| -(-\Delta_x)^{\frac{\alpha_1}{2}} \Phi_n(\frac{s}{\varepsilon},X_s^\varepsilon,Y_s^\varepsilon) \right|^m ds \right) \\
				&\leq C_{T,m} + C_{T,m} n^{m}+ C_{T,m} \varepsilon^{-\frac{m}{\alpha_2}}.\label{444}
			\end{aligned}
		\end{align}
		For $I_3$ and $I_4$, we can deduce from Burkholder-Davies-Gundy's inequality, the upper boundness of $\sigma_2$ and $\sigma_2$ in \ref{h1}, (\ref{well1}) and (\ref{3333}),
		\begin{align}
			\begin{aligned}
				I_3+I_4=\,&\mathbb{E}\left( \sup_{t \in [0,T]} \left|\int_{0}^{t}\langle \partial_x \Phi_n(\frac{s}{\varepsilon},X_s^\varepsilon,Y_s^\varepsilon), \sigma_1(X_s^\varepsilon)\rangle dW^1_ s \right|^m \right) \\
				&\quad+ \mathbb{E}\left( \sup_{t \in [0,T]} \left| \frac{1}{\sqrt{\varepsilon}}\int_{0}^{t}\langle \partial_y \Phi_n(\frac{s}{\varepsilon},X_s^\varepsilon,Y_s^\varepsilon), \sigma_2 (\frac{s}{\varepsilon},X_s^\varepsilon,Y_s^\varepsilon)\rangle dW^2_ s \right|^m \right)\\
				\leq& \,\mathbb{E}\left(  \int_{0}^{T}| \partial_x \Phi_n(\frac{s}{\varepsilon},X_s^\varepsilon,Y_s^\varepsilon) \sigma_1(X_s^\varepsilon)|^2 d s  \right)^{\frac{m}{2}} \\
				&\quad+ \varepsilon^{-\frac{m}{2}}\left(  \int_{0}^{T}|  \partial_y \Phi_n(\frac{s}{\varepsilon},X_s^\varepsilon,Y_s^\varepsilon) \sigma_2 (\frac{s}{\varepsilon},X_s^\varepsilon,Y_s^\varepsilon)|^2 d s  \right)^{\frac{m}{2}}  \\
				\leq& \,C_{T,m}\left(1 + |y|^m + \varepsilon^{-m/\alpha_2}\right)
				+C_T \varepsilon^{-\frac{m}{2}}.
			\end{aligned}
		\end{align}
		For $I_5$ and $I_6$, we can deduce from Burkholder-Davies-Gundy's inequality, the upper boundness of $g_1$ in \ref{h1}, (\ref{well2}) and (\ref{3333}),
		\begin{align}\label{555}
			I_5
			&= \mathbb E\left(
			\sup_{t\in[0,T]}
			\left|M_n^{1,\varepsilon}\right|^m
			\right)
			\notag\\
			&\leq C_{T,m}\,
			\mathbb E\left[
			\left(
			\int_0^T\int_{|z|\leq 1}
			\left|
			\Phi_n\left(
			\frac{s}{\varepsilon},
			X_s^\varepsilon+g_1\cdot z,
			Y_s^\varepsilon
			\right)
			-
			\Phi_n\left(
			\frac{s}{\varepsilon},
			X_s^\varepsilon,
			Y_s^\varepsilon
			\right)
			\right|^2
			\nu_1(dz)\,ds
			\right)^{m/2}
			\right]
			\notag	\\
			&\quad
			+ C_{T,m}\,
			\int_0^T\mathbb E\left[
			\int_{|z|>1}
			\left|
			\Phi_n\left(
			\frac{s}{\varepsilon},
			X_s^\varepsilon+g_1\cdot z,
			Y_s^\varepsilon
			\right)
			-
			\Phi_n\left(
			\frac{s}{\varepsilon},
			X_s^\varepsilon,
			Y_s^\varepsilon
			\right)
			\right|^m
			\nu_1(dz)
			\right]ds
			\notag\\
			&\leq C_{T,m}\,
			\mathbb E\left[
			\left(
			\int_0^T\int_{|z|\leq 1}
			|g_1\cdot z|^2
			\left(1+|Y_s^\varepsilon|^2\right)
			\nu_1(dz)\,ds
			\right)^{m/2}
			\right]
			\notag\\
			&\quad
			+ C_{T,m}\,
			\int_0^T\mathbb E\left[
			\int_{|z|>1}
			|g_1\cdot z|^m
			\left(1+|Y_s^\varepsilon|^m\right)
			\nu_1(dz)
			\right]ds
			\notag\\
			&\leq C_{T,m}\,
			\mathbb E\left[
			\left(
			\int_0^T
			\left(1+|Y_s^\varepsilon|^2\right)\,ds
			\right)^{m/2}
			\right]
			+ C_{T,m}\int_0^T
			\mathbb E\left(1+|Y_s^\varepsilon|^m\right)\,ds
			\notag\\
			&\leq C_{T,m}\left(
			1+\mathbb E\sup_{s\in[0,T]}|Y_s^\varepsilon|^m
			\right)
			\notag\\
			&\leq C_{T,m}\left(
			1+\varepsilon^{-m/\alpha_2}+|y|^m
			\right).
		\end{align}
		then, by $\nabla_y u_n = (\nabla_y u) * \rho^n$ again, the upper boundedness of $g_2$ in \ref{h1}, (\ref{3333}) in lemma \ref{yang},
		\begin{align}\label{666}
			I_6
			&= \mathbb E\left(
			\sup_{t\in[0,T]}
			\left|M_n^{2,\varepsilon}\right|^m
			\right)
			\notag\\
			&\leq C_{T,m}\,
			\mathbb E\left[
			\left(
			\int_0^T\int_{|z|\leq 1}
			\left|
			\Phi_n\left(
			\frac{s}{\varepsilon},
			X_s^\varepsilon,
			Y_s^\varepsilon
			+\varepsilon^{-1/\alpha_2}g_2\cdot z
			\right)
			-
			\Phi_n\left(
			\frac{s}{\varepsilon},
			X_s^\varepsilon,
			Y_s^\varepsilon
			\right)
			\right|^2
			\nu_2(dz)\,ds
			\right)^{m/2}
			\right]
			\notag\\
			&\quad
			+ C_{T,m}\,
			\int_0^T\mathbb E\left[
			\int_{|z|>1}
			\left|
			\Phi_n\left(
			\frac{s}{\varepsilon},
			X_s^\varepsilon,
			Y_s^\varepsilon
			+\varepsilon^{-1/\alpha_2}g_2\cdot z
			\right)
			-
			\Phi_n\left(
			\frac{s}{\varepsilon},
			X_s^\varepsilon,
			Y_s^\varepsilon
			\right)
			\right|^m
			\nu_2(dz)
			\right]ds
			\notag\\
			&\leq C_{T,m}\,\varepsilon^{-m/\alpha_2}
			\left[
			\mathbb E\left(
			\int_0^T\int_{|z|\leq 1}
			|g_2\cdot z|^2\,\nu_2(dz)\,ds
			\right)^{m/2}
			+
			\int_0^T\mathbb E\left(
			\int_{|z|>1}
			|g_2\cdot z|^m\,\nu_2(dz)
			\right)ds
			\right]
			\notag\\
			&\leq C_{T,m}\,\varepsilon^{-m/\alpha_2}
			\left[
			\left(
			T\int_{|z|\leq 1}|z|^2\,\nu_2(dz)
			\right)^{m/2}
			+
			T\int_{|z|>1}|z|^m\,\nu_2(dz)
			\right]
			\notag\\
			&\leq C_{T,m}\,\varepsilon^{-m/\alpha_2}.
		\end{align}
		Combining (\ref{111})-(\ref{666}) together, take $n = \varepsilon^{-\frac{1}{\alpha_2}}$, we have
		\begin{align*}
			&\mathbb{E}\left( \sup_{t \in [0,T]} \left| \int_0^t \big(b_1(X_s^\varepsilon,Y_s^\varepsilon) - \bar{b}_1(\frac{s}{\varepsilon},X_s^\varepsilon)\big) ds \right|^m \right) 
			\leq C_{T,m,x,y} \left( \varepsilon^{m(1 - \frac{1}{\alpha_2})} + \varepsilon^{\frac{m}{2}}\right),
		\end{align*}
		and so the proof is complete.
	\end{proof}
	
	\begin{proof}[Proof of Theorem 2.3]
		By \eqref{avco}, we get an averaging equation
		\[
		d\bar{X}^\varepsilon_t = \bar{b}_1(\frac{t}{\varepsilon}, \bar{X}_t)dt + \sigma_1(\bar{X}_t)dW_t^1+g_1(\bar{X}_t)dL^1_t,
		\]
		so that
		\begin{align*}
			X_t^\varepsilon - \bar{X}^\varepsilon_t &= \int_0^t \big( b( X_s^\varepsilon, Y_s^\varepsilon) - \bar{b}(\frac{s}{\varepsilon}, \bar{X}_s) \big) ds+\int_0^t(\sigma_1(X_s^\varepsilon)-\sigma_1(\bar{X}_s))dW^1_s+\int_0^t(g_1(X_s^\varepsilon)-g_1(\bar{X}_s))dL^1_s\\
			&= \int_0^t \big( b(X_s^\varepsilon, Y_s^\varepsilon) - \bar{b}(\frac{s}{\varepsilon}, X_s^\varepsilon) \big) ds + \int_0^t \big( \bar{b}(\frac{s}{\varepsilon}, X_s^\varepsilon) - \bar{b}(\frac{s}{\varepsilon}, \bar{X}_s) \big) ds\\
			&+\int_0^t(\sigma_1(X_s^\varepsilon)-\sigma_1(\bar{X}_s))dW^1_s+\int_0^t(g_1(X_s^\varepsilon)-g_1(\bar{X}_s))dL^1_s.
		\end{align*}
		
		By Lipschitz continuity of $\bar{b}$ and $\sigma_1$ , using B-D-G inequality, we have for $m \in [1, \alpha_1)$,
		\begin{align*}
			\mathbb{E} \left( \sup_{t \in [0,T]} |X_t^\varepsilon - \bar{X}^\varepsilon_t|^m \right)
			&\leq \mathbb{E} \left( \sup_{t \in [0,T]} \left| \int_0^t \big( b(X_s^\varepsilon, Y_s^\varepsilon) - \bar{b}(\frac{s}{\varepsilon}, X_s^\varepsilon) \big) ds \right|^m \right) \\
			&\quad + C_{T,m} \mathbb{E} \int_0^T |X_s^\varepsilon - \bar{X}^\varepsilon_s|^m ds\\
			&\quad+\mathbb{E} \left( \sup_{t \in [0,T]} \left| \int_0^t(g_1(X_s^\varepsilon)-g_1(\bar{X}^\varepsilon_s))dL^1_s \right|^m \right).
		\end{align*}
		By [\citenum{ChengHaoZhang2025}, Lemma A.4] and Jensen's inequality, Lipschitz continuity of $g_1$, we have
		\begin{align*}
			&\mathbb{E} \left( \sup_{t \in [0,T]} \left| \int_0^t(g_1(X_s^\varepsilon)-g_1(\bar{X}^\varepsilon_s))dL^1_s \right|^m \right)\\
			&\leq C_m \mathbb{E} \left( \sup_{t \in [0,T]} \left| \int_0^t\int_{|z|\leq1}(g_1(X_s^\varepsilon)-g_1(\bar{X}^\varepsilon_s))z\tilde{N}_1(ds,dz) \right|^m \right)\\
			&\quad+ C_m \mathbb{E} \left( \sup_{t \in [0,T]} \left| \int_0^t\int_{|z|\leq1}(g_1(X_s^\varepsilon)-g_1(\bar{X}^\varepsilon_s))zv_1(dz)ds \right|^m \right)\\
			&\quad+ C_m \mathbb{E} \left( \sup_{t \in [0,T]} \left| \int_0^t\int_{|z|>1}(g_1(X_s^\varepsilon)-g_1(\bar{X}^\varepsilon_s))zN_1(ds,dz) \right|^m \right)\\
			&\leq \frac{4-m}{2-m} C_m \mathbb{E} \Bigl( \sup_{t \in [0,T]} \Bigl( \int_0^t\int_{|z|\leq 1} \big|g_1(X_s^\varepsilon)-g_1(\bar{X}^\varepsilon_s)\big|^2 z^2 \tilde{N}_1(ds,dz) \Bigr)^{\frac{m}{2}} \Bigr)\\
			&\quad+C_m \mathbb{E} \Bigl( \sup_{t \in [0,T]} \Bigl( \int_0^t\int_{|z|\leq 1} \big|g_1(X_s^\varepsilon)-g_1(\bar{X}^\varepsilon_s)\big|^m z^m \tilde{N}_1(ds,dz) \Bigr) \Bigr)\\
			&\quad+C_m \mathbb{E} \Bigl( \sup_{t \in [0,T]} \Bigl( \int_0^t\int_{|z|\geq 1} \big|g_1(X_s^\varepsilon)-g_1(\bar{X}^\varepsilon_s)\big|^p z^p \tilde{N}_1(ds,dz) \Bigr) \Bigr)\\
			&\leq C_m \int_{0}^{T}|X_s^\varepsilon-\bar{X}^\varepsilon_s|^m ds.
		\end{align*}
		Then from Gronwall's inequality, and Theorem \ref{theorem5.4} , we know that
		\[
		\begin{aligned}
			\mathbb{E} \left( \sup_{t \in [0,T]} |X_t^\varepsilon - \bar{X}^\varepsilon_t|^m \right)
			&\leq \mathbb{E} \left( \sup_{t \in [0,T]} \left| \int_0^t \big( b_1( X_s^\varepsilon, Y_s^\varepsilon) - \bar{b}_1(\frac{s}{\varepsilon}, {X}^\varepsilon_s) \big) ds \right|^m \right) \\
			&\leq C_{T,m,x,y} \left( \varepsilon^{m(1 - \frac{1}{\alpha_2})} + \varepsilon^{\frac{m}{2}}\right). \qquad\qquad\qedhere
		\end{aligned}
		\]
	\end{proof}
	
	\section{The Second-Step Averaging Principle}
	
	We retain the standing assumptions on the slow--fast system and the
	frozen dynamics from the preceding sections. In particular,
	$1<\alpha_1,\alpha_2<2$, the spatial Lipschitz constants are uniform,
	and $g_1$ is bounded. The Brownian motions and stable processes are
	the same as in the original system. The argument below preserves the
	multiplicative coefficients $g_1$ and $g_2$.
	
	Write $\mathbb T^q=(\mathbb R/2\pi\mathbb Z)^q$, equipped with the
	distance
	\[
	d_{\mathbb T^q}(\theta,\theta')
	=\min_{\ell\in\mathbb Z^q}|\theta-\theta'+2\pi\ell|
	\]
	and normalized Haar measure $d\theta$.
	Thus $\int_{\mathbb T^q}d\theta=1$.
	All phase additions below are understood modulo $2\pi\mathbb Z^q$.
	\begin{hypo}[Quasi-periodic structure of the original coefficients]
		\label{h6}
		In addition to the standing assumptions, suppose that the following
		conditions hold.
		\begin{enumerate}
			\item There are $q\geq1$, $\omega\in\mathbb R^q$, and profiles
			$\widetilde b_2$, $\widetilde\sigma_2$, $\widetilde g_2$ such that
			\begin{align}
				b_2(t,x,y)&=\widetilde b_2(\omega t,x,y),\nonumber\\
				\sigma_2(t,x,y)&=\widetilde\sigma_2(\omega t,x,y),\nonumber\\
				g_2(t,x,y)&=\widetilde g_2(\omega t,x,y).
				\label{qpnew:profiles}
			\end{align}
			The standing bounds, spatial regularity, dissipativity, and
			non-degeneracy assumptions on the fast coefficients hold for these
			profiles uniformly over all $\theta\in\mathbb T^q$.
			\item For some $L_\theta<\infty$, uniformly in $x,y,\theta,\theta'$,
			\begin{align}
				&|\widetilde b_2(\theta,x,y)-\widetilde b_2(\theta',x,y)|
				+\|\widetilde\sigma_2(\theta,x,y)
				-\widetilde\sigma_2(\theta',x,y)\|\nonumber\\
				&\qquad
				+\|\widetilde g_2(\theta,x,y)-\widetilde g_2(\theta',x,y)\|
				\leq L_\theta d_{\mathbb T^q}(\theta,\theta').
				\label{qpnew:phase-lipschitz}
			\end{align}
			In particular, the constant in this condition is independent of $y$.
			\item There are $\gamma>0$ and $\nu\geq0$ such that
			\begin{equation}
				|k\cdot\omega|\geq\gamma |k|^{-\nu},
				\qquad k\in\mathbb Z^q\setminus\{0\}.
				\label{qpnew:diophantine}
			\end{equation}
			Here $|k|$ is the Euclidean norm.
		\end{enumerate}
	\end{hypo}
	For a fixed phase offset $\vartheta\in\mathbb T^q$ and a fixed
	$x\in\mathbb R^{d_1}$, consider
	\begin{align}
		dY_r={}&\widetilde b_2(\vartheta+\omega r,x,Y_r)\,dr
		+\widetilde\sigma_2(\vartheta+\omega r,x,Y_r)\,dW_r^2
		\nonumber\\
		&+\widetilde g_2(\vartheta+\omega r,x,Y_{r-})\,dL_r^2.
		\label{qpnew:frozen}
	\end{align}
	The offset $\vartheta$ is fixed, while the actual phase
	$\vartheta+\omega r$ evolves with time. Let
	$P_{s,t}^{\vartheta,x}$ denote its transition family and
	$P_{s,t}^{\vartheta,x,*}$ its action on probability measures.
	
	We use the uniform frozen-equation estimates established in the
	preceding sections. Their precise forms needed here are
	\begin{align}
		W_1\bigl(P_{s,t}^{\vartheta,x,*}\lambda_1,
		P_{s,t}^{\vartheta,x,*}\lambda_2\bigr)
		&\leq C_0e^{-\lambda_0(t-s)}W_1(\lambda_1,\lambda_2),
		\label{qpnew:frozen-contraction}\\
		\sup_{\vartheta\in\mathbb T^q}\sup_{s\leq t}
		\int |y|^p\,
		P_{s,t}^{\vartheta,x,*}\delta_0(dy)
		&\leq C_p(1+|x|^p),\qquad 1<p<\alpha_2.
		\label{qpnew:frozen-moment}
	\end{align}
	The constants $C_0,\lambda_0,C_p$ are independent of the phase
	offset and the time endpoints. These estimates apply uniformly
	because the profiles satisfy the standing fast-equation assumptions
	uniformly in phase. The contraction in
	\eqref{qpnew:frozen-contraction} is a statement about laws; no
	contraction under the synchronous coupling is required.
	
	For later use, we record a stochastic-integral estimate. If $L$ is
	a symmetric $\alpha$-stable process, $1\leq p<\alpha<2$, and $F$ is
	a predictable matrix-valued integrand, then
	\begin{align}
		\mathbb E\sup_{a\leq r\leq b}
		\left|\int_a^r F_u\,dL_u\right|^p
		\leq C_p\left\{
		\mathbb E\left(\int_a^b\|F_u\|^2\,du\right)^{p/2}
		+\int_a^b\mathbb E\|F_u\|^p\,du\right\}.
		\label{qpnew:stable-maximal}
	\end{align}
	Indeed, split the integral into $|z|\leq1$ and $|z|>1$.
	The compensated small-jump integral is estimated by the
	$L^p(L^2)$ martingale maximal inequality, and the compensated
	large-jump integral by the $L^p(L^p)$ inequality. The corresponding
	constants are finite since
	\begin{equation}
		\int_{|z|\leq1}|z|^2\,\nu_L(dz)<\infty,
		\qquad
		\int_{|z|>1}|z|^p\,\nu_L(dz)<\infty.
		\label{qpnew:levy-integrability}
	\end{equation}
	For $p=1$, the large-jump estimate also follows directly from
	the total variation and its compensator. Symmetry and $\alpha>1$
	allow the large jumps to be compensated without an additional drift.
	The estimates can first be applied to stopped integrands and then
	extended by localization.
	
	\begin{lemma}[Finite-time stability with respect to phase]
		\label{qpnew:finite-phase}
		Fix $S>0$ and $1<p<\alpha_2$. Two solutions of
		\eqref{qpnew:frozen} with offsets $\vartheta,\vartheta'$ can be
		constructed from the same initial value at time $t-S$ and driven
		by the same Brownian motion and stable process so that
		\begin{equation}
			\mathbb E\sup_{t-S\le r\leq t}|Y_r-Z_r|^p
			\leq C_{S,p}\,d_{\mathbb T^q}(\vartheta,\vartheta')^p.
			\label{qpnew:finite-path}
		\end{equation}
		The constant is independent of the common initial value and of $x,t$.
		Consequently, for every common initial law $\lambda$ with finite
		$p$-th moment,
		\begin{equation}
			W_1\bigl(P_{t-S,t}^{\vartheta,x,*}\lambda,
			P_{t-S,t}^{\vartheta',x,*}\lambda\bigr)
			\leq C_{S,p}\,d_{\mathbb T^q}(\vartheta,\vartheta').
			\label{qpnew:finite-law}
		\end{equation}
	\end{lemma}
	
	\begin{proof}
		Set $\delta=d_{\mathbb T^q}(\vartheta,\vartheta')$ and
		$D_r=Y_r-Z_r$. Translation invariance of the torus distance gives
		\[
		d_{\mathbb T^q}(\vartheta+\omega r,\vartheta'+\omega r)=\delta.
		\]
		Splitting each coefficient difference first in the state variable
		and then in the phase variable, the spatial Lipschitz condition and
		\eqref{qpnew:phase-lipschitz} imply
		\begin{equation}
			|\Delta b_2(r)|+\|\Delta\sigma_2(r)\|
			+\|\Delta g_2(r)\|
			\leq C(|D_r|+\delta),
			\label{qpnew:coefficient-difference}
		\end{equation}
		with left limits in the predictable jump integrand.
		
		On any subinterval $[a,a+h]\subset[t-S,t]$, put
		$A(a,h)=\mathbb E\sup_{a\leq r\leq a+h}|D_r|^p$.
		The drift inequality, the Brownian BDG inequality, and
		\eqref{qpnew:stable-maximal} give
		\begin{align}
			A(a,h)\leq{}&C\mathbb E|D_a|^p+C(h^p+h^{p/2}+h)\bigl(A(a,h)+\delta^p\bigr).
			\label{qpnew:local-phase}
		\end{align}
		For example, the Brownian contribution is at most
		$Ch^{p/2}(A(a,h)+\delta^p)$; the small- and large-jump
		contributions are bounded by
		$C(h^{p/2}+h)(A(a,h)+\delta^p)$.
		Choose $h_0\in(0,1]$ so that the coefficient of $A(a,h)$ on the
		right of \eqref{qpnew:local-phase} is at most $1/2$ whenever
		$h\leq h_0$. This therefore yields
		\[
		A(a,h)\leq K_p\bigl(\mathbb E|D_a|^p+\delta^p\bigr),
		\qquad h\leq h_0.
		\]
		Partition $[t-S,t]$ into finitely many intervals of length at most
		$h_0$. Since $D_{t-S}=0$, iteration gives
		\eqref{qpnew:finite-path}. Uniformity in the initial value follows
		from the unweighted bound \eqref{qpnew:phase-lipschitz}.
		Finally, the endpoint pair $(Y_t,Z_t)$ is a coupling of the two
		transition laws, and hence
		\[
		W_1\bigl(P_{t-S,t}^{\vartheta,x,*}\lambda,
		P_{t-S,t}^{\vartheta',x,*}\lambda\bigr)
		\leq\mathbb E|D_t|
		\leq\bigl(\mathbb E\sup_{t-S\leq r\leq t}|D_r|^p\bigr)^{1/p}.
		\]
		This proves \eqref{qpnew:finite-law}.
	\end{proof}
	
	\begin{proposition}[Phase representation and Lipschitz inheritance]
		\label{qpnew:inheritance}
		There is a canonical evolution system of measures
		$\{\mu_t^{\vartheta,x}\}_{t\in\mathbb R}$ for
		\eqref{qpnew:frozen}. With
		\begin{equation}
			\nu_\theta^x:=\mu_0^{\theta,x},
			\qquad
			B(\theta,x):=\int_{\mathbb R^{d_2}} b_1(x,y)\,\nu_\theta^x(dy),
			\label{qpnew:B-definition}
		\end{equation}
		one has
		\begin{equation}
			\mu_t^{\vartheta,x}=\nu_{\vartheta+\omega t}^x,
			\qquad \bar b_1(t,x)=B(\omega t,x).
			\label{qpnew:phase-representation}
		\end{equation}
		Moreover, uniformly in $\theta,\theta',x,x_1,x_2$,
		\begin{align}
			W_1(\nu_\theta^x,\nu_{\theta'}^x)
			&\leq C d_{\mathbb T^q}(\theta,\theta'),
			\label{qpnew:measure-phase-lip}\\
			|B(\theta,x)-B(\theta',x)|
			&\leq C d_{\mathbb T^q}(\theta,\theta'),
			\label{qpnew:B-phase-lip}\\
			|B(\theta,x_1)-B(\theta,x_2)|
			&\leq C|x_1-x_2|,
			\label{qpnew:B-x-lip}\\
			|B(\theta,x)|&\leq C(1+|x|).
			\label{qpnew:B-growth}
		\end{align}
	\end{proposition}
	
	\begin{proof}
		The canonical family is the pullback limit in $W_1$,
		\begin{equation}
			\mu_t^{\vartheta,x}
			=\lim_{s\to-\infty}P_{s,t}^{\vartheta,x,*}\delta_0.
			\label{qpnew:pullback}
		\end{equation}
		For $r<s<t$, the semigroup property and
		\eqref{qpnew:frozen-contraction} imply
		\begin{align*}
			W_1\bigl(P_{r,t}^{\vartheta,x,*}\delta_0,
			P_{s,t}^{\vartheta,x,*}\delta_0\bigr)\leq C_0e^{-\lambda_0(t-s)}
			W_1\bigl(P_{r,s}^{\vartheta,x,*}\delta_0,\delta_0\bigr)
			\leq C(1+|x|)e^{-\lambda_0(t-s)}.
		\end{align*}
		Thus the pullback laws form a Cauchy family in $W_1$.
		They inherit the uniform $p$-moment bounds, and their limit has the
		evolution property. Contraction also proves uniqueness among
		evolution systems with time-uniform first moments for fixed $x$.
		Time translation in \eqref{qpnew:frozen} gives
		\[
		P_{s,t}^{\vartheta,x}
		=P_{s-t,0}^{\vartheta+\omega t,x}.
		\]
		Taking pullback limits proves \eqref{qpnew:phase-representation}.
		
		Fix $x,\vartheta,\vartheta'$ and put
		$\delta=d_{\mathbb T^q}(\vartheta,\vartheta')$.
		Choose $S>0$ such that $C_0e^{-\lambda_0S}\leq1/2$.
		Using the evolution property, insert the intermediate law
		$P_{t-S,t}^{\vartheta,x,*}\mu_{t-S}^{\vartheta',x}$.
		The triangle inequality, \eqref{qpnew:frozen-contraction}, and
		Lemma~\ref{qpnew:finite-phase} yield
		\begin{align*}
			W_1(\mu_t^{\vartheta,x},\mu_t^{\vartheta',x})
			\leq\tfrac12 W_1(\mu_{t-S}^{\vartheta,x},
			\mu_{t-S}^{\vartheta',x})+C_S\delta.
		\end{align*}
		The supremum of the left-hand side over $t\in\mathbb R$ is finite
		by the first-moment bounds. Taking this supremum and absorbing its
		half therefore gives
		\[
		\sup_t W_1(\mu_t^{\vartheta,x},\mu_t^{\vartheta',x})
		\leq2C_S\delta.
		\]
		Setting $t=0$ proves \eqref{qpnew:measure-phase-lip}.
		Since $b_1(x,\cdot)$ is globally Lipschitz, taking the infimum
		over all couplings of $\nu_\theta^x$ and $\nu_{\theta'}^x$ gives
		\[
		|B(\theta,x)-B(\theta',x)|
		\leq \operatorname{Lip}_y(b_1)\,
		W_1(\nu_\theta^x,\nu_{\theta'}^x),
		\]
		which proves \eqref{qpnew:B-phase-lip}.
		
		The slow-parameter estimate proved for the frozen measures in the
		preceding sections applies to every phase-shifted equation,
		uniformly in the offset. At time zero it reads
		\[
		W_1(\nu_\theta^{x_1},\nu_\theta^{x_2})\leq C|x_1-x_2|.
		\]
		Combining it with the spatial Lipschitz continuity of $b_1$ proves
		\eqref{qpnew:B-x-lip}. The moment bound at $x=0$ gives
		$\sup_\theta|B(\theta,0)|<\infty$; hence
		\eqref{qpnew:B-growth} follows from \eqref{qpnew:B-x-lip}.
	\end{proof}
	
	Define
	\begin{equation}
		\widehat b_1(x):=\int_{\mathbb T^q}B(\theta,x)\,d\theta,
		\qquad H(\theta,x):=B(\theta,x)-\widehat b_1(x).
		\label{qpnew:mean-and-H}
	\end{equation}
	Then $\widehat b_1$ is globally Lipschitz with linear growth,
	$\int_{\mathbb T^q}H(\theta,x)\,d\theta=0$, and
	\begin{align}
		\begin{aligned}
		|H(\theta,x)|&\leq C(1+|x|),\\
		|H(\theta,x_1)-H(\theta,x_2)|&\leq C|x_1-x_2|,\\
		|H(\theta,x)-H(\theta',x)|
		&\leq C d_{\mathbb T^q}(\theta,\theta').
		\label{qpnew:H-bounds}
		\end{aligned}
	\end{align}
	The last constant is independent of $x$.
	\begin{lemma}[Positive trigonometric approximation]
		\label{qpnew:jackson}
		For every integer $N\geq1$, there is a trigonometric polynomial
		$H_N(\theta,x)$ with Fourier support in
		$\{k\in\mathbb Z^q:|k|\leq c_qN\}$, $c_q=2\sqrt q$, such that
		\begin{align}
			\int_{\mathbb T^q}H_N(\theta,x)\,d\theta&=0,\nonumber\\
			\sup_{\theta,x}|H_N(\theta,x)-H(\theta,x)|&\leq CN^{-1},
			\label{qpnew:jackson-error}\\
			|H_N(\theta,x)|&\leq C(1+|x|),\nonumber\\
			|H_N(\theta,x_1)-H_N(\theta,x_2)|&\leq C|x_1-x_2|.
			\label{qpnew:HN-spatial}
		\end{align}
		All constants are independent of $N$.
	\end{lemma}
	
	\begin{proof}
		On $\mathbb T=\mathbb R/(2\pi\mathbb Z)$, put
		\[
		D_N(u)=\sum_{j=0}^{N-1}e^{iju},\qquad
		J_N(u)=\frac{|D_N(u)|^4}{\int_{\mathbb T}|D_N(v)|^4\,dv}.
		\]
		Integrals here use normalized Haar measure.
		The kernel $J_N$ is nonnegative, has integral one, and has degree
		at most $2(N-1)$. For $|u|\leq\pi$,
		\[
		|D_N(u)|\leq\min\{N,C|u|^{-1}\},\qquad
		cN^3\leq\int_{\mathbb T}|D_N(v)|^4\,dv\leq CN^3.
		\]
		The lower bound follows by integrating over $|v|\leq1/N$,
		where $|D_N(v)|\geq cN$. Consequently,
		\[
		J_N(u)\leq\frac{CN}{(1+N|u|)^4},\qquad
		\int_{\mathbb T}|u|\,J_N(u)\,du\leq CN^{-1}.
		\]
		Take the product kernel
		$K_N(\varphi)=\prod_{j=1}^qJ_N(\varphi_j)$ and define
		\[
		H_N(\theta,x)
		=\int_{\mathbb T^q}K_N(\varphi)H(\theta-\varphi,x)\,d\varphi.
		\]
		This convolution uses nonnegative weights of total mass one, and
		$\int d_{\mathbb T^q}(\varphi,0)K_N(\varphi)\,d\varphi\leq C/N$.
		The phase Lipschitz bound in \eqref{qpnew:H-bounds} therefore
		proves \eqref{qpnew:jackson-error}. Positivity and normalization
		preserve the zero mean, growth bound, and spatial Lipschitz bound.
		The degree bound in each coordinate gives $|k|\leq c_qN$.
	\end{proof}
	
	\begin{lemma}[A finite small-divisor estimate]
		\label{qpnew:small-divisors}
		Under \eqref{qpnew:diophantine},
		\begin{equation}
			\sum_{0<|k|\leq c_qN}\frac1{|k\cdot\omega|^2}
			\leq CN^{2\nu},\qquad N\geq1.
			\label{qpnew:reciprocal-sum}
		\end{equation}
	\end{lemma}
	
	\begin{proof}
		For distinct $k,\ell$ with $|k|,|\ell|\leq c_qN$, including
		the zero vector,
		\[
		|k\cdot\omega-\ell\cdot\omega|
		\geq\gamma|k-\ell|^{-\nu}
		\geq a_N:=\gamma(2c_qN)^{-\nu}.
		\]
		Thus the retained frequencies, together with zero, are separated
		by at least $a_N$. List the positive frequencies as
		$0<\lambda_1<\cdots<\lambda_M$. Then $\lambda_j\geq ja_N$,
		and hence
		\[
		\sum_{j=1}^M\lambda_j^{-2}
		\leq a_N^{-2}\sum_{j=1}^\infty j^{-2}\leq 2a_N^{-2}.
		\]
		The same estimate holds for the absolute values of the negative
		frequencies. Their total contribution is at most
		$4a_N^{-2}\leq CN^{2\nu}$.
	\end{proof}
	
	Write
	\[
	H_N(\theta,x)=\sum_{0<|k|\leq c_qN}(H_N)_k(x)e^{ik\cdot\theta}
	\]
	and define the finite corrector
	\begin{equation}
		\chi_N(\theta,x)
		:=\sum_{0<|k|\leq c_qN}
		\frac{(H_N)_k(x)}{i\,k\cdot\omega}e^{ik\cdot\theta}.
		\label{qpnew:corrector-definition}
	\end{equation}
	The denominator is nonzero by \eqref{qpnew:diophantine}.
	Conjugate symmetry shows that $\chi_N$ is real-valued whenever
	$H_N$ is real-valued. Since the sum is finite,
	\begin{equation}
		\omega\cdot\nabla_\theta\chi_N=H_N.
		\label{qpnew:corrector-equation}
	\end{equation}
	
	\begin{lemma}[Size and spatial Lipschitz bounds]
		\label{qpnew:corrector-bounds}
		There exists a constant $C>0$ independent of $N, \theta, x, x_1, x_2$ such that
		\begin{align}
			|\chi_N(\theta,x)|&\leq CN^\nu(1+|x|),
			\label{qpnew:chi-growth}\\
			|\chi_N(\theta,x_1)-\chi_N(\theta,x_2)|
			&\leq CN^\nu|x_1-x_2|.
			\label{qpnew:chi-lip}
		\end{align}
	\end{lemma}
	
	\begin{proof}
		Cauchy--Schwarz, Lemma~\ref{qpnew:small-divisors}, and Parseval identity give
		\begin{align*}
			|\chi_N(\theta,x)|
			&\leq\left(\sum_k|(H_N)_k(x)|^2\right)^{1/2}
			\left(\sum_{0<|k|\leq c_qN}|k\cdot\omega|^{-2}\right)^{1/2}\\
			&\leq CN^\nu\left(\int_{\mathbb T^q}|H_N(\varphi,x)|^2
			\,d\varphi\right)^{1/2}
			\leq CN^\nu(1+|x|).
		\end{align*}
		For the difference at two slow variables, the same inequalities
		applied to the difference of the Fourier coefficients give
		\begin{align*}
			&|\chi_N(\theta,x_1)-\chi_N(\theta,x_2)|\\
			&\quad\leq CN^\nu
			\left(\int_{\mathbb T^q}|H_N(\varphi,x_1)
			-H_N(\varphi,x_2)|^2\,d\varphi\right)^{1/2}
			\leq CN^\nu|x_1-x_2|,
		\end{align*}
		which complete the proof.
	\end{proof}
	Let $\rho\in C_c^\infty(\mathbb R^{d_1})$ be nonnegative, supported
	in the unit ball, with integral one. For integers $n\geq1$, put
	$\rho_n(z)=n^{d_1}\rho(nz)$ and define
	\begin{align}
		\chi_{N,n}(\theta,x)&=
		\int\chi_N(\theta,x-z)\rho_n(z)\,dz,\nonumber\\
		H_{N,n}(\theta,x)&=
		\int H_N(\theta,x-z)\rho_n(z)\,dz.
		\label{qpnew:mollification}
	\end{align}
	
	\begin{lemma}[Smoothed corrector estimates]
		\label{qpnew:smoothing}
		For $N,n\geq1$,
		\begin{align}
			\omega\cdot\nabla_\theta\chi_{N,n}&=H_{N,n},
			\label{qpnew:smoothed-equation}\\
			\sup_{\theta,x}|H(\theta,x)-H_{N,n}(\theta,x)|
			&\leq C(N^{-1}+n^{-1}),
			\label{qpnew:total-approximation}\\
			|\chi_{N,n}(\theta,x)|&\leq CN^\nu(1+|x|),
			\label{qpnew:smoothed-growth}\\
			\|\nabla_x\chi_{N,n}\|_\infty&\leq CN^\nu,
			\label{qpnew:gradient}\\
			\|D_x^2\chi_{N,n}\|_\infty&\leq CN^\nu n.
			\label{qpnew:hessian}
		\end{align}
		For the slow jump generator
		\begin{align}
			J_1f(x):=\int_{\mathbb R^{d_1}\setminus\{0\}}
			\bigl[f(x+g_1(x)z)-f(x)
			-\mathbf1_{\{|z|\leq1\}}\langle\nabla f(x),g_1(x)z\rangle
			\bigr]\,\nu_1(dz),
			\label{qpnew:slow-jump-generator}
		\end{align}
		one also has
		\begin{equation}
			\sup_{\theta,x}|J_1\chi_{N,n}(\theta,\cdot)(x)|
			\leq CN^\nu n^{\alpha_1-1}.
			\label{qpnew:jump-generator-bound}
		\end{equation}
	\end{lemma}
	
	\begin{proof}
		Convolution in $x$ commutes with $\omega\cdot\nabla_\theta$, giving
		\eqref{qpnew:smoothed-equation}. By \eqref{qpnew:HN-spatial},
		$\|H_{N,n}-H_N\|_\infty\leq C/n$, and
		\eqref{qpnew:total-approximation} follows from
		\eqref{qpnew:jackson-error}. In particular, this smoothing error
		has no factor $N^\nu$, since it is the error for $H_N$.
		The growth bound is preserved by convolution.
		
		For $j=1,2$, the cancellation of the derivatives of $\rho_n$ gives
		\[
		D_x^j\chi_{N,n}(\theta,x)
		=\int\bigl[\chi_N(\theta,x-z)-\chi_N(\theta,x)\bigr]
		D^j\rho_n(z)\,dz.
		\]
		Using \eqref{qpnew:chi-lip} and scaling,
		\[
		|D_x^j\chi_{N,n}(\theta,x)|
		\leq CN^\nu\int |z|\,|D^j\rho_n(z)|\,dz
		\leq CN^\nu n^{j-1}.
		\]
		This proves \eqref{qpnew:gradient}--\eqref{qpnew:hessian}.
		
		For \eqref{qpnew:jump-generator-bound}, write
		$f(x)=\chi_{N,n}(\theta,x)$. Boundness of $g_1$ and the derivative
		boundness give
		\begin{align*}
			|f(x+g_1(x)z)-f(x)-\langle\nabla f(x),g_1(x)z\rangle|
			&\leq CN^\nu n|z|^2, &&|z|\leq n^{-1},\\
			|f(x+g_1(x)z)-f(x)-\langle\nabla f(x),g_1(x)z\rangle|
			&\leq CN^\nu|z|, &&n^{-1}<|z|\leq1,\\
			|f(x+g_1(x)z)-f(x)|&\leq CN^\nu|z|, &&|z|>1.
		\end{align*}
		The radial form of $\nu_1$ therefore yields
		\begin{align*}
			|J_1f(x)|
			\leq CN^\nu\left(
			n\int_0^{1/n}r^{1-\alpha_1}\,dr
			+\int_{1/n}^1r^{-\alpha_1}\,dr
			+\int_1^\infty r^{-\alpha_1}\,dr\right)\leq CN^\nu n^{\alpha_1-1},
		\end{align*}
		where $1<\alpha_1<2$ is used in the three integrals.
	\end{proof}
	Consider the intermediate averaged equation
	\begin{equation}
		d\bar X_t^\varepsilon
		=B(\omega t/\varepsilon,\bar X_t^\varepsilon)\,dt
		+\sigma_1(\bar X_t^\varepsilon)\,dW_t^1
		+g_1(\bar X_{t-}^\varepsilon)\,dL_t^1,
		\qquad\bar X_0^\varepsilon=x.
		\label{qpnew:intermediate}
	\end{equation}
	Global Lipschitz continuity and linear growth give existence and
	uniqueness and, for $1\leq m<\alpha_1$,
	\begin{equation}
		\mathbb E\sup_{0\leq t\leq T}|\bar X_t^\varepsilon|^m
		\leq C_{T,m}(1+|x|^m),\qquad 0<\varepsilon\leq1.
		\label{qpnew:slow-moment}
	\end{equation}
	This is the slow-equation moment estimate used earlier; it also
	follows directly from BDG, \eqref{qpnew:stable-maximal}, and
	short-interval absorption in the integral equation. The constants
	are uniform in $\varepsilon$ because the drift bounds are uniform
	in phase.
	
	For a process $Z$, write
	\[
	\|Z\|_{m,T}
	:=\left(\mathbb E\sup_{0\leq t\leq T}|Z_t|^m\right)^{1/m}.
	\]
	
	\begin{proposition}[Oscillatory integral estimate]
		\label{qpnew:oscillatory}
		Let $T>0$ and $1\leq m<\alpha_1$. Set
		\begin{equation}
			R_t^\varepsilon
			:=\int_0^tH(\omega s/\varepsilon,\bar X_s^\varepsilon)\,ds.
			\label{qpnew:residual}
		\end{equation}
		For all integers $N,n\geq1$,
		\begin{align}
			\|R^\varepsilon\|_{m,T}
			&\leq C_{T,m,x}\left[
			N^{-1}+n^{-1}
			+\varepsilon N^\nu(1+n+n^{\alpha_1-1})\right]\nonumber\\
			&\leq C_{T,m,x}
			\left[N^{-1}+n^{-1}+\varepsilon N^\nu n\right].
			\label{qpnew:three-errors}
		\end{align}
		In particular,
		\begin{equation}
			\|R^\varepsilon\|_{m,T}
			\leq C_{T,m,x}\varepsilon^{1/(\nu+2)}.
			\label{qpnew:residual-rate}
		\end{equation}
	\end{proposition}
	
	\begin{proof}
		By \eqref{qpnew:total-approximation},
		\begin{equation}
			\left\|\int_0^\cdot
			(H-H_{N,n})(\omega s/\varepsilon,\bar X_s^\varepsilon)\,ds
			\right\|_{m,T}
			\leq C_T(N^{-1}+n^{-1}).
			\label{qpnew:residual-approximation}
		\end{equation}
		Put $f=\chi_{N,n}$ and apply the It\^o--L\'evy formula componentwise
		to $\varepsilon f(\omega t/\varepsilon,\bar X_t^\varepsilon)$.
		Using \eqref{qpnew:smoothed-equation}, one obtains
		\begin{align}
			&\int_0^tH_{N,n}(\omega s/\varepsilon,\bar X_s^\varepsilon)\,ds
			\nonumber\\
			&\quad=\varepsilon f(\omega t/\varepsilon,\bar X_t^\varepsilon)
			-\varepsilon f(0,x)\nonumber\\
			&\qquad-\varepsilon\int_0^t
			D_xf(\omega s/\varepsilon,\bar X_s^\varepsilon)
			B(\omega s/\varepsilon,\bar X_s^\varepsilon)\,ds\nonumber\\
			&\qquad-\frac{\varepsilon}{2}\int_0^t
			\operatorname{Tr}\!\left[
			\sigma_1\sigma_1^*(\bar X_s^\varepsilon)
			D_x^2f(\omega s/\varepsilon,\bar X_s^\varepsilon)
			\right]\,ds\nonumber\\
			&\qquad-\varepsilon\int_0^t
			J_1f(\omega s/\varepsilon,\cdot)(\bar X_s^\varepsilon)\,ds
			-\varepsilon M_t^W-\varepsilon M_t^L.
			\label{qpnew:ito-identity}
		\end{align}
		For vector-valued $f$, the trace is taken componentwise. Here
		\begin{align*}
			M_t^W&=\int_0^t
			D_xf(\omega s/\varepsilon,\bar X_s^\varepsilon)
			\sigma_1(\bar X_s^\varepsilon)\,dW_s^1,\\
			M_t^L&=\int_0^t\int_{\mathbb R^{d_1}\setminus\{0\}}
			\bigl[f(\omega s/\varepsilon,
			\bar X_{s-}^\varepsilon+g_1(\bar X_{s-}^\varepsilon)z)
			-f(\omega s/\varepsilon,\bar X_{s-}^\varepsilon)\bigr]
			\,\widetilde N^1(ds,dz).
		\end{align*}
		
		We estimate the terms in \eqref{qpnew:ito-identity} separately.
		The endpoint terms are bounded in $\|\cdot\|_{m,T}$ by
		$C_{T,m,x}\varepsilon N^\nu$ using
		\eqref{qpnew:smoothed-growth} and \eqref{qpnew:slow-moment}.
		The drift term has the same bound by
		\eqref{qpnew:gradient}, \eqref{qpnew:B-growth}, and the slow moment
		estimate. Boundness of $\sigma_1$ and
		\eqref{qpnew:hessian} give the deterministic bound
		$C_T\varepsilon N^\nu n$ for the Brownian second-order term.
		The jump-generator term is bounded by
		$C_T\varepsilon N^\nu n^{\alpha_1-1}$ by
		\eqref{qpnew:jump-generator-bound}.
		
		For the Brownian martingale, BDG and boundness of $\sigma_1$
		give $\|M^W\|_{m,T}\leq C_{T,m}N^\nu$.
		For the jump martingale, its integrand has absolute value at most
		$CN^\nu|z|$, by \eqref{qpnew:gradient} and boundness of $g_1$.
		The small-jump $L^m(L^2)$ and large-jump $L^m(L^m)$ maximal
		inequalities, with
		\[
		\int_{|z|\leq1}|z|^2\nu_1(dz)<\infty,
		\qquad \int_{|z|>1}|z|^m\nu_1(dz)<\infty,
		\]
		therefore imply $\|M^L\|_{m,T}\leq C_{T,m}N^\nu$.
		Combining these bounds with
		\eqref{qpnew:residual-approximation} proves the first line of
		\eqref{qpnew:three-errors}. Since $n\geq1$ and $\alpha_1-1<1$,
		$1+n+n^{\alpha_1-1}\leq3n$, proving its second line.
		
		Finally, choose
		\begin{equation}
			N=n=\left\lceil\varepsilon^{-1/(\nu+2)}\right\rceil.
			\label{qpnew:parameter-choice}
		\end{equation}
		Then $N^{-1}+\varepsilon N^{\nu+1}
		\leq C\varepsilon^{1/(\nu+2)}$, which gives
		\eqref{qpnew:residual-rate}.
	\end{proof}
	Define $\widehat X$ on the same noise space by
	\begin{equation}
		d\widehat X_t=\widehat b_1(\widehat X_t)\,dt
		+\sigma_1(\widehat X_t)\,dW_t^1
		+g_1(\widehat X_{t-})\,dL_t^1,
		\qquad\widehat X_0=x.
		\label{qpnew:autonomous}
	\end{equation}
	The global Lipschitz and linear growth properties ensure a unique
	strong solution.
	
	\begin{proposition}[Second averaging rate]
		\label{qpnew:second-averaging}
		For $T>0$ and $1\leq m<\alpha_1$,
		\begin{equation}
			\|\bar X^\varepsilon-\widehat X\|_{m,T}
			\leq C_{T,m,x}\varepsilon^{1/(\nu+2)}.
			\label{qpnew:second-rate}
		\end{equation}
	\end{proposition}
	
	\begin{proof}
		Let $Z_t=\bar X_t^\varepsilon-\widehat X_t$. By adding and
		subtracting $\widehat b_1(\bar X_t^\varepsilon)$,
		\begin{align*}
			Z_t={}&R_t^\varepsilon
			+\int_0^t\bigl[\widehat b_1(\bar X_s^\varepsilon)
			-\widehat b_1(\widehat X_s)\bigr]\,ds\\
			&+\int_0^t\bigl[\sigma_1(\bar X_s^\varepsilon)
			-\sigma_1(\widehat X_s)\bigr]\,dW_s^1\\
			&+\int_0^t\bigl[g_1(\bar X_{s-}^\varepsilon)
			-g_1(\widehat X_{s-})\bigr]\,dL_s^1.
		\end{align*}
		On an interval $[a,a+h]\subset[0,T]$, set
		$A=\mathbb E\sup_{a\leq r\leq a+h}|Z_r|^m$.
		Lipschitz continuity, BDG, and \eqref{qpnew:stable-maximal} imply
		\[
		A\leq C\mathbb E|Z_a|^m
		+C\mathbb E\sup_{a\leq r\leq a+h}
		|R_r^\varepsilon-R_a^\varepsilon|^m
		+C(h^m+h^{m/2}+h)A.
		\]
		Choose a fixed $h_0>0$ small enough to absorb the last term.
		Since $Z_0=0$ and
		$\sup_{a\leq r\leq a+h}|R_r^\varepsilon-R_a^\varepsilon|
		\leq2\sup_{r\leq T}|R_r^\varepsilon|$,
		iteration over finitely many intervals yields
		\[
		\|\bar X^\varepsilon-\widehat X\|_{m,T}
		\leq C_{T,m}\|R^\varepsilon\|_{m,T}.
		\]
		Proposition~\ref{qpnew:oscillatory} now proves
		\eqref{qpnew:second-rate}.
	\end{proof}
	
	\begin{proof}[Proof of Theorem 2.4]
		The first averaging estimate from Section~5, followed by the
		slow-equation stability comparison, gives
		\[
		\|X^\varepsilon-\bar X^\varepsilon\|_{m,T}
		\leq C_{T,m,x,y}\varepsilon^{1-1/\alpha_2}.
		\]
		Minkowski's inequality and
		Proposition~\ref{qpnew:second-averaging} imply
		\[
		\|X^\varepsilon-\widehat X\|_{m,T}
		\leq\|X^\varepsilon-\bar X^\varepsilon\|_{m,T}
		+\|\bar X^\varepsilon-\widehat X\|_{m,T},
		\]
		which proves \eqref{qpnew:main-rate-moment}.
		For $0<\varepsilon\leq1$, the sum of the two powers is bounded
		by $2\varepsilon^\kappa$, giving
		\eqref{qpnew:main-rate-moment}.
	\end{proof}
	
	\begin{remark}[Range of the exponent]
		\label{qpnew:rate-scope}
		The proof applies throughout $1<\alpha_2<2$. The exponent obtained
		by this argument is $\kappa$, and no optimality is asserted.
		The first-step exponent $1-1/\alpha_2$ is retained when
		$1-1/\alpha_2\leq1/(\nu+2)$; otherwise the present second-step
		estimate determines the resulting bound.
		For $\omega=(1,\sqrt2)$ one may take $\nu=1$, so the estimate is
		\[
		\|X^\varepsilon-\widehat X\|_{m,T}
		\leq C_{T,m,x,y}
		(\varepsilon^{1-1/\alpha_2}+\varepsilon^{1/3}).
		\]
		Indeed, for $(k_1,k_2)\ne(0,0)$,
		\[
		|k_1+\sqrt2k_2|\,|k_1-\sqrt2k_2|
		=|k_1^2-2k_2^2|\geq1,
		\]
		and $|k_1-\sqrt2k_2|\leq\sqrt3\,|(k_1,k_2)|$.
		Thus \eqref{qpnew:diophantine} holds with
		$\gamma=1/\sqrt3$ and $\nu=1$.
	\end{remark}
	\section*{Acknowledgments}
	X. Su is supported by the National Natural Science Foundation of China (Grant No. 12371186, 12571207).

	\vspace{0.5cm}
	
\end{document}